\documentclass[11pt]{amsart}
\usepackage[a4paper,margin=1in]{geometry}
\usepackage{amsmath,amssymb,amsthm,mathtools}
\usepackage[T1]{fontenc}
\usepackage{lmodern,microtype,mathrsfs}
\usepackage{enumitem,booktabs,longtable,array}
\usepackage[colorlinks=true,linkcolor=blue!50!black,citecolor=blue!50!black,urlcolor=blue!50!black]{hyperref}
\usepackage{xcolor}
\hypersetup{pdftitle={Construction and exhaustion of enhanced supercuspidal Langlands parameters from Galois-side Yu data},pdfauthor={Manish Mishra}}
\numberwithin{equation}{section}
\theoremstyle{plain}
\newtheorem{theorem}{Theorem}[section]
\newtheorem{proposition}[theorem]{Proposition}
\newtheorem{lemma}[theorem]{Lemma}
\newtheorem{corollary}[theorem]{Corollary}
\theoremstyle{definition}
\newtheorem{definition}[theorem]{Definition}

\theoremstyle{remark}
\newtheorem{remark}[theorem]{Remark}
\DeclareMathOperator{\Irr}{Irr}
\DeclareMathOperator{\Hom}{Hom}
\DeclareMathOperator{\depth}{depth}
\DeclareMathOperator{\im}{im}
\DeclareMathOperator{\Ad}{Ad}

\DeclareMathOperator{\SL}{SL}
\DeclareMathOperator{\PGL}{PGL}

\newcommand{\C}{\mathbb C}
\newcommand{\Z}{\mathbb Z}

\newcommand{\Dpres}{\mathscr D_Z^{\mathrm{pres}}(G,\zeta_G)}

\newcommand{\Fr}{\mathrm{Fr}}
\newcommand{\hG}{\widehat G}

\newcommand{\cS}{\mathcal S}

\newcommand{\LG}[1]{{}^{L}#1}

\newtheorem{convention}[theorem]{Convention}
\newtheorem{maintheorem}{Theorem}

\DeclareMathOperator{\val}{val}

\newcommand{\OO}{\mathcal O_F}
\newcommand{\pp}{\mathfrak p_F}

\title[Enhanced supercuspidal parameters]{Construction and exhaustion of enhanced supercuspidal Langlands parameters from Galois-side Yu data}
\author{Manish Mishra}
\address{
Department of Mathematics \\
Indian Institute for Science Education and Research \\
Dr.\ Homi Bhabha Road, Pashan \\
Pune 411 008 \\
India
}
\email{manish@iiserpune.ac.in}
\thanks{This work was supported by the SERB Core Research Grant
(CRG/2022/000415) and the ANRF MATRICS Grant
(ANRF/ARGM/2025/000365/MTR)}
\date{}
\begin{document}
\begin{abstract}
Let \(G\) be a connected reductive group over a nonarchimedean local
field \(F\), and suppose that \(G\) splits over a tamely ramified
extension. Yu's construction of supercuspidal representations of \(G(F)\)
starts with a tower \(G^0\subsetneq\cdots\subsetneq G^d=G\) of tame
twisted Levi subgroups and depth-zero supercuspidal data for \(G^0(F)\).
At each step \(i<d\), it uses a positive-depth character of \(G^i(F)\)
that is generic relative to \(G^{i+1}\). Fintzen proved
that every irreducible supercuspidal representation arises in this way
when the residue characteristic \(p\) does not divide the order of
the absolute Weyl group of \(G\).

We prove a counterpart on the Galois side of the local Langlands
correspondence under the same hypothesis. We work with cuspidal
Langlands parameters equipped with a rigid enhancement, whose central
character records relevance to a prescribed rigid inner form.
Starting from a depth-zero enhanced parameter, we build such parameters
by adjoining positive-depth characters along a tower of twisted Levi
subgroups. Conversely, we recover data of this kind from every
relevant cuspidal rigidly enhanced parameter by successively removing
its positive-depth factors. We also determine when two choices of data
give the same enhanced parameter. Once their restrictions to wild
inertia are identified, the answer is determined by the conjugacy
classes of the resulting depth-zero parameters with their enhancements.
\end{abstract}
\maketitle
\setcounter{tocdepth}{1}
\tableofcontents
\clearpage

\section{Introduction}\label{sec:introduction}

Yu's construction \cite{Yu2001} starts from a depth-zero
supercuspidal representation and builds representations along a tower
of twisted Levi subgroups. For a connected reductive group \(G\)
that splits over a tamely ramified extension of \(F\), Fintzen
\cite{FintzenTypes} proved that every irreducible supercuspidal
representation of \(G(F)\) arises in this way provided that
\begin{equation}\label{eq:standing-prime}
                       p\nmid |W_{\mathrm{abs}}(G)|.
\end{equation}
Here \(p\) is the residue characteristic of \(F\). We work under the
same hypotheses and give a construction and exhaustion of enhanced
supercuspidal Langlands parameters.

Our starting point is a depth-zero cuspidal rigidly enhanced parameter
for \(G^0\) in an elliptic tame twisted Levi tower
\[
                  G^0\subsetneq G^1\subsetneq\cdots\subsetneq G^d=G.
\]
Consider a single step from \(G^i\) to \(G^{i+1}\). A positive-depth
character supplies a cocycle \(\kappa_i\) in the center of
\(\widehat{G^i}\). We use it to twist the parameter for \(G^i\) and
pass to \(G^{i+1}\) through a tame \(L\)-embedding \(\iota_i\).
What happens to the enhancement? Genericity ensures that the
centralizer in \(\widehat{G^{i+1}}\) of the leading wild image of
\(\kappa_i\) is exactly \(\widehat{G^i}\). Thus every element
centralizing the new parameter already lies in \(\widehat{G^i}\).
There the twist is central, so the centralizers of the two parameters
agree. Using the same finite central subgroup \(Z\) at both levels
then identifies the refined component groups on which the enhancements
are defined. Writing \(\alpha_i\) for this identification in a fixed
rational reference, the recursion is
\[
 \phi_{i+1}=\kappa_i\,\iota_i(\phi_i),\qquad
 \rho_{i+1}^+=\rho_i^+\circ\alpha_i^{-1},
 \quad
 \alpha_i:\mathcal S_{\phi_i}^{+,Z}
                   \xrightarrow{\sim}\mathcal S_{\phi_{i+1}}^{+,Z}.
\]

The term \emph{supercuspidal parameter} means a discrete parameter
whose rigid enhancement satisfies the generalized-Springer cuspidality
condition of Dillery--Schwein \cite{DillerySchwein}, as specified in
Definition~\ref{sc:rigid-cuspidality}. Under the same tameness and
residue-characteristic hypotheses, Kaletha constructs supercuspidal
\(L\)-packets for discrete Weil parameters \(W_F\to{}^LG\).
On the parameter side, he classifies torally wild parameters by
torally wild \(L\)-packet data \cite[Proposition~4.1.8]{KalethaPackets}.
Here we start with an intrinsic cuspidal rigidly enhanced depth-zero
parameter, which may have a nontrivial \(\SL_2\)-factor. We impose no
regularity assumption on this seed. The construction transports its
individual enhancement along the tower and determines when two such
presentations yield the same enhanced parameter.

\subsection{The construction and its input}

The character at level \(i\) is supplied on the full torus of
abelian quotients,
\[
 D_i=G^i/G^i_{\mathrm{der}},\qquad
 \bar\chi_i:D_i(F)\longrightarrow\C^\times.
\]
The torus LLC therefore determines a full cohomology class in
\(H^1(W_F,Z(\widehat{G^i})^\circ)\), not merely a leading
inertia map. The same finite central subgroup \(Z\) is used
throughout the tower. The seed includes its rigidification and its
marked relative realization in the prescribed ambient rigid form.
These rational data distinguish relevance sectors that an abstract
source group alone does not distinguish.

The use of a common cover matters when comparing enhancements. The
intrinsic component groups of Aubert--Moussaoui--Solleveld (AMS)
use the adjoint quotient and simply connected cover of each group
separately. Appendix~\ref{sc:rank-one} computes an elliptic torus
example in \(\SL_2\) where those component groups have different
orders, even though the ordinary parameter centralizers agree.

For a wild homomorphism \(\omega:P_F\to\hG\) extending to a Weil
parameter, write
\(L_\omega=Z_{\hG}(\omega(P_F))\) and let \(H_\omega\) be the
quasi-split group with this dual Levi and the induced outer Weil action.
We use the elliptic wild classes, the marked references
\(j_\omega:H_\omega\hookrightarrow G^*\), and the Weil sections
specified in Section~\ref{sc:wild-classification}.

\begin{maintheorem}\label{thm:main-overview}
Fix a tame ambient rigid inner form of class \(\zeta_G\), a finite central
subgroup \(Z\) representing its rigidification, and assume
\eqref{eq:standing-prime}.
\begin{enumerate}[label=\textup{(\roman*)},leftmargin=2.5em]
\item Every strict rationally marked presentation in the sense of
Definition~\ref{sc:classification-data}, with its chosen compatible
tame \(L\)-embeddings and a relevant depth-zero cuspidal rigid seed,
determines a relevant cuspidal rigidly enhanced parameter of the
ambient group. Full centralizer
cut-down, admissibility, and discreteness hold at every generic stage.
\item Every relevant cuspidal rigidly enhanced parameter occurs.
A reverse construction produces a finite tower of the specified
kind and reconstructs the prescribed enhanced parameter.
\item For each elliptic wild class \([\omega]\), the \(L\)-embedding defined
by its chosen Weil section and transport in the common cover give a bijection
\[
 \coprod_{\eta\in j_{\omega,*}^{-1}(\zeta_G)}
 \Phi_{0,\mathrm{cusp}}^{+,Z}(H_\omega,\eta)
 \xrightarrow{\ \sim\ }
 \Phi_{\mathrm{cusp}}^{+,Z}(G,\zeta_G)_{[\omega]} .
\]
The target consists of the classes with wild restriction conjugate to
\(\omega\). Two strict presentations give the same enhanced class
exactly when their wild restrictions are conjugate and their reduced
depth-zero enhanced seeds agree up to \(L_\omega\)-conjugacy.
\end{enumerate}
\end{maintheorem}

Corollary~\ref{sc:proved-construction-exhaustion} gives the surjection
from rationally marked presentations. Theorem~\ref{sc:wild-bijection}
proves the bijection in \textup{(iii)}. Theorem~\ref{sc:exact-fibers}
computes its two invariants directly from the tower factors and the
seed. Proposition~\ref{sc:rational-weyl-obstruction} explains why the
elementary rational moves in Definition~\ref{sc:exact-equivalence}
can generate a proper subrelation of equality of enhanced outputs.

\subsection{Reverse extraction and comparison of factorizations}

For a positive-depth parameter \(\phi\), let \(K\) be its
top-depth image and put \(L=Z_{\widehat G}(K)\). The prime
hypothesis makes this full centralizer a connected Levi. The isogeny
\(Z(L)^\circ\to L_{\mathrm{ab}}\) has degree prime to \(p\),
so its inverse on \(p\)-power torsion defines a central projection
of the finite wild image. That projection extends to a full Weil
cocycle. Dividing by it lowers depth and replaces the wild image
by a proper quotient. This proves termination without a separate
bound on ramification denominators. The parameter-side centralizer,
rational-reference, and tame-splitting results of Dat--Fintzen
\cite{DatFintzen} are used here.

The wild centralizer flag recovers the strict dual tower and its jump
depths. A conjugator of enhanced parameters may act on that tower by
a Weyl element which has no representative in the rational normalizer
of the marked embedding. For the fiber calculation we use the
reduction to depth zero at a fixed wild restriction for underlying
Weil parameters, as in Dat--Helm--Kurinczuk--Moss \cite[\S3]{DHKM}
and Dat--Fintzen \cite[\S4.3.1, diagram~(4.3.2)]{DatFintzen}.
We fix the homomorphism on \(P_F\) pointwise. Its stabilizer is then exactly
\(L_\omega\). Dividing out the chosen section carrying \(\omega\)
leaves a depth-zero parameter for \(H_\omega\), together with the
same common-cover label. The remaining conjugacy test takes place
entirely in \(L_\omega\).

\subsection{Relation to representation-side constructions}

The depth-zero input is itself a cuspidal enhanced parameter. Thus
the construction can be carried out entirely on the parameter side,
independently of finite Jordan decomposition or a depth-zero local
Langlands correspondence. A later representation-side application
must compare its depth-zero seed, character normalization, and
equivalence relation with those used here. A finite Heisenberg--Weil
representation alone does not provide an action of the refined
enhancement group and is not an input to this construction.

The input to the surjection includes the full quotient characters and
the marked rational realization. An application starting instead from
\((S,\theta,\phi_0,\rho_0^+)\), or retaining additional vertex labels,
requires a comparison with these data and with the fibers in
Theorem~\ref{sc:exact-fibers}.
Appendix~\ref{sec:quotient-normalization} records precisely what
positive-depth character normalization does provide.

\section{Conventions and quotient-character data}\label{sec:notation}
\label{sec:tame-yu-input}

At each proper step the character enters the construction in two ways.
Its restriction at the leading depth detects the smaller dual Levi
inside the larger group. The character on the full quotient torus
determines the Weil cocycle used to twist the parameter. The building
point specifies the filtration at which we test the genericity of
the leading restriction.

Let \(F\) be a nonarchimedean local field with ring of integers
\(\OO\), maximal ideal \(\pp\), residue characteristic \(p\),
and normalized valuation \(\val_F\). Write \(W_F\), \(I_F\),
and \(P_F\) for its Weil, inertia, and wild inertia groups, and
\(I_F^r\) for upper-numbered ramification groups. For a finite
extension \(E/F\), a valuation without further qualification
extends \(\val_F\); integer-normalized \(E\)-indices are
explicitly rescaled by the ramification index. Fix an additive
character \(\Psi\) trivial on \(\pp\) and nontrivial on
\(\OO\).

\begin{convention}[Parameters, depth, and tameness]
\label{conv:tame-MP}\label{conv:large-p}\label{conv:relevance}
For a connected reductive \(J\), use the \(\SL_2\)-model of
an admissible parameter
\(\phi:W_F\times\SL_2(\C)\to{}^LJ\): the Weil part is
continuous and Frobenius-semisimple, and the second factor is algebraic.
All groups in the positive-depth construction split tamely. Thus wild
inertia acts trivially on their pinned dual groups. The depth is
\[
 \depth(\phi)=\inf\{r\ge0:\phi(I_F^{r+})
                     \text{ has trivial dual-group part}\}.
\]
Depth zero means triviality on \(P_F\). We write
\(Z_{\widehat J}(\phi)\) for the centralizer of the full
parameter and \(Z_{\widehat J}(\phi|_{W_F})\) for the
Weil centralizer. Discreteness means that the former is finite
modulo \(Z(\widehat J)^{W_F}\). Relevance and cuspidality
refer to the fixed rigid inner form and are defined in Section
\ref{sc:rigid-comparison}. The standing prime assumption is
\eqref{eq:standing-prime}; individual algebraic or cohomological
lemmas stated more generally retain their stated scope.
\end{convention}

For a connected reductive group \(J\), put
\[
 D_J=J/J_{\mathrm{der}},\qquad q_J:J\longrightarrow D_J,
 \qquad A_J=Z(\widehat J)^\circ.
\]
The dual of \(q_J\) identifies \(\widehat D_J\) with
\(A_J\), including its Weil action. A character of \(J(F)\)
trivial on the simply connected derived image is not, for that
reason alone, a specified character of the full group \(D_J(F)\).

\begin{proposition}[The tame filtered torus correspondence]
\label{prop:torusLLC-filtered}
For a tame \(F\)-torus \(S\), the torus LLC
\[
 \mathcal L_S:\Hom_{\mathrm{sm}}(S(F),\C^\times)
                       \xrightarrow{\sim}H^1(W_F,\widehat S)
\]
is functorial and preserves depth. In particular the quotients of
classes of depth at most \(r\) by those of depth less than \(r\)
agree, and restriction of a full class of positive depth \(r\)
to \(I_F^r/I_F^{r+}\) gives its leading inertia map.
\end{proposition}
\begin{proof}
Functoriality is that of the torus correspondence. The filtration
statement is \cite[Theorem~7.10]{YuTori}; see also
\cite[Lemma~2.4.1]{DatFintzen} for its use for group characters.
Wild inertia acts trivially on \(\widehat S\), so coboundaries
do not change the positive-depth restriction.
\end{proof}

\subsection{Genericity and its geometric witnesses}\label{sec:genericity}

Write \(J(F)_{x,t}\) and \(\mathfrak j(F)_{x,t}\) for
Moy--Prasad groups and lattices at a compatible building point.
For positive \(t\), use the usual graded group--Lie algebra
identification. Normalize the dual lattices by
\[
 \mathfrak j^*(F)_{x,s}
  =\{X^*:X^*(\mathfrak j(F)_{x,(-s)+})\subset\pp\}.
\]
\begin{lemma}[Graded duality]\label{lem:MP-duality}
For \(t>0\), the pairing induced by \(\Psi\) identifies
\[
 \mathfrak j^*(F)_{x,-t}/\mathfrak j^*(F)_{x,(-t)+}
\]
with the character group of
\(\mathfrak j(F)_{x,t}/\mathfrak j(F)_{x,t+}\).
\end{lemma}
\begin{proof}
The annihilator of \(\mathfrak j(F)_{x,t+}\) is the first
dual lattice, and that of \(\mathfrak j(F)_{x,t}\) is the
second. The resulting finite additive quotients are in perfect
duality through \(\Psi\).
\end{proof}

\begin{definition}[Generic quotient character]
\label{def:generic-element}\label{def:generic-character}
Let \(H\subset G\) be a tame twisted Levi and let
\(\chi=\bar\chi\circ q_H\) have positive depth \(r\).
It is \(G\)-generic if its leading character is represented by
an \(H\)-invariant central functional \(X^*\) and, for a
tame maximal torus \(S\subset H\),
\[
 \val_F(X^*(H_\alpha))=-r
 \quad\bigl(\alpha\in R(S,G)\setminus R(S,H)\bigr),
 \qquad H_\alpha=d\alpha^\vee(1),
\]
together with the absolute Weyl stabilizer condition \(\mathrm{GE2}\)
on its leading reduction. In the standing prime range Lemma
\ref{sc:graded-genericity} proves that the displayed root
condition already implies \(\mathrm{GE2}\), and identifies
it with nontriviality of the leading dual root characters.
\end{definition}

Ellipticity of a rational twisted Levi \(H\subset G\) means
that \(Z(H)^\circ/Z(G)^\circ\) is anisotropic. A tower is
elliptic if its successive inclusions are elliptic. No condition
of regularity is imposed on an auxiliary pair \((S,\theta)\).

\subsection{The factorized input}
\begin{definition}[Quotient-character datum]\label{sc:def-data}\label{def:galois-side-yu-datum}
Fix an elliptic tame twisted Levi tower
\[
 G^0\subsetneq G^1\subsetneq\cdots\subsetneq G^d=G,
\]
the geometric data needed for genericity, a finite central
\(F\)-subgroup scheme \(Z\subset Z(G)\), and the compatible
dual covers of Section~\ref{sc:rigid-comparison}.
Use this same \(Z\), contained in every \(G^i\), at every stage.
Supply an intrinsic depth-zero rigidly enhanced seed
\[
 \Phi_0^+=(\phi_0,\rho_0^+;\mathfrak r_0),\qquad
 \rho_0^+\in\Irr(\cS_{\phi_0}^{+,Z}),
\]
cuspidal in the common-cover convention of
Definition~\ref{sc:rigid-cuspidality}. The compatible relative
rigidification record \(\mathfrak r\), whose seed part is
\(\mathfrak r_0\), is part of the datum as in
Definition~\ref{sc:rigid-seed}; it is not quotiented out by
forgetting the ambient embedding or rigidification.
A \emph{quotient-character datum} additionally supplies full smooth
characters
\[
 \bar\chi_i:D_i(F)\longrightarrow\C^\times,
 \qquad D_i=G^i/G^i_{\mathrm{der}},\qquad 0\le i\le d,
\]
and sets \(\chi_i=\bar\chi_i\circ q_{G^i}\).
For \(i<d\), require \(\chi_i\) to be \(G^{i+1}\)-generic
of depth \(r_i\), with \(0<r_0<\cdots<r_{d-1}\). The final
\(\bar\chi_d\) is an arbitrary central twist. The positive-depth
filtration comparison identifies the depth of \(\bar\chi_i\)
with that of \(\chi_i\), as in Lemma~\ref{sc:graded-genericity}.

The \(\bar\chi_i\), not merely their pullbacks or leading restrictions,
are part of the datum. Compatible \(L\)-embeddings and the full
centralizer conditions are explicit inputs in general; Theorem
\ref{sc:generic-stage-theorem} supplies them for elliptic tame
data with \(p\nmid|W_{\mathrm{abs}}(G)|\). Once they are supplied,
the refined enhancement operation is the transport proved in
Theorem~\ref{sc:rigid-preservation} and
Proposition~\ref{sc:rigid-relevance}.
\end{definition}

The induction starts at \(G^0\) with the intrinsic enhanced seed.
The first proper transition uses \(\bar\chi_0\); the final
\(\bar\chi_d\) is an ambient-central twist. The case \(d=0\)
has no proper transition. Tori and compatible building points are
genericity witnesses; they are not inferred from the depth-zero seed.
Definition~\ref{sc:classification-data} specifies the presentations,
Definition~\ref{sc:exact-equivalence} their elementary rational moves,
and Theorem~\ref{sc:exact-fibers} the global fiber criterion.

\subsection{Full central cocycles and representatives}
\label{sc:full-classes}\label{def:toral-cocycle}
\begin{proposition}[Full classes from full characters]
\label{sc:full-cocycle}
For a quotient-character datum, the usual LLC for tori gives a specified
full class
\[
 [\kappa_i]=\mathcal L_{D_i}(\bar\chi_i)
       \in H^1(W_F,A_i),\qquad A_i=Z(\widehat{G^i})^\circ.
\]
Assume the tame filtration comparison and choose any representative
\(\kappa_i\). Its leading restriction is the one obtained from
\(\chi_i\) by restriction to a compatible maximal torus and the
torus LLC. No lift from a prescribed graded class is required in this
forward construction.
\end{proposition}
\begin{proof}
The full class is the torus LLC of the supplied full character.
For a maximal torus \(S\subset G^i\), the restriction of \(\chi_i\)
is the pullback of \(\bar\chi_i\) along \(S\to D_i\).
Functoriality of the torus LLC identifies its class with the image of
\([\kappa_i]\) under \(A_i\to\widehat S\). Restricting this
identity to the depth quotient gives the leading map, using the stated
filtration compatibility.
\end{proof}

\begin{lemma}[Change of cocycle representative]\label{sc:coboundary}
Fix a compatible \(L\)-embedding \(\iota:{}^LH\to{}^LG\), a
parameter \(\phi_H\), and a central cocycle \(\kappa\) for the
Weil action on \(A_H\). For \(a\in A_H\), set
\[
 \kappa'(w)=a\kappa(w)w(a)^{-1}.
\]
Then
\[
 \kappa'(w)\iota(\phi_H(w))
     =a\bigl(\kappa(w)\iota(\phi_H(w))\bigr)a^{-1},
\]
with the same conjugation of the \(\SL_2\)-part. Thus the underlying
parameter class is independent of the representative of \([\kappa]\).
\end{lemma}
\begin{proof}
Conjugation by \(\iota(\phi_H(w))\) on \(A_H\) is the prescribed
Weil action, since the dual-group part in \(\widehat H\) centralizes
\(A_H\). Move \(a^{-1}\) across this factor. The element \(a\)
also centralizes the \(\SL_2\)-image.
\end{proof}

The cocycle is central in the source dual group, not necessarily in
the target. Its unramified part is retained; no condition
\(\kappa(\Fr)=1\) is imposed. Enhancement covariance under these
conjugacies follows from Theorem~\ref{sc:rigid-preservation}.

\section{Refined component groups and rigid seeds}
\label{sc:rigid-comparison}

Even for an elliptic torus in split \(\SL_2\), fixing the ambient
rigid inner form can leave more than one possible rigidification of
the torus. In the example of Appendix~\ref{sc:rank-one}, the split
ambient form admits two relevant enhancements of the parameter.
Each of the two torus rigidifications above the split class selects
one of them. To recover both enhancements, the seed must therefore
record how the torus is realized in the prescribed ambient form.

We use Kaletha's refined component groups with the same finite
central subgroup \(Z\) at every level. In the resulting common
cover, the centralizer comparison transports the enhancement, and
its central character records the rigid inner form to which it is
relevant. We also specify the generalized-Springer cuspidality
condition on the enhancement of the depth-zero seed.

\subsection{The common cover}

Fix a finite central \(F\)-subgroup scheme \(Z\subset Z(G)\),
and use its induced inclusion in every member of the twisted Levi
tower. It is chosen sufficiently large to realize the ambient
rigidification under consideration. Write
\[
 \overline{G^i}=G^i/Z,\qquad
 \pi_i^Z:\widehat{\overline{G^i}}\longrightarrow\widehat{G^i}.
\]
These are finite central isogenies of connected complex reductive
groups. The notation \(\pi_i^Z\) is distinct from the quotient
map \(q_{G^i}:G^i\to D_i\) used for full characters.
Compatible dual realizations mean that the chosen inclusions identify
\begin{equation}\label{sc:eq-cover-cartesian}
 \widehat{\overline{G^i}}
   =(\pi_{i+1}^Z)^{-1}(\widehat{G^i}),\qquad
 \pi_i^Z=\pi_{i+1}^Z|_{\widehat{\overline{G^i}}}.
\end{equation}
Thus every cover is a full inverse-image subgroup of the final cover.
For a dual Levi inclusion induced by compatible root data of
\(G^i\subset G^{i+1}\) and their quotients by \(Z\), this is
the usual dual central-isogeny construction: the inverse image contains
the inverse-image maximal torus and has the lifted Levi roots.
We require these compatible realizations; an abstract inclusion of
unrelated dual groups does not supply them. No lift of the Weil part
of an \(L\)-embedding to \(\widehat{\overline G}\) is required
for taking the inverse images below.

For a parameter of a group \(J\) in the tower, set
\begin{equation}\label{sc:eq-rigid-groups}
 S_{\phi_J}=Z_{\widehat J}(\phi_J),\qquad
 S_{\phi_J}^{+,Z}=(\pi_J^Z)^{-1}(S_{\phi_J}),\qquad
 \cS_{\phi_J}^{+,Z}=\pi_0(S_{\phi_J}^{+,Z}).
\end{equation}
This is the refined group of \cite[\S5.4]{KalethaRigid}, with the
dependence on \(Z\) displayed, and the enhancement group of
\cite[Definition~4.2.1]{DillerySchwein} in the \(\SL_2\)-model.
In particular, we do not first
divide \(S_{\phi_J}\) by the whole center of \(\widehat J\).

\begin{definition}[Rigid enhancement]\label{def:enhanced}
A rigid enhancement of \(\phi_J\) at level \(Z\) is an
irreducible representation \(\rho_J^+\) of
\(\mathcal S_{\phi_J}^{+,Z}\). For a fixed rigid inner twist
we require the relevance condition \eqref{sc:eq-rigid-relevance}.
Cuspidality is Definition~\ref{sc:rigid-cuspidality}.
\end{definition}

\begin{theorem}[Refined component groups are preserved in one step]
\label{sc:rigid-preservation}\label{prop:component-groups-yu-step}
Let \(H\subset G\) be a step with compatible covers
\eqref{sc:eq-cover-cartesian}. Suppose
\(\phi_G=\kappa\iota(\phi_H)\) is an admissible parameter,
where \(\kappa\) is central in the source dual group, and suppose
the full ordinary-centralizer equality holds:
\[
 \iota\bigl(Z_{\widehat H}(\phi_H)\bigr)
       =Z_{\widehat G}(\phi_G).
\]
Then the lifted dual inclusion induces an isomorphism of algebraic
groups \(S_{\phi_H}^{+,Z}\simeq S_{\phi_G}^{+,Z}\), and hence
\begin{equation}\label{sc:eq-rigid-preservation}
 \boxed{\;\cS_{\phi_H}^{+,Z}
       \xrightarrow{\;\sim\;}\cS_{\phi_G}^{+,Z}.\;}
\end{equation}
These identifications compose along a tower using the same \(Z\).
They are compatible with dual conjugacy; the induced map on
irreducible representation classes is independent of the lift of
a conjugating element to the central cover.
\end{theorem}
\begin{proof}
Identify \(\widehat H\) with its image and put
\(C=Z_{\widehat H}(\phi_H)=Z_{\widehat G}(\phi_G)\).
The Cartesian cover condition gives, inside
\(\widehat{\overline G}\),
\[
 (\pi_H^Z)^{-1}(C)=(\pi_G^Z)^{-1}(C).
\]
These are exactly the two groups in \eqref{sc:eq-rigid-groups}.
Taking connected components proves
\eqref{sc:eq-rigid-preservation}. All identifications along a
tower are restrictions inside the final cover, so their composites
agree. A dual conjugating element has a lift through the surjective
central isogeny. Any two such lifts differ by a central kernel element
and induce the same conjugation. If two conjugacies between the same
parameter realizations differ by a parameter-centralizer element,
their induced maps differ by an inner automorphism of the refined
group. They therefore give the same map on irreducible classes.
\end{proof}

Lemma~\ref{sc:ordinary-centralizers} supplies the centralizer
equality under the cut-down hypotheses of this paper. The toral
instance of the inverse-image argument is precisely the mechanism in
\cite[Lemma~5.3.4]{Kaletha}. The theorem above uses the already
established full centralizer equality; it does not infer that equality
from rootwise genericity alone.

\subsection{Relative rigidification and relevance}

Define the central relevance group and its natural map by
\begin{align}
 Z(\widehat{\overline J})^{+,Z}
   &:=(\pi_J^Z)^{-1}\bigl(Z(\widehat J)^{\Gamma_F}\bigr),
   &\mathcal Z_J^{+,Z}&:=\pi_0\bigl(Z(\widehat{\overline J})^{+,Z}\bigr),
       \label{sc:eq-relevance-center}\\
 a_J&:\mathcal Z_J^{+,Z}\longrightarrow\cS_{\phi_J}^{+,Z}.
       \notag
\end{align}
The image of \(a_J\) is central; \(a_J\) need not be injective.
The centrality follows also directly: a lift of an element central
in \(\widehat J\) has commutators in the finite central kernel,
and the resulting continuous map from the connected covering group
to that finite kernel is constant.

Write \(H^1_{\mathrm{rig}}(F,Z\to J^*)\) for the rigid
cohomology of a quasi-split reference form. In characteristic zero
this denotes \(H^1(u\to W,Z\to J^*)\) as in
\cite{KalethaRigid}; for local function fields use the rigid-gerbe
theory of \cite{Dillery}. The algebraic comparison theorem itself
has no characteristic-zero restriction. A rigid inner twist class
\(\zeta_J\), realized by the chosen \(Z\), specifies a
character \(\chi_{\zeta_J}\) of \(\mathcal Z_J^{+,Z}\)
via rigid duality. A representation \(\rho_J^+\) is relevant
to this rigid twist when
\begin{equation}\label{sc:eq-rigid-relevance}
 \rho_J^+(a_J(t))=\chi_{\zeta_J}(t)\operatorname{Id}
       \qquad(t\in\mathcal Z_J^{+,Z}).
\end{equation}
This is the relevance convention of
\cite[\S5.4, diagram~(5.7)]{KalethaRigid} and
\cite[Definition~5.2.1]{DillerySchwein}.

\begin{definition}[Rigid seed and relative realization record]
\label{sc:rigid-seed}
Fix an ambient rigid inner twist, its class \(\zeta_G\), and
the common \(Z\). The seed of Definition~\ref{sc:def-data}
consists of \((\phi_0,\rho_0^+)\) and a record
\(\mathfrak r_0\) specifying its rigid inner twist class
\(\zeta_0\) and its relative realization in this ambient twist.
For the tower, a compatible record \(\mathfrak r\) specifies
the intermediate rigid twists and the embeddings/comparison
identifications realizing them in the ambient group. It retains the
rational conjugacy information of these embeddings.

Whenever the reference groups have compatible \(F\)-embeddings,
the record can be expressed by lifts \(\zeta_i\) of
\(\zeta_G\) through the induced rigid-cohomology maps, together
with identifications of the induced twists. If quasi-split reference
forms are chosen separately, such maps are not assumed to exist
merely from the notation: the relative embedding and comparison
data must be supplied. At the dual level the record must satisfy
the character compatibility in
Proposition~\ref{sc:rigid-relevance}.

Equivalent rational embedding data may replace a rigidification
record only after specifying their relative invariant and checking
that rigid duality gives the same central characters and transport.
An abstract isomorphism type of \(G^0\), or an unmarked stable
embedding class, is not this record. Existence of a compatible
record for every proposed input is not part of the definition's assertion.
Theorem~\ref{sc:one-step-realization} solves this realization problem
over a supplied compatible rational reference embedding.
\end{definition}

\begin{proposition}[Relevance under refined transport]
\label{sc:rigid-relevance}
In Theorem~\ref{sc:rigid-preservation}, suppose the actions on
the centers are compatible with the \(L\)-embedding. The lifted
inclusion then gives a map
\[
 t_{G,H}:\mathcal Z_G^{+,Z}\longrightarrow\mathcal Z_H^{+,Z}.
\]
Let \(\alpha:\cS_{\phi_H}^{+,Z}\xrightarrow{\sim}
\cS_{\phi_G}^{+,Z}\) be the refined comparison. Suppose the
relative rigidification record satisfies
\begin{equation}\label{sc:eq-character-compatibility}
 \chi_{\zeta_H}\circ t_{G,H}=\chi_{\zeta_G}.
\end{equation}
Then
\[
 \rho_G^+=\rho_H^+\circ\alpha^{-1}
\]
is relevant to \(\zeta_G\) whenever \(\rho_H^+\) is
relevant to \(\zeta_H\). For rigidifications related by an
actual reference embedding, \eqref{sc:eq-character-compatibility}
is the functoriality of rigid duality.
\end{proposition}
\begin{proof}
An ambient central fixed element belongs to the source center and
is fixed by its Weil action, since the extra conjugating factors
of the \(L\)-embedding act trivially on ambient central elements.
Taking inverse images yields the asserted map. The diagram of
natural component maps satisfies
\[
 a_G=\alpha\circ a_H\circ t_{G,H}.
\]
Consequently
\[
 \rho_G^+(a_G(t))
  =\rho_H^+(a_H(t_{G,H}(t)))
  =\chi_{\zeta_H}(t_{G,H}(t))\operatorname{Id}
  =\chi_{\zeta_G}(t)\operatorname{Id}.
\]
This proves relevance without assuming that either central
component map is injective.
\end{proof}

\begin{remark}[One fixed source sector is not the whole target sector]
The theorem identifies the full refined groups and hence all their
irreducible classes. It need not identify the relevance sectors of
one fixed pair of rigid twists bijectively. The \(H\)-central
character can impose more conditions than the \(G\)-central
character. Exhaustion on a fixed ambient rigid form must allow the
appropriate source rigidifications over it, or equivalent relative
embedding data. Conversely, an enhanced target determines a
character on the source central group after refined transport.
Theorem~\ref{sc:enhancement-realization} realizes that character by
a rigid inner twist and a marked rational embedding when a compatible
rational reference embedding has been fixed. Producing the reference
embedding and the other seed data is supplied in the prime-to-Weyl
range by Theorem~\ref{sc:reverse-step-theorem}.
\end{remark}

\subsection{Cuspidal local systems in the same cover}

\begin{definition}[Cuspidal rigid enhancements in the common cover]
\label{sc:rigid-cuspidality}\label{def:cuspidal}
For an admissible parameter \(\phi_J\), put
\[
 Q_{J,\phi_J}^{+,Z}
   =(\pi_J^Z)^{-1}\bigl(Z_{\widehat J}(\phi_J|_{W_F})\bigr).
\]
Let \(\gamma_J:\SL_2(\C)\to\widehat J\) be its algebraic
factor. It lifts uniquely to \(\gamma_J^Z\) in
\(\widehat{\overline J}\), since \(\SL_2\) is simply
connected; set
\[
 u_J^+=\gamma_J^Z
   \left(\begin{smallmatrix}1&1\\0&1\end{smallmatrix}\right).
\]
The canonical identification of the component groups of the
centralizers of \(\gamma_J^Z\) and \(u_J^+\) realizes
\(\rho_J^+\) as an equivariant local system on the
\(Q_{J,\phi_J}^{+,Z}\)-orbit of \(u_J^+\).
We call \((\phi_J,\rho_J^+)\) cuspidal if the full
\(W_F\times\SL_2(\C)\)-parameter \(\phi_J\) is discrete in the sense
of Convention~\ref{conv:tame-MP} and the local system satisfies the
cuspidality condition of \cite[Definitions~2.6.1 and~4.4.1]{DillerySchwein}.
Thus its constituents for the identity component of \(Q_{J,\phi_J}^{+,Z}\)
are cuspidal local systems.

In the Weil--unipotent model, put \(u_J=\pi_J^Z(u_J^+)\) and
\(C=Z_{\widehat J}(\phi_J|_{W_F},u_J)\).
The full parameter centralizer is a Levi factor of \(C\)
by \cite[Fact~4.2.2]{DillerySchwein}. Thus discreteness is equivalent
to finiteness of the maximal reductive quotient \(C/R_u(C)\) modulo
the image of \(Z(\widehat J)^{W_F}\). The unipotent radical \(R_u(C)\)
can have positive dimension.
\end{definition}

\begin{proposition}[Preservation of the local-system data]
\label{sc:rigid-cuspidality-preserved}\label{prop:cuspidal-output}
In addition to Theorem~\ref{sc:rigid-preservation}, assume the
ordinary Weil-centralizer equality of
Lemma~\ref{sc:ordinary-centralizers}, and identify the
\(\SL_2\)-factors through the chosen embedding. Then the
common-cover groups \(Q_{H,\phi_H}^{+,Z}\) and
\(Q_{G,\phi_G}^{+,Z}\), their unipotent elements, and their
equivariant local systems under refined transport agree. In
particular, if the output is discrete, a cuspidal refined seed
gives a cuspidal refined output.
\end{proposition}
\begin{proof}
Apply the inverse-image argument to the common ordinary
Weil-centralizer. This identifies the two groups \(Q^{+,Z}\).
Uniqueness of the lift from \(\SL_2\) identifies their
homomorphisms \(\gamma^Z\) and hence their unipotent elements.
The inverse image of the full parameter centralizer equals
the centralizer of \(\gamma^Z\) in \(Q^{+,Z}\): an element
centralizing its projection has commutators with \(\SL_2\)
in the finite covering kernel, and these commutators are trivial
by connectedness. The centralizer of the \(\SL_2\)-homomorphism
is a Levi subgroup of the unipotent centralizer, whose complementary
unipotent radical is connected. Thus their component groups agree;
compare \cite[Fact~4.2.2 and Corollary~4.2.3]{DillerySchwein}.
Under these identifications the representation \(\rho^+\) and
the equivariant local system are the same on both sides. The
cuspidality condition on the common identity component is therefore
unchanged. Discreteness is the separately stated condition needed
to obtain a cuspidal parameter.
\end{proof}

This proves transport of cuspidality as defined in
Definition~\ref{sc:rigid-cuspidality}.
For the AMS interface, \cite[Remark~5.2.3]{DillerySchwein} provides
a comparison in prescribed central-character sectors with the stated
rigidifying choices. Any application of that comparison must identify its choices with the common
\(Z\) and relative records used here, rather than assert an
isomorphism of intrinsic AMS groups across the tower. Any local-system
descent still requires the appropriate kernel to act trivially.

\section{Generic stages and the forward construction}\label{sec:onestep}\label{sc:centralizers}

We prove the full centralizer cut-down, the fixed-center criterion for
discreteness, and the generic forward step. The positive-depth
root calculation is separated from transport in the common cover.

Let \(\widehat H\subset\widehat G\) be a dual Levi containing
\(\widehat S\), and let \(K\subset Z(\widehat H)^\circ\) be
a finite subgroup. Write \(W(\widehat G,\widehat S)_{(K)}\) for
the subgroup fixing every element of \(K\).

\begin{proposition}[Full cut-down under explicit root and stabilizer conditions]
\label{sc:cutdown}
Suppose
\[
 \widehat\alpha|_K=1\ \Longleftrightarrow\
       \widehat\alpha\in\Phi(\widehat H,\widehat S),\qquad
 W(\widehat G,\widehat S)_{(K)}=W(\widehat H,\widehat S).
\]
Then \(Z_{\widehat G}(K)=\widehat H\).
\end{proposition}
\begin{proof}
The first identity computes the connected centralizer root by root as
\(\widehat H\). Every component of the full centralizer can be
represented by an element normalizing \(\widehat S\): conjugate
its transported maximal torus back by an element of \(\widehat H\).
Its Weyl class fixes \(K\) pointwise, so the second identity places
that representative in \(\widehat H\).
\end{proof}

\begin{lemma}[Equality of ordinary parameter centralizers]
\label{sc:ordinary-centralizers}
Let \(\phi_G=\kappa\iota(\phi_H)\) be an admissible step with
\(\depth(\phi_H)<r\), a tame embedding, and leading image
\(K\) satisfying Proposition~\ref{sc:cutdown}. Then
\[
 Z_{\widehat G}(\phi_G|_{W_F})
       =Z_{\widehat H}(\phi_H|_{W_F}),\qquad
 Z_{\widehat G}(\phi_G)=Z_{\widehat H}(\phi_H).
\]
\end{lemma}
\begin{proof}
A centralizer element on the left centralizes the top-depth image, hence
belongs to \(\widehat H\). Within \(\widehat H\), the factors
\(\kappa(w)\) are central and cancel from the centralizer equation.
Intersecting with the centralizer of the common \(\SL_2\)-image
gives the second equality.
\end{proof}

These are ordinary-centralizer statements, not intrinsic AMS group
comparisons. In unrestricted characteristic the root and pointwise
stabilizer conditions must remain separate. The standing prime hypothesis implies the required full stabilizer condition below.

\subsection{The prime-to-Weyl-order range}
\label{sc:generic-stages}

We use the standing hypothesis, recorded again for reference:
\begin{equation}\label{sc:eq-prime-range}
                 p\nmid |W_{\mathrm{abs}}(G)|.
\end{equation}
In particular the range \(p>|W_{\mathrm{abs}}(G)|\)
satisfies every prime hypothesis below. The same condition holds for
all geometric Levi subgroups and for the dual root datum.

\begin{lemma}[Finite \(p\)-images and their full centralizers]
\label{sc:p-centralizers}
Let \(\widehat J\) be a connected complex reductive group with
\(p\nmid|W(\widehat J)|\). Every finite \(p\)-subgroup \(B\)
of \(\widehat J\) is toral, and its full centralizer is a connected
Levi subgroup. Consequently, if a toral \(p\)-subgroup \(K\)
has connected centralizer \(\widehat H\), then
\[
 Z_{\widehat J}(K)=\widehat H,\qquad
 W(\widehat J,\widehat S)_{(K)}=W(\widehat H,\widehat S).
\]
\end{lemma}
\begin{proof}
The full-centralizer assertion is
\cite[Lemma~2.2.9(iv)]{DatFintzen}, using Steinberg's
prime-to-torsion centralizer results. Since \(B\) centralizes its
own centralizer and the centralizer of a Levi is its center, it is
contained in the center of that Levi and hence in a maximal torus.
The asserted Weyl identity follows by taking the normalizer of
\(\widehat S\) in the equality of full centralizers.
\end{proof}

\begin{lemma}[Positive filtrations and the generic leading map]
\label{sc:graded-genericity}
Let \(H\subset G\) be a tame twisted Levi inclusion satisfying
\eqref{sc:eq-prime-range}. Let \(\bar\chi\) be a full
character of \(D_H(F)\) of positive depth \(r\), let
\(\chi=\bar\chi\circ q_H\), and let
\(\kappa\in Z^1(W_F,Z(\widehat H)^\circ)\) represent its
torus-LLC class. Fix a tame maximal torus \(S\subset H\), a
finite tame Galois splitting field \(E/F\) for \(S\), and a
compatible building point \(x\). Put
\(K=\kappa(I_F^r)\), so \(\kappa(I_F^{r+})=1\).
Then \(\chi\) has depth \(r\) at \(x\). For a root
\(\alpha\in R(S,G)\setminus R(S,H)\), the following are equivalent:
\begin{enumerate}[label=\textup{(\roman*)},leftmargin=2.5em]
\item \(\widehat\alpha|_K\ne1\);
\item \(\chi\circ N_{E/F}\circ\alpha^\vee\) is nontrivial
on \(E_r^\times/E_{r+}^\times\), where the valuation on \(E\)
extends that on \(F\);
\item a central leading functional \(X^*\) for \(\chi\)
satisfies \(\operatorname{val}_F(X^*(H_\alpha))=-r\).
\end{enumerate}
If these conditions hold for every such root, the absolute Weyl
stabilizer condition \(\mathrm{GE2}\) holds as well. Thus the
character is generic in the full Yu sense in this prime range.
\end{lemma}
\begin{proof}
For \(t>0\), applying \cite[Lemma~3.3.2(3)]{Kaletha} to the tame group
\(H\) and its connected derived subgroup \(H_{\mathrm{der}}\)
gives the exact sequence
\[
 1\longrightarrow H_{\mathrm{der}}(F)_{x,t}
 \longrightarrow H(F)_{x,t}\xrightarrow{q_H}D_H(F)_t
 \longrightarrow 1.
\]
In particular \(q_H(H(F)_{x,t})=D_H(F)_t\). Taking the union
over depths strictly larger than \(t\) gives the same equality at
\(t+\). Hence \(\chi\) is nontrivial on \(H(F)_{x,r}\) and trivial
on \(H(F)_{x,r+}\). By compatibility of the graded Moy--Prasad
maps with \(q_H\), a leading functional for \(\bar\chi\) pulls back
along \(dq_H\) to an \(H\)-invariant central functional for \(\chi\).
Depth preservation for the tame torus LLC is
\cite[Theorem~7.10]{YuTori}; its version for group characters is
also recorded in \cite[Lemma~2.4.1]{DatFintzen}.

Functoriality and base change for the torus LLC identify the
parameter of the character in (ii) with
\(\widehat\alpha\circ\kappa|_{W_E}\). If
\(e=e(E/F)\), tame upper numbering gives
\(I_E^{er}=I_F^r\) and \(I_E^{er+}=I_F^{r+}\).
Local reciprocity therefore proves (i)\(\Leftrightarrow\)(ii).
Here \(er\), rather than \(r\), is the depth measured with
the integer-normalized valuation of \(E\).

For (ii)\(\Leftrightarrow\)(iii), choose an additive character
of \(F\) trivial on \(\mathfrak p_F\) and nontrivial on
\(\mathcal O_F\). Its composition with \(\operatorname{Tr}_{E/F}\)
has conductor \(\mathfrak p_E\), because the different of a
tame extension has exponent \(e-1\). Linearizing the normed
coroot on the leading quotient gives the pairing
\(u\mapsto\psi_F(\operatorname{Tr}_{E/F}(X^*(H_\alpha)u))\).
It is nonzero on the quotient at integer depth \(er\) precisely
when \(\operatorname{val}_E(X^*(H_\alpha))=-er\).
This is also the comparison in \cite[Lemma~3.6.8]{Kaletha}.
For roots of \(H\) the scalar is zero, since the functional
annihilates the derived Lie algebra.

For completeness, the stabilizer assertion follows by averaging.
Let \(\bar X\) be the reduction of a depth-zero rescaling of
\(X^*\) on the split torus, and let \(W_{\bar X}\) be its
stabilizer in the absolute Weyl group \(W\). It contains
\(W_H\). Lift \(\bar X\) to the integral character-lattice
representation over a characteristic-zero discrete valuation ring
with the same residue field, and average the lift over
\(W_{\bar X}\). Division by \(|W_{\bar X}|\) is allowed by
\eqref{sc:eq-prime-range}; the averaged lift still reduces
to \(\bar X\). It is fixed by \(W_H\), so it annihilates
the coroots of \(H\), and it annihilates no other coroot, by
reduction. The characteristic-zero Weyl stabilizer theorem, proved
by choosing a real chamber, makes its stabilizer exactly \(W_H\).
Thus \(W_{\bar X}\subset W_H\), proving \(\mathrm{GE2}\).
One may take a finite residue field containing all coefficients of
\(\bar X\), so the characteristic-zero lift used here exists
also when \(F\) has positive characteristic.
\end{proof}

\begin{proposition}[Genericity gives full cut-down]
\label{sc:generic-full-cutdown}\label{prop:centralizer-cutdown-generic-cocycle}\label{lem:genericity-implies-cocycle-regularity}
In Lemma~\ref{sc:graded-genericity}, if \(\chi\) is
\(G\)-generic of depth \(r\), then
\[
 Z_{\widehat G}\bigl(\kappa(I_F^r)\bigr)=\widehat H.
\]
In particular the absolute pointwise Weyl stabilizer condition in
Proposition~\ref{sc:cutdown} is a conclusion, not a further
hypothesis, under \eqref{sc:eq-prime-range}.
\end{proposition}
\begin{proof}
The roots of \(\widehat H\) vanish on the leading image,
which lies in its connected center. Every other root is nontrivial
there by Lemma~\ref{sc:graded-genericity}. The connected
centralizer is therefore \(\widehat H\). The leading image is
a finite \(p\)-group, so Lemma~\ref{sc:p-centralizers}
identifies the full centralizer with its identity component.
\end{proof}

\begin{lemma}[Ellipticity and fixed connected dual centers]
\label{sc:elliptic-fixed-centers}\label{lem:discrete-regular}
For a rational twisted Levi inclusion \(H\subset G\) with its
full compatible dual root-data marking, one has
\[
 (Z(\widehat G)^{W_F})^\circ
       \subset (Z(\widehat H)^{W_F})^\circ.
\]
Equality holds if and only if \(Z(H)^\circ/Z(G)^\circ\)
is anisotropic. Consequently a step with equal ordinary parameter
centralizers preserves discreteness in both directions whenever
the rational twisted Levi inclusion is elliptic.
\end{lemma}
\begin{proof}
The inclusion follows from the compatible \(L\)-embedding;
conjugation by its Weil part acts on \(Z(\widehat G)\) as
the ambient Galois action. For any connected reductive \(J\),
\[
 \dim (Z(\widehat J)^{W_F})^\circ
     =\operatorname{rank}_F Z(J)^\circ.
\]
Indeed, over \(\mathbb Q\) the character module of
\(Z(\widehat J)^\circ\) is the quotient of the absolute
cocharacter module by the coroot span. It is isogenous to the
cocharacter module of \(Z(J)^\circ\). A finite Galois group
has equally dimensional invariants and coinvariants on a rational
representation. This proves the displayed dimension formula.
The two included tori have equal dimension exactly when the
quotient of connected centers has split rank zero. Finally a
parameter is discrete exactly when the identity component of its
centralizer is the fixed connected dual center. Apply this
criterion to the equal ordinary centralizers.
\end{proof}

\begin{theorem}[The generic forward stage in the prime-to-Weyl range]
\label{sc:generic-stage-theorem}\label{thm:onestep-target}\label{prop:L-embedding-kappa}
Let \(H\subset G\) be an elliptic tame twisted Levi inclusion
satisfying \eqref{sc:eq-prime-range}, with the full compatible
dual marking. Let \(\bar\chi\) be a full quotient character
whose pullback to \(H(F)\) is generic of positive depth \(r\),
and let \(\phi_H\) be an admissible discrete parameter of depth
strictly smaller than \(r\). A compatible tame \(L\)-embedding
\(\iota:{}^LH\to{}^LG\) exists. For any such embedding and
any representative of the quotient-torus class \(\kappa\),
\[
                  \phi_G=\kappa\,\iota(\phi_H)
\]
is admissible, discrete, and of depth \(r\). Its full ordinary
and Weil centralizers equal those of \(\phi_H\).
The common-\(Z\) enhancement transport of Section
\ref{sc:rigid-comparison} therefore applies, including
cuspidality and the prescribed relevance condition.
\end{theorem}
\begin{proof}
Use simultaneous quasi-split references from
Theorem~\ref{sc:simultaneous-references}. For odd \(p\), choose
minimally ramified \(\chi\)-data as in \cite[\S5.5]{KalethaCovers}.
For a symmetric root \(\alpha\), the restriction of \(\chi_\alpha\)
to \(F_{\pm\alpha}^{\times}\) is the quadratic character attached to
\(F_\alpha/F_{\pm\alpha}\). With these tame data,
\cite[\S6.1, Corollary~6.7]{KalethaCovers} gives a compatible
\(L\)-embedding whose Weil part has trivial dual-group component on
wild inertia. If \(p=2\)
under \eqref{sc:eq-prime-range}, the root system is empty
and there is no proper twisted Levi step to construct.

The central cocycle identity proves multiplicativity, and the
\(\SL_2\)-image commutes with the central factors. We record
why Frobenius semisimplicity is preserved. Put
\(y=\iota(\phi_H(\Fr))\) and \(a=\kappa(\Fr)\).
The embedding preserves semisimplicity: after a power killing the
finite pinned actions, its image of \((1,\Fr)\) has power in
\(Z(\widehat H)\), since the centralizer of a dual Levi is
its center. Thus the image of any semisimple source Frobenius has
a semisimple power. Next take a power \(N\) for which
\(\Ad(y)^N\) is trivial on \(Z(\widehat H)^\circ\).
Then
\[
 (ay)^N=\Bigl(\prod_{j=0}^{N-1}\Ad(y)^j(a)\Bigr)y^N
\]
is a product of commuting semisimple elements. Therefore \(ay\)
is semisimple. Inertia has finite image by continuity and compactness:
a profinite subgroup of a complex Lie group is finite. The parameter
is admissible. The same argument works for inverse central twisting.

On \(I_F^r\) the smaller-depth parameter and the tame embedding
contribute no dual-group factor, so the restriction equals the
nontrivial leading map of \(\kappa\); on \(I_F^{r+}\) it is
trivial. Proposition~\ref{sc:generic-full-cutdown} and
Lemma~\ref{sc:ordinary-centralizers} give both full centralizer
equalities. Lemma~\ref{sc:elliptic-fixed-centers} proves
discreteness. The common-cover assertions now follow from
Theorem~\ref{sc:rigid-preservation} and Propositions
\ref{sc:rigid-relevance} and
\ref{sc:rigid-cuspidality-preserved}.
\end{proof}

The generic stage theorem applies to elliptic
quotient-character data in \eqref{sc:eq-prime-range}.
Ellipticity and relevance of the supplied rigid seed remain genuine
conditions on the data. The theorem does not assert full cut-down
from rootwise genericity at excluded primes.

\subsection{Iteration along the tower}\label{sec:induction}

\begin{theorem}[Forward construction with full characters and rigid seeds]
\label{sc:forward-theorem}\label{thm:inductive-target}\label{hyp:stage-inputs}
Let an elliptic tame quotient-character datum as in Definition
\ref{sc:def-data} be supplied in the standing prime range, with
one common \(Z\), compatible covers, and an intrinsic depth-zero
cuspidal rigid seed. Retain compatible relative records whose central
characters satisfy \eqref{sc:eq-character-compatibility}.
Choose compatible tame \(L\)-embeddings and express enhancement
transport in a fixed rational reference presentation.
Then the full torus classes of Proposition~\ref{sc:full-cocycle}
give the recursive underlying parameters. The enhancement is defined by
\[
 \rho_{i+1}^+=\rho_i^+\circ\alpha_i^{-1},\qquad
 \alpha_i:\cS_{\phi_i}^{+,Z}\xrightarrow{\sim}
                    \cS_{\phi_{i+1}}^{+,Z}.
\]
The output is relevant and cuspidal in the sense of
Definition~\ref{sc:rigid-cuspidality}, with depth \(r_i\) after the
\(i\)-th positive-depth transition. It is invariant under the specified central
refactorizations and intrinsic seed twists with fixed compatible
realizations and retained rigidification data. It is also invariant
under the compensated tame embedding changes of
Theorem~\ref{sc:tame-embedding-compensation}, covariant under actual
reference isomorphisms, and compatible with the endpoint comparisons
of Theorem~\ref{sc:endpoint-factorization-comparison}.
Here the raw transport is written in a fixed reference presentation;
changes to geometrically matched marked references use the coherent
normalization of Theorem~\ref{sc:coherent-origin-transport}.
Theorem~\ref{sc:exact-fibers} characterizes equality of enhanced outputs
for arbitrary strict presentations.
\end{theorem}
\begin{proof}
The full classes are supplied by Proposition~\ref{sc:full-cocycle}.
The central cocycle identity gives multiplicativity after the compatible
embedding. Since the old parameter and the tame embedding have trivial
dual-group part on the new highest wild-inertia quotient, the leading
restriction is exactly that of \(\kappa_i\), which proves the depth
statement. Theorem~\ref{sc:generic-stage-theorem} gives compatible tame
embeddings, admissibility, full centralizer equality, and discreteness.
Theorem~\ref{sc:rigid-preservation} defines \(\alpha_i\);
Proposition~\ref{sc:rigid-relevance} and the prescribed relative records give relevance;
Proposition~\ref{sc:rigid-cuspidality-preserved} gives cuspidality.
Iterate along the strictly increasing depths, retaining the relative
record. The invariance assertion follows from
Proposition~\ref{sc:rigid-compensation} and
Corollary~\ref{sc:enhanced-refactorization}, together with the
embedding and endpoint comparison results just cited.
\end{proof}

\section{Relative realization and simultaneous rational references}
\label{sec:relative-realization}

The refined enhancement determines a source relevance character.
We realize that character by a rigid inner twist and a marked
rational embedding in the prescribed ambient form, prove composition
along a tower, and construct simultaneous quasi-split references.
The source inner form is recovered from the enhancement; it is
not prescribed independently of that relevance character.
\subsection{One-step realization over a rational reference embedding}
\label{sc:relative-realization}

In this subsection, let \(j:H\hookrightarrow G\) be an actual
\(F\)-embedding of connected reductive groups whose geometric image is
a Levi subgroup. These are reference groups, and \(Z\subset Z(G)\)
is identified with the same finite central subgroup of \(H\). The
cohomological statements do not require either reference group to be
quasi-split. For their application to the absolute relevance convention
of Section~\ref{sc:rigid-comparison}, we will explicitly require
quasi-split references. All fields here are nonarchimedean local fields.

Write
\[
 \mathcal H_J=H^1_{\mathrm{rig}}(F,Z\to J),\qquad
 \mathcal Z_J=\mathcal Z_J^{+,Z},\qquad
 \mathcal Z_J^\vee=\Hom(\mathcal Z_J,\C^\times).
\]
Rigid duality gives functorial bijections
\begin{equation}\label{sc:eq-perfect-duality}
 D_J:\mathcal H_J\xrightarrow{\sim}\mathcal Z_J^\vee,
 \qquad
 D_G(j_*\eta)=D_H(\eta)\circ t_j.
\end{equation}
Here \(t_j:\mathcal Z_G\to\mathcal Z_H\) is dual to the
map induced by \([Z\to H]\to[Z\to G]\). The inputs are
\cite[Corollary~5.4]{KalethaRigid} for \(p\)-adic fields and
\cite[Corollary~7.11]{Dillery} for local function fields. The correction
to the real-field version of the former statement in
\cite{KalethaErrata} does not affect its nonarchimedean use here.

For compatible root-data realizations, \(t_j\) is exactly the map on
relevance centers already used above. To check the identification, choose
a maximal \(F\)-torus \(S\subset H\), also maximal in \(G\).
The character groups of the full inverse images of the fixed dual centers
are
\[
 X_*(S/Z)\big/\bigl(I_\Gamma X_*(S)+Q_J^\vee\bigr),\qquad J=H,G,
\]
where \(I_\Gamma X_*(S)\) denotes the subgroup generated by
\(\sigma\lambda-\lambda\) for \(\sigma\in\Gamma_F\),
and \(Q_J^\vee\) is the coroot
lattice in \(X_*(S)\). The inclusion
\(Q_H^\vee\subset Q_G^\vee\) gives the quotient map between these
groups. Taking torsion subgroups gives the character-group map for
\(t_j\); this is also the map used in rigid duality. We henceforth
write \(t_{G,H}=t_j\) when the chosen dual embedding realizes these
same root data. This verification concerns the central map and does not
construct the Weil part of an \(L\)-embedding.

Fix a rigid cocycle \(z_G\) with class \(\zeta_G\), and let
\(G_{z_G}\) be its inner twist. A \emph{relative reduction} is a
rigid \(H\)-torsor \(P_H\) on the rigid gerbe together with an
isomorphism
\[
 \beta:P_H\times^H G\xrightarrow{\sim}P_G.
\]
The isomorphism is part of the record. Isomorphisms of reductions must
respect it. The group \(G_{z_G}(F)=\operatorname{Aut}(P_G)\) acts
by postcomposition. Denote by \(\operatorname{Rel}_j(\zeta_G)\)
the set of isomorphism classes of reductions modulo this ambient action.
This is a classification of marked reductions, not of unmarked subgroup
images modulo arbitrary normalizer actions.

In characteristic zero, use the rigid extension
\(u\to\mathcal W\to\Gamma_F\); \(\mathcal W\) is not
the Weil group \(W_F\). A reduction is represented by
\((z_H,g)\), with \(g\in G(\overline F)\), satisfying
\begin{equation}\label{sc:eq-relative-cocycle}
 j(z_H(w))=g^{-1}z_G(w)\sigma(g),\qquad
       w\in\mathcal W,\quad \sigma=\overline w.
\end{equation}
This equality includes restriction to the band \(u\), on which both
sides take values in the common \(Z\). Source changes of trivialization
and ambient rational conjugacy are respectively
\begin{align}
 (z_H,g)&\longmapsto
       (h^{-1}z_H\sigma(h),\,g j(h)),\quad h\in H(\overline F),
       \label{sc:eq-relative-equivalence}\\
 (z_H,g)&\longmapsto(z_H,a g),\quad
       a^{-1}z_G\sigma(a)=z_G.
       \notag
\end{align}
For positive characteristic, the torsor definition is primary. These
formulas mean the same identities after a common fpqc trivialization,
with pullbacks in place of \(\sigma\) and equality on the gerbe
band. In particular no identification of a non-smooth band with its
geometric points is used.

\begin{theorem}[One-step relative realization]
\label{sc:one-step-realization}
With the preceding reference embedding fixed, there is a bijection
\begin{equation}\label{sc:eq-relative-fiber}
 \operatorname{Rel}_j(\zeta_G)
   \xrightarrow{\sim}j_*^{-1}(\zeta_G)
   \xrightarrow{\sim}
 \{\chi\in\mathcal Z_H^\vee:
                    \chi\circ t_j=\chi_{\zeta_G}\}.
\end{equation}
Consequently a prescribed character \(\chi_H\) is realized by a
relative reduction if and only if
\(\chi_H\circ t_j=\chi_{\zeta_G}\). When it is realized,
its source rigidification class and its class of marked reductions
modulo ambient rational conjugacy are unique.

A reduction \((z_H,g)\) gives an \(F\)-embedding
\[
 j_g:H_{z_H}\hookrightarrow G_{z_G},\qquad
       j_g=\Ad(g)\circ j
\]
in the twisting coordinates. It remains a geometric Levi embedding.
If \(j\) is elliptic, this embedding is elliptic as well.

Without a prescribed source character, the only obstruction to a
relative reduction is
\begin{equation}\label{sc:eq-relative-obstruction}
 \chi_{\zeta_G}|_{\ker(t_j)}=1.
\end{equation}
When nonempty, the set in \eqref{sc:eq-relative-fiber} is a torsor
under \(\Hom(\mathcal Z_H/\im(t_j),\C^\times)\).
\end{theorem}
\begin{proof}
Forgetting \(\beta\) maps a reduction to a class in
\(j_*^{-1}(\zeta_G)\). Conversely, equality of these cohomology
classes is precisely the existence of an isomorphism
\(P_H\times^H G\simeq P_G\). For fixed \(P_H\), any two
such isomorphisms differ by an automorphism of \(P_G\). Replacing
\(P_H\) by an isomorphic torsor has exactly the effect stipulated
in the definition of a reduction. This proves the first bijection.
The second bijection follows from both bijectivity and functoriality
in \eqref{sc:eq-perfect-duality}. This also proves existence and
uniqueness for a prescribed character; character compatibility is
sufficient here because the actual reference morphism \(j_*\) exists.

For clarity, the descent calculation for the embedding is explicit.
The twisted action on \(J_{z_J}\) is
\(x\mapsto z_J(w)\sigma(x)z_J(w)^{-1}\). Equation
\eqref{sc:eq-relative-cocycle} implies
\[
 z_G(w)\sigma(gj(x)g^{-1})z_G(w)^{-1}
   =g j\bigl(z_H(w)\sigma(x)z_H(w)^{-1}\bigr)g^{-1}.
\]
The band acts trivially by conjugation, so this is descent over \(F\).
Equivalently, in every characteristic the closed immersion is obtained
by descent of \(j\) through \(P_H\) and \(\beta\).
The geometric Levi property is preserved by this descent. Inner twisting
does not change centers or their \(F\)-actions, and \(j_g\) induces
the original center inclusion in these coordinates. The anisotropy of
\(Z(H)^\circ/Z(G)^\circ\), when assumed, is therefore unchanged.

One can also see the first bijection on cocycles: for two transporters
\(g,g'\) with the same \(z_H\), the element \(a=g'g^{-1}\)
satisfies \(a^{-1}z_G\sigma(a)=z_G\), and is an element of
\(G_{z_G}(F)\) in twisting coordinates. The remaining changes are
\eqref{sc:eq-relative-equivalence}.

Finally, a character of \(\mathcal Z_G\) factors through
\(\im(t_j)\) precisely when it kills \(\ker(t_j)\). Every
character of a subgroup of a finite abelian group extends to that
group: this follows by adding one cyclic generator at a time and
choosing a root in \(\C^\times\) of its prescribed value.
Two extensions differ by a unique character of
\(\mathcal Z_H/\im(t_j)\). Apply this to
\(\im(t_j)\subset\mathcal Z_H\).
\end{proof}

\begin{corollary}[Elliptic reference embeddings]
\label{sc:elliptic-realization}
If \(H\) contains a maximal \(F\)-torus that is elliptic in \(G\),
then \(j_*:\mathcal H_H\to\mathcal H_G\) is surjective and
\(t_j\) is injective. In particular every ambient rigidification
has a relative reduction. This applies to an elliptic twisted-Levi
embedding \(j\).
\end{corollary}
\begin{proof}
For such a torus \(S\), the composite
\(\mathcal H_S\to\mathcal H_H\to\mathcal H_G\) is surjective
by \cite[Corollary~3.7(2)]{KalethaRigid} or
\cite[Proposition~5.7(2)]{Dillery}. Thus \(j_*\) is surjective,
and perfect duality gives injectivity of \(t_j\).
For the last assertion choose an elliptic maximal torus in \(H\).
Its maximal split subtorus is the maximal split subtorus of \(Z(H)\),
which is that of \(Z(G)\) when \(H\) is elliptic in \(G\).
It is therefore elliptic in \(G\).
\end{proof}

\subsection{Recovering the record from marked rational embeddings}

The marking in the preceding theorem has a concrete interpretation.
Fix a rigid inner twist \((\xi_G,z_G):G\to G'\). Consider an
\(F\)-embedding \(k:H'\hookrightarrow G'\), an inner twist
\(\xi_H:H\to H'\), and a geometric conjugacy
\begin{equation}\label{sc:eq-marked-embedding}
 \xi_G^{-1}\circ k\circ\xi_H=\Ad(g)\circ j.
\end{equation}
Here \(\xi_H\) is initially an ordinary inner twist; a source
rigidification has not been supplied. This condition specifies the
inner marking of the embedding. A bare geometric conjugacy class of
subgroups does not specify it.

\begin{proposition}[Extraction from an actual marked realization]
\label{sc:marked-extraction}
Data \eqref{sc:eq-marked-embedding} determine a source rigidification
realized by the same \(Z\), and a relative reduction. In cocycle
coordinates it is
\begin{equation}\label{sc:eq-extracted-rigid-cocycle}
 z_H(w)=j^{-1}\bigl(g^{-1}z_G(w)\sigma(g)\bigr).
\end{equation}
This rigidifies the given inner twist \(\xi_H\). Changing \(g\)
without changing \eqref{sc:eq-marked-embedding}, or changing the
inner marking by an inner isomorphism, gives an isomorphic source
rigidification. Ambient rational conjugacy gives the corresponding
equivalent relative record. In particular its relevance character is
the one furnished by \eqref{sc:eq-perfect-duality}.
\end{proposition}
\begin{proof}
We first justify that the argument of \(j^{-1}\) belongs to \(j(H)\).
Write \(d_G(\sigma)=\xi_G^{-1}\sigma(\xi_G)=\Ad(z_G(w))\)
and \(d_H(\sigma)=\xi_H^{-1}\sigma(\xi_H)\), the latter an
inner automorphism of \(H\). Apply \(\sigma\) to
\eqref{sc:eq-marked-embedding} and use that \(j,k\) are defined
over \(F\). With \(n(w)=g^{-1}z_G(w)\sigma(g)\), this gives
\[
 \Ad(n(w))\circ j=j\circ d_H(\sigma).
\]
A geometric Levi subgroup has scheme-theoretic centralizer equal to its center: choose
a maximal torus of \(H\), also maximal in \(G\), and then impose
commutation with its root subgroups. Since \(d_H(\sigma)\) is inner,
the displayed equality implies
\(n(w)\in j(H)Z_G(j(H))=j(H)\).
The function \(n\) is a coboundary conjugate of \(z_G\), so it
is a cocycle. On the band it equals \(z_G|_u\), since that map has
central image in \(Z\). Hence \eqref{sc:eq-extracted-rigid-cocycle}
is a rigid \(H\)-cocycle with band in the same \(Z\); its adjoint
action is \(d_H\), so it rigidifies \(\xi_H\).
These identities are descent identities and hold on an fpqc
trivialization of the gerbe as well. Locally an inner automorphism
lifts to an element of \(H\); the scheme-theoretic centralizer identity
then places \(n\) in \(j(H)\), and this membership descends.
This proves the function-field case as well.

If \(g'\) has the same property with the same \(\xi_H,k\), then
\(g'=g j(h)\) for \(h\in Z(H)(\overline F)\).
The extracted cocycles differ by \(z_H\mapsto h^{-1}z_H\sigma(h)\).
If \(\xi_H\) is replaced by \(\xi_H\Ad(h)\), the same calculation
holds with \(h\in H(\overline F)\). Postcomposition of \(k\)
by ambient rational conjugacy is the second change in
\eqref{sc:eq-relative-equivalence}. These are precisely the allowed
isomorphisms of records. Formula \eqref{sc:eq-perfect-duality}
then determines the central character and its compatibility.
\end{proof}

Thus the relative invariant is the extracted class \([z_H]\),
together with a realization representative when an actual embedding
is needed. The theorem does not quotient by source outer automorphisms
or by the geometric normalizer of \(H\). Such operations can change
the marked parameter, its character, and its enhancement; their treatment
belongs to a later comparison theorem.

\subsection{Realization of the enhancement and fixed-tower composition}

\begin{theorem}[Recovery of the source relevance sector]
\label{sc:enhancement-realization}
Assume now that \(j:H\hookrightarrow G\) is a rational embedding
between quasi-split reference groups, with the compatible dual root data
above. Fix a parameter step satisfying Theorem~\ref{sc:rigid-preservation}
and Proposition~\ref{sc:rigid-relevance}, so that
\[
 \alpha:\cS_{\phi_H}^{+,Z}\xrightarrow{\sim}\cS_{\phi_G}^{+,Z},
 \qquad a_G=\alpha\circ a_H\circ t_j.
\]
For every \(\zeta_G\)-relevant irreducible \(\rho_G^+\), the
representation \(\rho_H^+=\rho_G^+\circ\alpha\) determines a
unique \(\zeta_H\in j_*^{-1}(\zeta_G)\) to which it is relevant.
This class has a marked rational realization in \(G_{z_G}\), unique
up to the equivalences of Theorem~\ref{sc:one-step-realization}.
In particular there is a bijection
\begin{equation}\label{sc:eq-sector-bijection}
 \coprod_{\zeta_H\in j_*^{-1}(\zeta_G)}
 \Irr\bigl(\cS_{\phi_H}^{+,Z},\chi_{\zeta_H}\bigr)
 \xrightarrow{\sim}
 \Irr\bigl(\cS_{\phi_G}^{+,Z},\chi_{\zeta_G}\bigr),
 \qquad \rho_H^+\longmapsto\rho_H^+\circ\alpha^{-1}.
\end{equation}
Here the central-character condition is imposed through \(a_J\);
these maps need not be injective. If both parameters are discrete and
the Weil-centralizer hypotheses of
Proposition~\ref{sc:rigid-cuspidality-preserved} hold, the bijection
restricts to cuspidal enhancements.
\end{theorem}
\begin{proof}
The image of \(a_H\) is central. By Schur's lemma there is a unique
character \(\chi_H\) of \(\mathcal Z_H\) such that
\[
 \rho_G^+\bigl(\alpha(a_H(t))\bigr)
       =\chi_H(t)\operatorname{Id}.
\]
Using the central-map identity and relevance of \(\rho_G^+\) gives
\(\chi_H\circ t_j=\chi_{\zeta_G}\). Theorem
\ref{sc:one-step-realization} now gives the unique source class
and its relative realization. Notice that \(\chi_H\) is automatically
trivial on \(\ker(a_H)\); no injectivity is used.
Conversely each term of the coproduct maps into the target relevance
sector by Proposition~\ref{sc:rigid-relevance}. An irreducible
source representation has only one central character, and perfect duality
identifies that character with only one \(\zeta_H\). Thus the terms
are disjoint and the preceding construction is the inverse. The final
assertion follows from the equality of the local-system data in
Proposition~\ref{sc:rigid-cuspidality-preserved}.
\end{proof}

In the unramified torus example of Appendix~\ref{sc:rank-one}, \(t_j:\mu_2\hookrightarrow\mu_4\).
The characters extending the trivial ambient character are exactly
\(\lambda_0\) and \(\lambda_2\). Theorem
\ref{sc:enhancement-realization} therefore recovers the two source
rigidifications separately, with their marked rational embeddings.
It does not collapse them to the single intrinsic AMS torus enhancement.

\begin{proposition}[Composition, covariance, and enlargement of the cover]
\label{sc:relative-composition}
For a fixed chain of rational reference embeddings using one common
\(Z\), compatible relative reductions compose associatively. Their
central characters and the enhancement transports above compose with
them. The construction is invariant under the source and ambient
isomorphisms of records. It is also compatible with enlarging the
common \(Z\) and inflating the enhancement, with the reference
embeddings and dual realizations fixed.
\end{proposition}
\begin{proof}
It suffices to treat \(H\xrightarrow{j}M\xrightarrow{k}G\).
If the successive records are \((z_H,g_H)\) over \(z_M\) and
\((z_M,g_M)\) over \(z_G\), then their composite transporter is
\begin{equation}\label{sc:eq-relative-composition}
 g_{H,G}=g_M\,k(g_H).
\end{equation}
Substitution in \eqref{sc:eq-relative-cocycle} verifies this
formula. In the torsor formulation it is simply composition of the
isomorphisms after extension of structure group, so associativity
holds in every characteristic. Changing the intermediate representative
to \(z_M^h=h^{-1}z_M\sigma(h)\) changes the two transporters to
\(h^{-1}g_H\) and \(g_M k(h)\); these factors cancel in
\eqref{sc:eq-relative-composition}. Source and ambient changes
give exactly \eqref{sc:eq-relative-equivalence} at the endpoints.

On centers, \(t_{G,H}=t_{M,H}\circ t_{G,M}\), and on rigid
cohomology the maps compose in the other direction. Thus compatibility
of characters composes. The maps \(\alpha\) are restrictions in
the same dual cover and compose by Theorem~\ref{sc:rigid-preservation}.
Changing dual representatives through the specified conjugacies has
no effect on irreducible classes by that theorem. This is covariance
for the given comparison data, not a construction of comparisons
between arbitrary \(L\)-embeddings.

For \(Z\subset Z'\subset Z(G)\), regard the same cocycles as
\(Z'\)-rigid cocycles; their transporters and all the preceding
identities are unchanged. The map of dual covers induces maps
\(\cS_{\phi_J}^{+,Z'}\to\cS_{\phi_J}^{+,Z}\), and the
new enhancement is pullback along this map. Functoriality of rigid
duality identifies its central character with that of the enlarged
rigidification. This proves compatibility for the inflated objects;
it does not assert descent for arbitrary enhancements at level \(Z'\).
\end{proof}

\subsection{Simultaneous quasi-split references for a rational Yu tower}
\label{sc:general-reference-input}

The relevant general input is \cite[Lemma~6.4]{KalethaCovers}: if
\(\xi:K\to J\) is an inner twist with \(K\) quasi-split and
\(M\subset J\) is a rational twisted Levi subgroup, an inner
modification of \(\xi\) makes its inverse image of \(M\) rational
and quasi-split, with inner restriction to \(M\).
We apply this successively inside the already constructed terms.
There is no assertion that the original twisted Levi subgroup of a
quasi-split group is itself quasi-split.

\begin{theorem}[A simultaneous reference tower with one rigid cocycle]
\label{sc:simultaneous-references}
Let \(F\) be a nonarchimedean local field, and let
\[
 S\subset J_0\subsetneq J_1\subsetneq\cdots\subsetneq J_d=J
\]
be an actual rational twisted Levi tower, with \(S\) a maximal
\(F\)-torus of \(J_0\). Fix a quasi-split inner form \(K\) of
\(J\), a rigid inner twist \((\xi,z):K\to J\), and a finite
central subgroup scheme \(Z\subset Z(K)\) through which its band
map factors. Identify \(Z\) with its image under \(\xi\).
There exists \(g\in K(F^s)\) with the following properties. Put
\[
 \Xi=\xi\circ\Ad(g),\qquad
 K_i=\Xi^{-1}(J_i),\qquad S^*=\Xi^{-1}(S).
\]
\begin{enumerate}[label=\textup{(\roman*)},leftmargin=2.5em]
\item The \(K_i\) form a rational twisted Levi tower in \(K_d=K\),
every \(K_i\) is quasi-split, and
\(\Xi_i=\Xi|_{K_i}:K_i\to J_i\) is an inner twist.
The restriction \(\Xi_S:S^*\to S\) is an \(F\)-isomorphism.
\item The conjugated rigid cocycle
\begin{equation}\label{sc:eq-toral-reference-cocycle}
 z^g(w)=g^{-1}z(w)\sigma(g)
\end{equation}
takes values in \(S^*\), and its band still takes values in the
same \(Z\). Viewed in any \(K_i\), it rigidifies \(\Xi_i\).
Its pushforwards along all reference inclusions agree literally.
\item With the original ambient marking \((\xi,z)\) retained,
use \((\Xi_i,z^g)\) for \(i<d\). The adjacent relative
transporters can be chosen to be \(1\) below the final step and
\(g\) at the final step. Thus the construction presents the
original rational tower in the prescribed ambient rigid inner form.
\end{enumerate}
No enlargement of \(Z\), residue-characteristic bound, or
ellipticity hypothesis is required for this reference theorem.
\end{theorem}
\begin{proof}
First ignore the rigid cocycle. Apply the cited lemma to
\(J_{d-1}\subset J\) and \(\xi\). It gives a modification of
\(\xi\) for which \(K_{d-1}\) is rational and quasi-split and
the restriction to \(J_{d-1}\) is inner. Now apply the same lemma
inside \(K_{d-1}\), to that restricted inner twist and
\(J_{d-2}\subset J_{d-1}\). This further modification is
conjugation by an element of \(K_{d-1}(F^s)\). It preserves
\(K_{d-1}\) and all higher reference subgroups as subgroups;
their restrictions remain inner twists. Continuing downwards gives
one product \(g\) and a nested rational quasi-split reference tower.
Apply the lemma one last time to \(S\subset J_0\), inside
\(K_0\). A maximal torus is a geometric Levi subgroup, and an
inner twist of a torus is an \(F\)-isomorphism, since a torus
has no nontrivial inner automorphisms. This last modification again
preserves every reference subgroup. It proves (i). For \(d=0\)
only this last application is needed.

The cocycle in \eqref{sc:eq-toral-reference-cocycle} is the
usual change of representative of \(z\), and its adjoint action
is \(\Xi^{-1}\sigma(\Xi)\). The rationality of \(\Xi_S\)
implies that this action fixes \(S^*\) pointwise. The
scheme-theoretic centralizer of a maximal torus of a connected
reductive group is that torus. Hence \(z^g\) takes values in
\(S^*\). Its band equals the original band: conjugation by
\(g\) fixes \(Z\), and the band has trivial Galois projection.
Every \(K_i\) contains \(S^*\) and \(Z\), so the same
cocycle is a rigid cocycle for every term, with adjoint action
\(\Xi_i^{-1}\sigma(\Xi_i)\). This proves (ii). In positive
characteristic these are identities of descent data on a flat
trivialization of the rigid gerbe. The scheme-theoretic centralizer
identity proves membership in \(S^*\) there, and membership descends;
no smoothness of \(Z\) is used.

Write \(j_i:K_i\hookrightarrow K_{i+1}\) and
\(k_i:J_i\hookrightarrow J_{i+1}\). For \(i<d-1\),
\(k_i\Xi_i=\Xi_{i+1}j_i\); at the final step,
\[
 k_{d-1}\Xi_{d-1}=\xi\Ad(g)j_{d-1},\qquad
 j_{d-1}(z^g(w))=g^{-1}z(w)\sigma(g).
\]
These are exactly the marked embedding and cocycle identities of
Proposition~\ref{sc:marked-extraction}. Their composition has
transporter \(g\). When \(d=0\), the pair \((\Xi,z^g)\)
is simply an isomorphic representative of the given ambient rigid
twist. This proves (iii).
\end{proof}

\begin{proposition}[Full root data, covers, and quotient characters]
\label{sc:reference-root-data}
The reference presentation of
Theorem~\ref{sc:simultaneous-references} identifies the full
absolute root-data diagram of the actual tower, computed on \(S\),
with that of the reference tower, computed on \(S^*\), equivariantly
for \(\Gamma_F\). It identifies their common-\(Z\) quotient diagrams
as well. In these coordinates one can transport the connected dual
diagram and the compatible \(L\)-embeddings of the given actual tower.

The induced maps
\[
 D_{K_i}=K_i/(K_i)_{\mathrm{der}}
       \xrightarrow[\ \Xi_{i,\mathrm{ab}}\ ]{\ \sim\ }
 D_{J_i}=J_i/(J_i)_{\mathrm{der}}
\]
are \(F\)-isomorphisms. They transport the supplied full quotient
characters and their full torus-LLC classes. If \(S\) is tame,
then \(S^*\) is tame and all reference groups split over a tame
extension. Ellipticity of a tower inclusion is preserved.
\end{proposition}
\begin{proof}
The single geometric isomorphism \(\Xi\) carries every root
subsystem \(R(K_i,S^*)\) and its coroots to
\(R(J_i,S)\) and its coroots. Its restriction to the common
torus is rational. The character and cocharacter lattice maps,
including all root and coroot subsets, are therefore
\(\Gamma_F\)-equivariant and commute with the tower inclusions.
Inner modifications act trivially on the ambient center; the
identification with \(Z\) is the fixed one throughout. Passing to
the appropriate sublattices for quotient by \(Z\), and then
dualizing, gives the complete diagram of central covers and Levi
inclusions, not just its maps on dual centers.

Choose pinnings of the quasi-split reference groups. Their based
Galois actions are those obtained from the preceding root data by
the usual Weyl corrections. These are exactly the outer actions
for the actual inner forms. The connected dual diagram and its
common covers can thus be marked by these identifications, and
compatible Weil parts of already supplied \(L\)-embeddings are
carried along by the resulting \(L\)-isomorphisms. Simultaneously
Galois-stable compatible Borels for the entire tower are not needed.

An inner automorphism acts trivially on an abelian quotient, so
\(\Xi_{i,\mathrm{ab}}\) is rational. It commutes with the maps
from the common maximal torus and with the quotient maps between
successive terms. Pull back \(\bar\chi_i\) along this map;
functoriality of the torus LLC identifies its full class with the
original class in \(H^1(W_F,Z(\widehat{J_i})^\circ)\).
Torus filtrations and the toral leading restrictions are preserved
under these \(F\)-isomorphisms. This does not assert an
identification of buildings or transfer an independently chosen
depth-zero representation.

The splitting field of \(S\) also splits \(S^*\). A reductive
group containing a split maximal torus is split, so it splits
every \(K_i\). Finally the isomorphisms on centers commute with
their inclusions. They preserve the condition that the connected
center of a source modulo the ambient connected center is
anisotropic, which is the ellipticity condition here.
\end{proof}

\begin{corollary}[Reference existence for the forward Yu input]
\label{sc:general-reference-closure}
For quotient-character data of Definition~\ref{sc:def-data}
based on an actual rational tame twisted Levi tower in a specified
ambient rigid inner form, compatible quasi-split references with
full root-data and common-cover compatibility exist. A supplied change
of reference for that tower transports its relative records and edge
maps by Theorem~\ref{sc:coherent-origin-transport}.
\end{corollary}
\begin{proof}
Apply Theorem~\ref{sc:simultaneous-references} to the actual
tower and any chosen maximal \(F\)-torus of its bottom term,
using the tame torus when it is part of the datum. Apply
Proposition~\ref{sc:reference-root-data} to obtain the full
compatible dual presentation. Express the supplied rigid seed and
relative record in this presentation, retaining their inner markings.
Proposition~\ref{sc:marked-extraction} gives the cocycle formula
for marked realizations; Theorem~\ref{sc:enhancement-realization}
tests their prescribed relevance sectors. Apply Theorem
\ref{sc:coherent-origin-transport} to the marked change of reference
to obtain the asserted edge maps.
\end{proof}

This corollary supplies simultaneous references for
the \emph{forward} rational Yu data used in this paper. Its existence
statement does not select a distinguished rational origin or identify
different marked reductions: the rigidly enhanced seed and its
relative record remain part of the datum. In particular the two
source classes in Appendix~\ref{sc:rank-one} remain distinct.
An arbitrarily prescribed enhanced seed may fail the relevance
condition for a prescribed marked realization; reference existence
does not remove that condition.

For the reverse problem one may start only with a complex dual flag
and a final enhanced parameter. Producing a rational twisted Levi
tower of the required type, the full quotient characters, and the
depth-lowering factorization requires more than this reference theorem.
These are constructed in Section~\ref{sc:reverse} in the
prime-to-Weyl range. Once such a rational tower is available, its
simultaneous reference presentation and the
common-\(Z\) rigid restrictions are supplied by the theorem above.
Nor does the theorem identify an independently prescribed dual
diagram with the actual tower from equality of center maps alone:
the full dual root-data marking is essential.

\begin{remark}[Scope of a rational realization]\label{sc:reference-problem}
If a particular source inner form is prescribed in advance, its class
must equal the one recovered from the enhancement. Reference existence
does not remove this relevance condition. An actual reference embedding
supplies the map on rigid cohomology; an abstract dual inclusion alone
does not. Changes of rational origin are handled in Section
\ref{sec:reference-transport}, and extraction of the underlying
parameter tower is proved in Section~\ref{sc:reverse}.
\end{remark}
\section{Central seed twists and refactorization}
\label{sc:refactorization}

Different factorizations can give the same enhanced parameter. For
example, take a depth-zero character of the first quotient torus.
Twist the seed by its Weil cocycle and multiply the first character
factor by the inverse character. These two changes cancel in the
forward construction. We first describe central twists of the seed
and then show how the twists arising from full quotient characters
cancel along a tower.

\subsection{A seed operation requiring no depth-zero LLC}
\label{sc:seed-twist}

For a central cocycle \(z\in Z^1(W_F,Z(\widehat H))\), write
\[
 T_z(\phi,\rho^+;\mathfrak r)
       =(z\phi,\rho^+;\mathfrak r),
\]
whenever the twist is admissible, using the same \(Z\) and
relative rigidification record. The analogous notation without
superscripts applies when discussing the intrinsic AMS convention.

\begin{lemma}[Intrinsic central-twist covariance]\label{sc:central-twist}
Central twisting preserves the refined group, the common-cover
cuspidality condition, and the rigid central-character relevance
condition. It also preserves the corresponding AMS objects when
that convention is used. If both
\(\phi\) and \(z\) are trivial on wild inertia, then the twist is
again depth zero.
\end{lemma}
\begin{proof}
For every \(g\in\widehat H\), the central factors cancel in the
equation expressing that \(g\) centralizes \(z\phi\). Hence its
ordinary centralizer and its inverse image in the fixed
\(\widehat{H/Z}\) are unchanged, both for the Weil restriction
and for the full parameter. The same holds for the adjoint image
and its simply connected inverse image in the AMS convention.
The lifted unipotent element and relative rigidification are
unchanged. Containment in a proper
\(L\)-Levi is unchanged because every such Levi contains the dual
center. These facts preserve discreteness, the cuspidal local system,
and the relevance character. The last assertion follows by restriction
to wild inertia.
\end{proof}

This seed operation is defined directly on intrinsic parameters.
When comparing two data below, the relation \(\Phi'_0=T_{z_0}\Phi_0\)
is part of the intrinsic refactorization move. Proving that it mirrors
a twist of a representation belongs to a later LLC comparison.

\subsection{Central cohomology and abelian quotients}

An inclusion \(G^i\subset G^{i+1}\) induces a morphism
\[
 f_i:D_i\longrightarrow D_{i+1}.
\]
Its dual goes in the opposite direction:
\[
 \widehat f_i:A_{i+1}\longrightarrow A_i.
\]
Require compatibility of the chosen \(L\)-embedding with this map:
under \(\widehat{G^i}\hookrightarrow\widehat{G^{i+1}}\), the
image of \(\widehat f_i(a)\) is \(a\), with the correct Weil
actions. This is an explicit compatibility, not an identification of
every source-central cocycle with a target-central cocycle.

Take two tuples of full quotient-torus characters on the same tower.
For \(j\ge i\), let \(f_{i,j}:D_i\to D_j\) be the composite,
and set
\begin{align}
 \bar\xi_i&=\prod_{j=i}^{d}f_{i,j}^{*}
              (\bar\chi_j(\bar\chi_j')^{-1}),
      & [z_i]&=\mathcal L_{D_i}(\bar\xi_i),\label{sc:eq-xi}\\
 [\kappa_i]&=\mathcal L_{D_i}(\bar\chi_i),
      & [\kappa_i']&=\mathcal L_{D_i}(\bar\chi_i').\notag
\end{align}
The formula for \(\bar\xi_i\) uses full characters of the quotient
tori, so it does not involve choosing an extension from their rational
point images.

\begin{proposition}[Telescoping identity for quotient characters]
\label{sc:telescoping}
For \(0\le i<d\), one has
\begin{equation}\label{sc:eq-telescope}
 \boxed{\quad [\kappa_i'][z_i]
       =[\kappa_i]\,\widehat f_{i,*}[z_{i+1}]
       \quad\text{in }H^1(W_F,A_i).\quad}
\end{equation}
Assume compatible \(L\)-embeddings, admissible forward recursions, and
\(\phi'_0\sim z_0\phi_0\). Then
\[
 \phi'_i\sim z_i\phi_i\qquad(0\le i\le d).
\]
After the final factors \(\kappa_d,\kappa_d'\) are applied, the
two underlying parameters are conjugate.
\end{proposition}
\begin{proof}
The character identity
\[
 \bar\chi_i'\bar\xi_i
      =\bar\chi_i f_i^*\bar\xi_{i+1}
\]
and functoriality of the torus LLC prove
\eqref{sc:eq-telescope}. Choose representatives at the current
stage; its cohomological equality is an equality of cocycles up to a
coboundary in \(A_i\), which can be removed by
Lemma~\ref{sc:coboundary}. If \(\phi'_i\sim z_i\phi_i\), then
\begin{align*}
 \phi'_{i+1}
 &\sim\kappa_i'\,\iota_i(z_i\phi_i)
  =\kappa_i'z_i\,\iota_i(\phi_i)\\
 &\sim\kappa_i\widehat f_i(z_{i+1})\,\iota_i(\phi_i)
  =z_{i+1}\phi_{i+1}.
\end{align*}
All multiplications before the last equality take place in the image of
the source group; the last uses the compatibility of \(\widehat f_i\)
with \(\iota_i\). This proves the induction. Finally
\([z_d]=[\kappa_d][\kappa_d']^{-1}\), so the final twists cancel
the remaining discrepancy. The \(\SL_2\)-parts are transported by
the same embeddings and conjugacies throughout.
\end{proof}

All central cohomology products in \eqref{sc:eq-telescope} lie
in the source center. A cocycle central in a smaller group is not
thereby central in each larger group.

The enhancement compatibility corresponding to the telescoping
identity is the square
\begin{equation}\label{sc:eq-enhancement-square}
 \mathcal E_i^{\kappa_i'}\circ T_{z_i}
       =T_{z_{i+1}}\circ\mathcal E_i^{\kappa_i}.
\end{equation}
The bottom twist has depth zero so that it preserves the seed class.
The following proposition proves this square for common-cover transport.

\begin{proposition}[Compensation for the specified rigid transport]
\label{sc:rigid-compensation}
In Proposition~\ref{sc:telescoping}, use the same common
\(Z\), compatible dual inclusions, and relative rigidification
record for both factorizations. Assume all ordinary centralizer
and discreteness hypotheses needed for the two constructions.
Define enhancement transport by
\(\rho_{i+1}^+=\rho_i^+\circ\alpha_i^{-1}\) from
Theorem~\ref{sc:rigid-preservation}. Then the compensation
identity \eqref{sc:eq-enhancement-square} holds on refined
irreducible classes, with rigidification records retained.
\end{proposition}
\begin{proof}
A central twist at a fixed stage leaves its ordinary centralizers
and their inverse images in the fixed cover unchanged. It therefore
acts as the identity on refined enhancement labels and on the
rigidification record. The telescoping identity identifies the
two outputs up to the central twists and coboundary conjugacies
of Lemma~\ref{sc:coboundary}. These conjugacies lift to the
cover, and all enhancement identifications are restrictions of
these same lifted maps. Theorem~\ref{sc:rigid-preservation}
shows that choices of lift do not matter and that any ambiguity
from a parameter-centralizer conjugacy is inner. Hence the two
paths induce the same map on irreducible classes. Relevance and
cuspidality follow from Propositions~\ref{sc:rigid-relevance}
and~\ref{sc:rigid-cuspidality-preserved}.
\end{proof}

\begin{corollary}[Enhanced refactorization]
\label{sc:enhanced-refactorization}
Under Proposition~\ref{sc:telescoping} and
Proposition~\ref{sc:rigid-compensation}, with
\(\Phi'_0=T_{z_0}\Phi_0\), the final enhanced parameters are conjugate.
\end{corollary}
\begin{proof}
Apply the compensation identity at each stage and cancel the final
central twist as in the proposition.
\end{proof}

For changes of tame \(L\)-embedding, the difference cocycle may
lie in the full dual center rather than its identity component.
Section~\ref{sc:embedding-comparison} treats that compensation.
Section~\ref{sc:refactorization-bounds} proves the precise depth bounds
for refactorizations in a fixed presentation.

\section{Changes of reference and coherent enhancement transport}
\label{sec:reference-transport}

The recursion in the introduction is written in a fixed rational
reference. There are two ways its presentation can change. If we
change the Weil part of a tame \(L\)-embedding, a central twist of the
parameter compensates for the change. If we change the rational
reference embedding, we change the origin used to label the possible
source rigidifications. The label of the enhancement can then change
as well. We compute the character that corrects this label and show
that the corrections are compatible along the tower.

\subsection{Tame \texorpdfstring{\(L\)}{L}-embedding changes and compensation}
\label{sc:embedding-comparison}

\begin{lemma}[Difference of two \(L\)-embeddings]
\label{sc:L-embedding-difference}
Let \(\iota,\iota':{}^LH\to{}^LG\) have the same restriction
\(\widehat j:\widehat H\hookrightarrow\widehat G\), whose image
is a Levi subgroup. There is a unique continuous cocycle
\(\delta\in Z^1(W_F,Z(\widehat H))\) such that
\begin{equation}\label{sc:eq-L-embedding-difference}
 \iota'(h,w)=\widehat j(\delta(w))\,\iota(h,w).
\end{equation}
If the two embeddings agree on wild inertia, then \(\delta|_{P_F}=1\).
In particular this holds for embeddings with trivial dual-group part
on wild inertia in the conventions of this paper.
\end{lemma}
\begin{proof}
The ratio \(\iota'(1,w)\iota(1,w)^{-1}\) belongs to
\(\widehat G\). Both elements induce the prescribed action of
\(w\) on \(\widehat j(\widehat H)\), so their ratio centralizes
that Levi subgroup. Its centralizer is \(\widehat j(Z(\widehat H))\).
This defines \(\delta(w)\) uniquely. Multiplicativity of the two
embeddings gives
\(\delta(ww')=\delta(w)\,{}^w\delta(w')\), and continuity is
inherited. The assertion about wild inertia follows from the ratio.
\end{proof}

The center here is the \emph{full} center, not just
\(A_H=Z(\widehat H)^\circ\). Thus it is not legitimate to claim
that every embedding change comes from a character of \(D_H(F)\).
The intrinsic central twisting operation of Lemma~\ref{sc:central-twist}
uses the full center and is sufficient for the following comparison.

\begin{theorem}[Compensation for tame embedding changes on a fixed tower]
\label{sc:tame-embedding-compensation}
Fix a rational reference tower, its compatible dual inclusions, the
common \(Z\), all full quotient-character cocycles \(\kappa_i\),
and the relative rigidification records. For \(0\le i<d\), replace
\(\iota_i\) by \(\iota_i'\) with the same connected-dual
restriction and with tame difference cocycle
\(\delta_i\in Z^1(W_F,C_i)\), where
\(C_i=Z(\widehat{G^i})\). Keep the final central factor fixed.

The inclusions of full centers give equivariant maps
\(\nu_i:C_{i+1}\to C_i\). Define tame central cocycles backwards by
\begin{equation}\label{sc:eq-embedding-backward-twist}
 v_d=1,\qquad v_i=\delta_i^{-1}\nu_i(v_{i+1}).
\end{equation}
Replace the depth-zero seed by
\[
 (\phi_0',\rho_0^{+\prime};\mathfrak r_0')
       =(v_0\phi_0,\rho_0^+;\mathfrak r_0).
\]
Then the recursive parameters satisfy \(\phi_i'=v_i\phi_i\).
In particular the final enhanced parameters coincide, and all rigid
classes and relative realization records are unchanged. If the original
steps satisfy the full centralizer, admissibility, and discreteness
conditions, so do the compensated steps. The construction depends only
on the indicated cocycle classes up to the prescribed dual conjugacy.
Successive embedding changes compose with their compensating seed twists.
\end{theorem}
\begin{proof}
The maps \(\nu_i\) are equivariant because conjugation by the Weil
part of \(\iota_i\) has the specified action on \(C_i\), and its
restriction to the ambient center is the ambient Weil action.
Consequently \eqref{sc:eq-embedding-backward-twist} defines
cocycles; they are tame because all the \(\delta_i\) are tame.
Thus \(v_0\phi_0\) is again a depth-zero seed. In the source
center one has the exact identity
\[
 \kappa_i\delta_i v_i=\kappa_i\nu_i(v_{i+1}).
\]
Using \eqref{sc:eq-L-embedding-difference} in the recursion proves
\(\phi_{i+1}'=v_{i+1}\phi_{i+1}\) by induction. At the top
\(v_d=1\), and the final factors agree.

At every stage, central twisting leaves both the full and Weil
centralizers, their inverse images in the common cover, and the
\(\SL_2\)-factor unchanged. The connected-dual inclusions also
remain unchanged, so the maps \(\alpha_i\) on refined groups
are the same. This proves equality of the transported enhancement
labels. Their central characters and rigid classes are unchanged,
so the same relative records apply. Discreteness is preserved as in
Lemma~\ref{sc:central-twist}. Admissibility is preserved as well:
the inertia image remains finite, and for each Weil element a power
killing its action on the dual center expresses the twisted power as
a product of commuting semisimple elements, namely a central norm
and the corresponding power of the original parameter. Thus the
twisted Weil element is semisimple. The algebraic \(\SL_2\)-map
and its commutation relation are unchanged. Tameness of the twists
preserves all positive-depth restrictions.
The original leading-image cut-down therefore still applies.

Changing a cocycle representative introduces a central coboundary,
which is removed by the corresponding source dual conjugacy. The
proof of Lemma~\ref{sc:coboundary} works in the full center;
the common-cover conjugacy is independent of its lift. Finally the
difference cocycle for a composite change is the product of the
successive difference cocycles. Backward recursion is multiplicative,
so the seed twists and the unchanged enhancement identifications compose.
\end{proof}

This is a comparison of \emph{compensated data}. Holding the seed
fixed while changing an \(L\)-embedding can change the output.
The familiar change-of-\(\chi\)-data formula for tori is an instance
of this compensation; compare the proof of
\cite[Proposition~5.2.7]{Kaletha}. The argument above does not
assume regularity or invoke a representation-side LLC.

Embedding changes and the quotient-character refactorizations can also
be combined. In the cohomology of the full source center the required
equation is
\begin{equation}\label{sc:eq-general-embedding-compensation}
 [\kappa_i'][\delta_i][v_i]
       =[\kappa_i]\,\nu_{i,*}[v_{i+1}].
\end{equation}
At the top take \([v_d]=[\kappa_d][\kappa_d']^{-1}\).
If the resulting bottom twist is depth zero and both factorizations
satisfy the stage hypotheses, the same proof gives equality of final
enhanced outputs. Unlike the case with fixed characters in the theorem,
depth zero of this bottom twist is an additional condition.

\subsection{Covariance and endpoint comparison for different towers}
\label{sc:tower-endpoint-comparison}

\begin{proposition}[Covariance under actual reference isomorphisms]
\label{sc:reference-isomorphism-covariance}
Suppose two rational reference towers are related by \(F\)-isomorphisms
\(f_i\) that commute with all inclusions and identify the common
\(Z\). Use the corresponding dual isomorphisms in the same direction,
namely the duals of \(f_i^{-1}\), and transport the full characters,
parameters, and \(L\)-embeddings through them. Then the construction
and the recovered relative records are covariant under these maps.
If the transported \(L\)-embeddings are subsequently changed tamely
with their connected-dual parts fixed, apply
Theorem~\ref{sc:tame-embedding-compensation} as well.
\end{proposition}
\begin{proof}
On rigid torsors, extension of structure group through \(f_i\)
commutes with extension through the tower inclusions. It therefore
carries each relative reduction to a relative reduction in the other
tower, including its isomorphism with the ambient torsor. In cocycle
coordinates the rule is
\(z_i\mapsto f_i(z_i)\) and
\(g_i\mapsto f_{i+1}(g_i)\). Applying \(f_{i+1}\) to
\eqref{sc:eq-relative-cocycle} proves the required identity.
Functoriality of rigid duality gives exactly the transported central
characters. At the dual level, the isomorphisms carry the full
centralizers and their common covers to the corresponding groups;
pullback of irreducibles gives the same enhancement transport.
The torus LLC is functorial for the quotient-torus isomorphisms, so
the underlying recursive parameters transform in the same way.
The last assertion follows from the stated theorem.
\end{proof}

This is covariance, including for outer isomorphisms. An outer
automorphism is not thereby declared to act trivially on the seed or
its enhancement. Invariance applies after identifying the actual
transported data, and ambient rational conjugacy is already part of
the equivalence of relative records.

\begin{theorem}[Endpoint transport is independent of a supplied factorization]
\label{sc:endpoint-factorization-comparison}
Consider two rational reference chains from the same marked bottom
group \(H\) to the same ambient group \(G\). Suppose their
composite reference embedding \(j:H\hookrightarrow G\) is the
same, their dual chains have the same composite connected-dual
inclusion, and they use the same common \(Z\). Suppose valid
parameter constructions along the chains have identified bottom
enhanced seeds and, in these same coordinates, equal final parameter
maps after any specified compensated twists. Use the raw transports in
the common reference, or transport both systems of edge maps to it by
Theorem~\ref{sc:coherent-origin-transport}. Assume the full centralizer
hypotheses at each step and compatible relative reductions.

Then both chains induce the same map from bottom to final refined
enhancement classes. They also give the same class of composite marked
relative realization in \(\operatorname{Rel}_j(\zeta_G)\), when
their bottom rigid classes have been identified. The assertion is
covariant under an actual isomorphism carrying all the stated
identifications, as in Proposition~\ref{sc:reference-isomorphism-covariance}.
\end{theorem}
\begin{proof}
For each chain the composite map on refined groups is the restriction
of the same lifted bottom-to-ambient dual inclusion. This follows
because each step is an inverse-image comparison inside the final
cover. Thus the two maps agree, and so do their transports of
irreducible classes. Central twists do not change any of the relevant
centralizers, so the same argument applies after the stipulated
compensations. If a final conjugacy is used, any two choices that
identify the same full parameter maps differ by its centralizer;
their lifts act by an inner automorphism on the final refined group.

For the rational records, compose the transporters as in
\eqref{sc:eq-relative-composition}, obtaining a reduction along
the same \(j\) in each case. Both have the same bottom rigid class.
The bijection \eqref{sc:eq-relative-fiber} therefore identifies
their classes modulo ambient rational conjugacy and source
isomorphisms. Covariance under actual isomorphisms was proved above.
\end{proof}

This comparison applies after identifying both endpoint maps and the
marked realization. The successive character factors restrict on wild
inertia to the center of the bottom dual group. Their product can have
several distinct jumps in its centralizer flag, each of which a strict
factorization must realize by a proper stage
(Proposition~\ref{sc:wild-flag-recovery}).

\subsection{Affine comparison of rational reference embeddings}
\label{sc:affine-reference-change}

Now let \(j,j':H\hookrightarrow G\) be two twisted-Levi \(F\)-embeddings
of the same marked reference group, related by
\(j'=\Ad(b)\circ j\) over \(\overline F\). Define
\begin{equation}\label{sc:eq-reference-defect}
 c_b(\sigma)=j^{-1}\bigl(b^{-1}\sigma(b)\bigr).
\end{equation}
Because both embeddings are defined over \(F\), the expression
\(b^{-1}\sigma(b)\) centralizes \(j(H)\), hence belongs to
\(j(Z(H))\). Thus \(c_b\) is a central Galois cocycle. In
positive characteristic the same expression is read after an fpqc
trivialization; it gives a class in \(H^1(F,Z(H))\) without any
smoothness assumption on the center. Denote its image in
\(\mathcal H_H\) by \([c_b]_H\), and put
\(\epsilon_b=D_H([c_b]_H)\).

\begin{proposition}[Change of rational reference is affine]
\label{sc:affine-reference-comparison}
For a fixed ambient rigid torsor, the change from \(j\) to \(j'\)
is given on relative cocycle records by
\begin{equation}\label{sc:eq-affine-record-change}
 (z_H,g)\longmapsto (c_b^{-1}z_H,\,g b^{-1}).
\end{equation}
It leaves the realized embedding unchanged, with the identical
underlying source inner form. On source classes and central characters
it induces
\begin{equation}\label{sc:eq-affine-character-change}
 \zeta_H\longmapsto[c_b^{-1}z_H],\qquad
 \chi_{\zeta_H}\longmapsto\epsilon_b^{-1}\chi_{\zeta_H}.
\end{equation}
Here \(\epsilon_b\circ t_j=1\), and the two embeddings induce
the same map on relevance centers. These changes depend only on the
marked reference comparison up to the indicated equivalences and
compose as changes of origin. If \(b\in G(F)\), the defect is trivial.
\end{proposition}
\begin{proof}
Writing \(\sigma(b)=b j(c_b(\sigma))\), compute
\[
 (gb^{-1})^{-1}z_G(w)\sigma(gb^{-1})
   =b j(z_H(w))\sigma(b)^{-1}
   =j'\bigl(c_b(\sigma)^{-1}z_H(w)\bigr).
\]
The central cocycle \(c_b\) is inflated from the base field, so the
band of the rigid cocycle is unchanged. Its adjoint image is trivial,
so the source inner form is unchanged as well. Finally
\(\Ad(gb^{-1})j'=\Ad(g)j\). These identities prove the first
assertion and its fpqc version.

For completeness, central multiplication corresponds to multiplication
of rigid-duality characters. Choose an elliptic maximal \(F\)-torus
\(S\subset H\). The map \(\mathcal H_S\to\mathcal H_H\)
is surjective, and \(Z(H)\subset S\). Lift the class of \(z_H\)
to \(\mathcal H_S\), and map the central class \([c_b]\) to
\(H^1(F,S)\). Centrality ensures that their product maps to
\([c_bz_H]\), regardless of the chosen source trivialization.
Rigid duality for the torus is a homomorphism, and its functoriality
for \(S\to H\) proves the asserted character formula.
The image of \([c_b]_H\) in \(\mathcal H_G\) is trivial,
since \(j(c_b)=b^{-1}\sigma(b)\). Equation
\eqref{sc:eq-perfect-duality} gives \(\epsilon_b\circ t_j=1\).
The two maps on dual centers agree because their root-data maps
differ by geometric inner conjugation.

If another \(b'\) gives the same marked comparison, then
\(b'=b j(h)\) for \(h\in Z(H)(\overline F)\), and the
defect changes by the central coboundary \(h^{-1}\sigma(h)\).
The records in \eqref{sc:eq-affine-record-change} change by
the corresponding source isomorphism. For consecutive comparisons
\(j_1=\Ad(b_1)j\) and \(j_2=\Ad(b_2)j_1\), their defects
satisfy \(c_{b_2b_1}=c_{b_1}c_{b_2}\), where the second defect
is computed relative to \(j_1\). Substitution shows that the record
maps compose. All these assertions are identities of descent data.
\end{proof}

This affine formula is necessary even when both reference groups are
quasi-split. It compares the same physical embedding with different
reference origins; it does not authorize identifying distinct rational
embeddings while holding the origin fixed. In Appendix~\ref{sc:rank-one},
changing between the two marked rational reference embeddings has
\(\epsilon_b=\lambda_2\). The same physical embedding is then
recorded by the two different source rigid classes in the two coordinate
systems. Simply retaining its old central-character label would fail
the new relevance condition.

\subsection{The additional enhancement comparison for a nontrivial origin shift}
\label{sc:origin-label-comparison}

An affine change in rigidification is different from the central
\emph{parameter} twists in Theorem~\ref{sc:tame-embedding-compensation}.
The latter preserve the rigidification and its relevance character.
For \eqref{sc:eq-affine-character-change}, keeping the same
parameter and transporting the enhancement by the identity on the
refined group cannot work when \(\epsilon_b\ne1\).

Here is an exact algebraic criterion for comparison by a linear
character. It separates existence of a correction from its normalization.

\begin{proposition}[Criterion and ambiguity for an origin correction]
\label{sc:origin-character-criterion}
Set \(A=\cS_{\phi_H}^{+,Z}\) and
\(a=a_H:\mathcal Z_H\to A\). There exists a linear character
\(\lambda_b:A\to\C^\times\) with
\begin{equation}\label{sc:eq-origin-character}
 \lambda_b\circ a=\epsilon_b^{-1}
\end{equation}
if and only if \(\epsilon_b\) is trivial on
\(a^{-1}([A,A])\). When it exists, the choices form a torsor under
\[
 \Hom\bigl(A/(a(\mathcal Z_H)[A,A]),\C^\times\bigr).
\]
Such a choice sends the old relevance sector to the new one by
\(\rho_H^+\mapsto\lambda_b\otimes\rho_H^+\).
To preserve cuspidality it suffices in addition that \(\lambda_b\)
be pulled back from \(\pi_0(Q_{H,\phi_H}^{+,Z})\) through the
unipotent-centralizer description of \(A\).
\end{proposition}
\begin{proof}
A linear character of \(A\) factors through its abelianization.
Thus \eqref{sc:eq-origin-character} is possible only if the
prescribed character is trivial on the inverse image of \([A,A]\).
Under that condition it defines a character on the image of
\(\mathcal Z_H\) in \(A^{\mathrm{ab}}\), and this character
extends to the finite abelian group \(A^{\mathrm{ab}}\). Its
extensions differ precisely by the displayed character group.
The relevance assertion follows by multiplication of central characters.
For the last assertion, tensoring by a character of \(\pi_0(Q)\)
does not change the constituents of the local system restricted to
\(Q^\circ\). Its cuspidality is therefore unchanged.
\end{proof}

Suppose a comparison character with the required cuspidality property
has been chosen. With \(\alpha\) fixed under the dual identifications,
the raw common-cover transport changes the ambient label by
\begin{equation}\label{sc:eq-origin-raw-output}
 (\lambda_b\otimes\rho_H^+)\circ\alpha^{-1}
   =(\lambda_b\circ\alpha^{-1})\otimes\rho_G^+.
\end{equation}
Its ambient relevance character is unchanged since
\(\epsilon_b\circ t_j=1\), but its irreducible label need not
be unchanged. A comparison preserving the old ambient label must
include the compensating normalization
\begin{equation}\label{sc:eq-normalized-origin-transport}
 \mathcal E_{j'}(\rho_H^{+\prime})
   =(\lambda_b^{-1}\otimes\rho_H^{+\prime})\circ\alpha^{-1},
 \qquad \rho_H^{+\prime}=\lambda_b\otimes\rho_H^+.
\end{equation}
Then \(\mathcal E_{j'}(\rho_H^{+\prime})=\mathcal E_j(\rho_H^+)\).
For consecutive reference changes the corresponding characters must
multiply after transport of the groups. Any residual character trivial
on \(a(\mathcal Z_H)\) must be shown to act trivially on the
relevant irreducible labels, or retained in the normalization data.
Equality only of the central characters does not prove this coherence.

The standard input for resolving this normalization is as follows.
For \(p\)-adic fields,
\cite[\S6.3, construction preceding Lemma~6.2]{KalethaIsocrystals}
constructs a linear comparison character for a change of rigidifying
cocycle, using the negative differential and a duality pairing; its
construction is parameter-side. The compatibility of that construction
with an LLC is the separate, conditional assertion of Lemma~6.2
of that paper.
Its application to the present affine defect must match the cocycle,
finite cover, sign, cuspidality, and tower composition with
\eqref{sc:eq-affine-record-change}--\eqref{sc:eq-normalized-origin-transport}.
An auxiliary source cover is not automatically the common \(Z\) of
the Yu tower. The following subsections establish the required descent
for ordinary central defects, using the all-characteristic finite
coefficient construction of \cite[\S5.4]{DillerySchwein}, and then
prove coherence along each such marked change of reference. The
resulting character specifies how to transport the individual label.

\subsection{Canonical characters for ordinary central changes of origin}
\label{sc:canonical-origin-character}

The defects in Proposition~\ref{sc:affine-reference-comparison} are
\emph{ordinary} classes in \(H^1(F,Z(H))\): their restriction to the
rigid gerbe band is trivial. This is stronger than an arbitrary change
of rigidification. It allows the comparison character to descend all
the way to the ordinary Weil-centralizer component group. In
particular, the auxiliary finite subgroup used to compute the character
need not be central in the ambient group, and it need not replace the
common \(Z\) in the construction.

We use the ordinary local Tate pairing
\[
 \langle\ ,\ \rangle_F:
 H^1(F,C)\times H^1(\Gamma_F,X^*(C))\longrightarrow\C^\times
\]
for a group of multiplicative type \(C\). In positive characteristic
the first factor means fppf cohomology. Normalize this pairing by
compatibility with the ordinary-cohomology restriction of the finite
coefficient pairing in \cite[Proposition~5.4.1]{DillerySchwein}.
We use the negative differential convention of its \S5.4 and of
\cite[\S6.3]{KalethaIsocrystals}. These conventions fix the inverse
in the formulas below. The pairing is functorial in \(C\); it can
also be expressed by the evaluation cup product and the local invariant
map, with this Tate--Nakayama convention.

For a connected complex reductive group \(\widehat H\), let
\(\widetilde{\widehat H}\) be its universal covering complex Lie
group. It is the product of the simply connected covering of its
derived group with the additive universal cover of its connected
central torus. Its discrete central kernel is canonically
\[
 M_H=\pi_1(\widehat H)=X^*(Z(H)).
\]
For example, this identification follows from
\(X_*(\widehat T)/\Z\Phi^\vee(\widehat H,\widehat T)
 =X^*(T)/\Z\Phi(H,T)\).
This universal cover is only a device for defining a boundary class;
all enhancements still use the finite algebraic common cover.

\begin{lemma}[The ordinary central boundary]\label{sc:universal-boundary}
Put \(Q_\phi=Z_{\widehat H}(\phi|_{W_F})\). The action
\(\operatorname{Ad}\circ\phi\) lifts uniquely to
\(\widetilde{\widehat H}\). For \(q\in Q_\phi\), choose a lift
\(\widetilde q\) and set
\[
 d_{\phi,\widetilde q}(w)
   =\widetilde q^{-1}\,{}^{\phi(w)}\widetilde q\in M_H.
\]
Its cohomology class defines a homomorphism
\begin{equation}\label{sc:eq-universal-boundary}
 \partial_{H,\phi}:\pi_0(Q_\phi)
       \longrightarrow H^1(\Gamma_F,M_H).
\end{equation}
It is unchanged by a central parameter twist and is covariant under
actual reference isomorphisms and the corresponding dual identifications.
\end{lemma}
\begin{proof}
Every automorphism of a connected Lie group lifts uniquely to its
universal covering group as an automorphism fixing the identity.
The lifted automorphisms consequently form an action. Inner
conjugation acts trivially on its central kernel, so the action on
\(M_H\) is its usual finite Galois action. The displayed expression
is a continuous Weil cocycle. Changing the lift changes it by a
coboundary. Centrality of the kernel proves multiplicativity in \(q\).

A path in \(Q_\phi\) lifts to the covering group, and for each
\(w\) its boundary along the path is a continuous function with values
in the discrete group \(M_H\), hence constant. Thus the boundary
class is trivial on \(Q_\phi^\circ\) and factors through the finite
group \(\pi_0(Q_\phi)\). In particular its class in
\(H^1(W_F,M_H)\) has finite order. Choose a finite Galois extension
\(E/F\) trivializing the action on \(M_H\). On \(W_E\) the
cocycle is a homomorphism. If its class has order dividing \(n\),
its restriction takes values in the finite group \(M_H[n]\), since
coboundaries vanish for this trivial action. It follows, by using the
finitely many cosets of \(W_E\), that the original cocycle has finite
image. The associated homomorphism to the semidirect product with the
finite action group has finite image and extends to \(\Gamma_F\).
This gives the asserted Galois class. Restriction of such classes to
\(W_F\) is injective by density, so the class is unambiguous.
Central twists have the same adjoint action; functoriality of universal
covers proves the remaining covariance assertions.
\end{proof}

The boundary and Tate-pairing construction on ordinary component groups
also occurs in \cite[\S9, desideratum~(3)]{GGP}; here it is used solely
as a parameter-side construction, with the rigidification sign convention
specified above.

\begin{theorem}[Canonical origin correction on every common cover]
\label{sc:canonical-origin-correction}
For \(c\in H^1(F,Z(H))\), define
\begin{equation}\label{sc:eq-canonical-origin-character}
 \Lambda_{c,\phi}(q)
   =\langle c,-\partial_{H,\phi}(q)\rangle_F,
       \qquad q\in\pi_0(Q_\phi).
\end{equation}
Pull this character back to \(\cS_{\phi}^{+,Z}\), and denote the
result by \(\lambda_{c,\phi}^{Z}\). Then:
\begin{enumerate}[label=\textup{(\roman*)},leftmargin=2.5em]
\item Tensoring by \(\lambda_{c,\phi}^{Z}\) takes the relevance
sector of \(z_H\) to that of \(cz_H\). In particular
\[
 \lambda_{c,\phi}^{Z}\circ a_H=D_H([c]_H).
\]
\item The twist preserves cuspidality, is compatible with enlargement
of \(Z\), and needs no auxiliary enlargement of the common ambient
central subgroup.
\item The character depends only on \(c\), is multiplicative in
\(c\), equals \(1\) on a central coboundary, and is unchanged by
central parameter twists. It has the covariance of
Lemma~\ref{sc:universal-boundary}.
\end{enumerate}
Consequently the correction in Proposition
\ref{sc:origin-character-criterion} exists canonically for every
actual reference defect: take
\begin{equation}\label{sc:eq-reference-canonical-lambda}
 \lambda_b=\lambda_{[c_b]^{-1},\phi_H}^{Z}.
\end{equation}
Its cuspidality condition is automatic. The finite-group obstruction
in that proposition therefore vanishes for these defects.
\end{theorem}
\begin{proof}
The boundary lemma and bilinearity prove that \(\Lambda_{c,\phi}\)
is a character with the properties in (iii). We verify both its
normalization and its descent by comparison with the finite coefficient
construction; this also treats non-smooth centers in positive
characteristic.

There is a finite central subgroup scheme \(C\subset Z(H)\)
containing \(Z\) and a class \(\widetilde c\in H^1(F,C)\)
mapping to \(c\). Indeed, let \(T\) be the maximal torus in the
group of multiplicative type \(Z(H)\). Choose an integer \(n_0\)
killing the finite quotient \(Z(H)/T\). Multiplication by \(n_0\)
maps \(Z(H)\) onto \(T\) in the fppf topology and has finite
kernel. The group \(H^1(F,T)\) has finite exponent. Multiplying
\(n_0\) by this exponent, and further if necessary to kill the
prescribed finite \(Z\), gives an isogeny \(Z(H)\to T\) whose
map on \(H^1\) is zero and whose kernel \(C\) contains \(Z\).
The cohomology sequence supplies \(\widetilde c\).

The dual isogeny \(\widehat{H/C}\to\widehat H\) has kernel
\(K_C\) canonically identified with \(X^*(C)\) as a finite
Galois module. The map from the universal cover to
\(\widehat{H/C}\) identifies its boundary in \(H^1(\Gamma_F,K_C)\)
with the image of \(\partial_{H,\phi}\) under
\(X^*(Z(H))\to X^*(C)\).
The pairing of the inflated ordinary class \(\widetilde c\) with
\(Z^1(\Gamma_F,K_C)\) annihilates coboundaries and is the ordinary
Tate pairing on cohomology. Functoriality therefore gives, for every
lift \(q_C\) of \(q\),
\begin{equation}\label{sc:eq-finite-origin-comparison}
 \langle\widetilde c,-d_C(q_C)\rangle
      =\langle c,-\partial_{H,\phi}(q)\rangle_F.
\end{equation}
The left side is precisely the finite coefficient comparison character
of \cite[\S5.4]{DillerySchwein}. Its restriction to the refined
central relevance group is the ratio of the characters for
\(\widetilde c z_H\) and \(z_H\), by compatibility of that
pairing with rigid Tate--Nakayama duality; in the function-field setting
see also \cite[\S5.4, diagram~(32)]{DilleryLimits}. This compatibility is
parameter-side and does not use the conjectural LLC discussed there.

Equation~\eqref{sc:eq-finite-origin-comparison} shows directly
that this character factors through \(\pi_0(Q_\phi)\). In
particular it is independent of the auxiliary \(C\) and its lift,
and descends to every prescribed common cover. Equivalently, on the
kernel of the map to a smaller cover its differential is a coboundary,
which an ordinary class annihilates. Inflate an enhancement from the
\(Z\)-cover to the \(C\)-cover and apply the finite coefficient
relevance identity. Surjectivity of the maps on the corresponding
central component groups then descends that identity to the
\(Z\)-cover. Since \(\widetilde c\) and \(c\) have the same
image in central cohomology, the resulting rigid class is \([cz_H]\).
This proves (i) with the indicated sign.

The same character pulls back from \(\pi_0(Q_\phi)\) to
\(\pi_0(Q_\phi^{+,Z})\). Thus it is trivial on the connected
Weil centralizer and leaves unchanged the cuspidal constituents of
the unipotent local system on that connected group. This proves (ii),
as in the parameter-side argument of
\cite[Proposition~5.4.2]{DillerySchwein}. Cover compatibility also
follows immediately from its factorization through the ordinary group.
Finally \(c=[c_b]^{-1}\) gives the relevance ratio
\(\epsilon_b^{-1}\) of Proposition
\ref{sc:affine-reference-comparison}, proving the last assertions.
\end{proof}

In the unramified torus example of Appendix~\ref{sc:rank-one}, \(Q_{\phi_H}=\widehat H^{\Gamma_F}=\mu_2\).
The nontrivial class in \(H^1(F,H)\) gives its nontrivial character.
Pullback along \(\mu_4\to\mu_2\) is \(\lambda_2\), exactly
the required interchange of the two source sectors. Notice that the
correction can be nontrivial even if its image in ambient rigid
cohomology is neutral. Ambient relevance alone does not determine the
individual enhancement label.

\subsection{Coherent normalization along a changed reference tower}
\label{sc:coherent-reference-tower}

The preceding characters remove the arbitrary extension choices of
Proposition~\ref{sc:origin-character-criterion}. We retain an
initial reference presentation and its enhancement transport. The
following theorem gives its coherent expression in other marked
reference presentations. It does not select the initial presentation
from a pinning alone.

\begin{theorem}[Coherent change of marked rational reference]
\label{sc:coherent-origin-transport}
Let \(j_i:G^i\hookrightarrow G^{i+1}\) be a compatible rational
reference chain, and let
\(j_i'=\Ad(b_i)j_i\) be another chain on the same marked groups.
Use the corresponding fixed connected-dual inclusions and the common
\(Z\); compensate any tame changes of Weil part as in
Theorem~\ref{sc:tame-embedding-compensation}. Put
\[
 c_i(\sigma)=j_i^{-1}(b_i^{-1}\sigma(b_i)),\qquad
 y_d=1,\qquad y_i=j_i^{-1}(y_{i+1})c_i^{-1}.
\]
The central inclusion \(Z(G^{i+1})\subset j_i(Z(G^i))\) makes
this recursion meaningful. Then
\[
 z_i'=y_i z_i,\qquad g_i'=g_i b_i^{-1}
\]
is a compatible chain of relative records with the same ambient rigid
cocycle and the same realized marked embeddings. All cocycle formulas
have the fpqc interpretation of Section~\ref{sc:relative-realization}
in positive characteristic. Define
\[
 T_i(\rho)=\lambda_{[y_i],\phi_i}^{Z}\otimes\rho,
 \qquad
 E_i(\rho)=\rho\circ\alpha_i^{-1}.
\]
The normalized enhancement transition for the new reference is
\begin{equation}\label{sc:eq-coherent-edge}
 E_i'=T_{i+1}\,E_i\,T_i^{-1}
     =E_i\circ\bigl(\lambda_{[c_i],\phi_i}^{Z}\otimes-\bigr).
\end{equation}
It preserves relevance and cuspidality, and
\begin{equation}\label{sc:eq-coherent-tower}
 E_{d-1}'\cdots E_0'\,T_0=E_{d-1}\cdots E_0.
\end{equation}
These comparisons compose for successive changes of the same marked
reference data. A closed comparison returning to the identical marked
reference has trivial correction. They are compatible with actual
reference isomorphisms, central parameter compensation, and enlargement
of the common cover. No comparison is asserted after forgetting the
source markings or quotienting by an unspecified normalizer action.
\end{theorem}
\begin{proof}
Centrality of \(y_{i+1}\) gives
\[
 (g_i b_i^{-1})^{-1}z_{i+1}'\sigma(g_i b_i^{-1})
 =y_{i+1}\,j_i'(c_i^{-1}z_i)=j_i'(y_i z_i).
\]
The adjoint cocycles are unchanged, and
\(\Ad(g_i b_i^{-1})j_i'=\Ad(g_i)j_i\), proving the record
assertion. All \(y_i\) are ordinary central cocycles, so no band or
common-cover change is involved.

We need the naturality of the canonical character for one step. The
lifted dual inclusion between universal covers carries the source
boundary to the target boundary through
\[
 X^*(Z(G^i))\longrightarrow X^*(Z(G^{i+1})).
\]
Indeed, on the source dual group the action through \(\phi_{i+1}\)
is the action through \(\phi_i\); the extra source-central factor
acts trivially there. Uniqueness of the lift gives the same identity
on the universal covers. For a class \(y\) in the ambient center,
functoriality of the Tate pairing therefore gives
\begin{equation}\label{sc:eq-origin-step-naturality}
 \lambda_{[j_i^{-1}(y)],\phi_i}^{Z}
     =\lambda_{[y],\phi_{i+1}}^{Z}\circ\alpha_i.
\end{equation}
Multiplicativity and the recursion for \(y_i\) now prove the second
formula in \eqref{sc:eq-coherent-edge}. The first formula takes
the new source relevance sector back to the old one, applies the old
transition, and then gives the new target sector. Each operation
preserves cuspidality. On multiplying the edge identities, the
intermediate \(T_i\)'s cancel, and \(T_d=1\), proving
\eqref{sc:eq-coherent-tower}.

For successive comparisons, the central defects multiply, as proved
in Proposition~\ref{sc:affine-reference-comparison}. The
\(y_i\)'s and their canonical characters therefore multiply after
the prescribed identifications. If a composite comparison fixes the
marked reference embedding, its transporter centralizes that embedding;
its defect is then a central coboundary. The character is trivial by
Theorem~\ref{sc:canonical-origin-correction}. This proves
coherence, not just equality of relevance characters. Naturality of
the boundary proves compatibility with actual reference isomorphisms;
central parameter twists do not change it. Cover compatibility was
proved in the same theorem.
\end{proof}

Within a fixed reference presentation, equality of the refined groups
gives the linear transport. Under a change of origin, the marked
reference comparison determines the character in
\eqref{sc:eq-coherent-edge}. Its restriction to the central subgroup
alone need not determine its action on the enhancement.

\section{Recovery of the dual tower and comparison in a fixed presentation}
\label{sc:constructed-tower-recovery}

The wild centralizer flag recovers every valid strict dual tower.
We use it to compare two supplied constructions, before invoking
reverse extraction from an arbitrary parameter. A valid strict construction here has the tower of
Definition~\ref{sc:def-data}, strictly increasing positive depths
\(0<r_0<\cdots<r_{d-1}\), a wild-trivial seed, embeddings whose
dual-group parts on wild inertia are trivial, and, at each stage,
\begin{equation}\label{sc:eq-leading-full-cut}
 Z_{\widehat{G^{i+1}}}\bigl(\kappa_i(I_F^{r_i})\bigr)
       =\widehat{G^i}.
\end{equation}
Here \(\kappa_i\) is trivial on \(I_F^{r_i+}\), and the other
validity conditions of the forward theorem are retained. This is the
full leading-image cut-down hypothesis used by this comparison argument.
Proposition~\ref{sc:generic-full-cutdown} proves it from
genericity in the prime-to-Weyl range.

\begin{proposition}[The wild centralizer flag recovers a valid strict tower]
\label{sc:wild-flag-recovery}
Let \(\Psi\) be the final underlying parameter of such a construction,
including its possible final ambient-central twist. Identify the dual
tower inside \(\widehat G\) by its supplied connected-dual inclusions.
For \(s\ge0\), put
\[
 k(s)=\#\{i<d:r_i\le s\}.
\]
Then the \emph{full} algebraic centralizer satisfies
\begin{equation}\label{sc:eq-wild-flag}
 Z_{\widehat G}\bigl(\Psi(I_F^{s+})\bigr)
       =\widehat{G^{k(s)}}.
\end{equation}
In particular \(\Psi\) determines the embedded strict dual tower and
all its positive-depth transition values. An ambient-central twist
cannot change this flag or its jumps.
\end{proposition}
\begin{proof}
On wild inertia the recursive formula is the product of the embedded
\(\kappa_i\)'s and the final ambient-central factor: the seed and
all embedding Weil parts contribute the identity. If \(i<k(s)\),
then \(\kappa_i\) is trivial on \(I_F^{s+}\). For
\(i\ge k(s)\), its values lie in
\(Z(\widehat{G^i})\subset Z(\widehat{G^{k(s)}})\).
Thus the right side of \eqref{sc:eq-wild-flag} centralizes the
indicated image.

Conversely suppose \(k(s)<d\) and an element centralizes that image.
The subgroup \(I_F^{r_{d-1}}\) is contained in \(I_F^{s+}\).
On it every earlier factor is trivial, and the final factor is central
in \(\widehat G\); hence \eqref{sc:eq-leading-full-cut}
places the element in \(\widehat{G^{d-1}}\). Inside this group
the factor \(\kappa_{d-1}\), and every later factor, is central.
If \(k(s)<d-1\), restrict to \(I_F^{r_{d-2}}\) and apply the
same argument inside \(\widehat{G^{d-1}}\). Descending induction
places the element in \(\widehat{G^{k(s)}}\). If \(k(s)=d\),
the wild tail is already ambient-central, so its centralizer is all of
\(\widehat G\). This proves the formula. Strictness of the tower
makes its changes occur exactly at the \(r_i\)'s.
\end{proof}

\begin{theorem}[Comparison after alignment in a fixed presentation]
\label{sc:realized-comparison-generation}
Consider two valid strict constructions. Fix an identification of their
actual marked rational towers and relative records, including the common
\(Z\), and express the two constructions in the same rational reference
with its common-cover transports. The coherent changes of
Theorem~\ref{sc:coherent-origin-transport} may be used to reach these
coordinates. Suppose their connected-dual inclusions agree and that an
element \(h\) of the embedded bottom dual group \(\widehat{G^0}\)
conjugates the first final parameter map to the second.

Then the differences of full central factors and tame Weil embeddings
satisfy the compensation identity \eqref{sc:eq-recovered-compensation}.
The required bottom twist has depth at most zero. If a lift of \(h\)
also carries the first endpoint label to the second, inverse transport
recovers the corresponding bottom enhancement. The common reference then identifies the relevance sectors and the
relative records attached to the two bottom labels.
\end{theorem}
\begin{proof}
Conjugate the first seed by the preimage of \(h\) in its dual group
and transport the label by a lift in its common cover. Each stage is
composition with a homomorphism over \(W_F\) whose connected restriction
is the specified inclusion, followed, if necessary, by a central twist.
Thus this seed conjugacy induces precisely conjugation by \(h\) at the
endpoint. The final maps now agree literally, including their
\(\SL_2\)-parts, with the rational tower and its factors still in the
given coordinates. Proposition~\ref{sc:wild-flag-recovery} gives the
same transition depths.
Lemma~\ref{sc:L-embedding-difference} writes
\(\iota_i'=\delta_i\iota_i\), with \(\delta_i\) tame and
source-central. Equality of the final full parameters gives
\(\phi_d'=v_d\phi_d\) with
\(v_d=\kappa_d\kappa_d'^{-1}\). Descending through the recursions,
define
\begin{equation}\label{sc:eq-recovered-compensation}
 v_i=\kappa_i'^{-1}\delta_i^{-1}\kappa_i\nu_i(v_{i+1}).
\end{equation}
Each factor is a cocycle in the full source center. Injectivity of
\(\iota_i\) shows successively that \(\phi_i'=v_i\phi_i\).
Since the dual parts of both \(\phi_0\) and \(\phi_0'\) are trivial
on wild inertia, \(v_0\) is tame. This gives the asserted compensation identity on parameters.

The enhancement maps on the fixed dual flag are isomorphisms, and
central parameter twists leave them unchanged. Theorem
\ref{sc:coherent-origin-transport} compares the rational origins
with their canonical characters. Hence equality of the final enhanced
labels determines the required comparison of bottom labels by inverse
transport. At each quasi-split reference term its central character
determines its unique rigid class by rigid duality. Theorem
\ref{sc:enhancement-realization} realizes these compatible classes.
For the record comparison, a compatible chain of reductions is
isomorphic to the chain obtained by extending its bottom torsor along
all reference inclusions, equipped with its composite identification
with the ambient torsor. Theorem~\ref{sc:one-step-realization}
for the composite reference embedding therefore gives the asserted
uniqueness modulo the stipulated equivalences. The affine formulas
supply the same conclusion after changing references.
\end{proof}

The hypothesis on \(h\) is additional to the recovery of the embedded
flag. The following geometric observation describes the information
about quasi-split terms available from the flag and its markings.

\begin{proposition}[Geometric comparison of marked reference diagrams]
\label{sc:quasisplit-reference-matching}
Conjugation of final parameter maps identifies their embedded dual
flags and the induced outer Weil action on each term. If, in addition,
an identification of their full marked root-data diagrams is given,
compatible with the based Galois actions and common \(Z\) and fixing
the ambient marking, it identifies the quasi-split reference groups
over \(F\) and their corresponding embeddings up to geometric
conjugacy.
\end{proposition}
\begin{proof}
Proposition~\ref{sc:wild-flag-recovery} identifies the flags. Normality
of \(I_F^{s+}\) makes every term stable under the final Weil parameter.
Its action on the based root datum is the action of the reference
group. Indeed the later character factors centralize that term, and
the stage parameter acts through the given outer action. Conjugation
of final maps therefore intertwines these actions.

For the second assertion choose Galois-stable pinnings of the
quasi-split groups. The supplied maps of based root data induce
pinned isomorphisms equivariant for the descent actions.
Applying this at each term gives the required \(f_i\), with the
ambient marking fixed. The common-central-subgroup condition on the
root-data diagrams ensures that the \(f_i\)'s identify the same
\(Z\). After these isomorphisms, corresponding marked embeddings
have the same absolute root-data map. Over \(\overline F\),
align maximal tori and compatible Borels of the Levi subgroups.
The remaining discrepancies on root-group parametrizations are removed
by torus conjugation. Hence the embeddings are conjugate by an element
\(b_i\in G^{i+1}(\overline F)\). This proves the asserted geometric
comparison. A rational normalizer representative inducing a prescribed
action on the embedded characters is a further condition, as
Proposition~\ref{sc:rational-weyl-obstruction} illustrates.
\end{proof}

\section{Reverse construction and terminating factorization}
\label{sc:reverse}\label{sec:exhaustion-large-p}

In Section~\ref{sc:constructed-tower-recovery} we started with a
parameter already constructed from a tower and recovered the dual
groups in that tower from its ramification. We now start with an
arbitrary relevant cuspidal rigidly enhanced parameter. If it has
positive depth \(r\), the centralizer of its image on \(I_F^r\)
determines the dual group in which to seek a parameter of smaller
depth. A central projection of the entire wild image extends to a
Weil cocycle. Removing this factor lowers the depth, and the
centralizer comparison transports the enhancement to the group at
the preceding step. Its relevance character determines the source
rigidification, which we realize in the prescribed ambient form.

At each step the remaining wild image is a proper quotient of the
preceding finite wild image. Its order therefore strictly decreases
until the parameter has depth zero. Reading the steps in the opposite
order constructs the original parameter with its given enhancement.

The target is the class of Definition~\ref{sc:rigid-cuspidality},
with relevance \eqref{sc:eq-rigid-relevance} for the prescribed
ambient rigid form. It is specified independently of the forward
construction. The parameter-side inputs from Dat--Fintzen are used
without a representation-side correspondence.
\subsection{A full central factor with controlled wild restriction}
\label{sc:reverse-central-factor}

\begin{lemma}[The central projection on finite \(p\)-images]
\label{sc:central-p-projection}
Let \(L\) be a connected complex reductive group with
\(p\nmid|W(L)|\), put \(A=Z(L)^\circ\), and write
\(q:L\to L_{\mathrm{ab}}\). The isogeny
\(q_A:A\to L_{\mathrm{ab}}\) induces an isomorphism
\[
       A[p^\infty]\xrightarrow{\sim}L_{\mathrm{ab}}[p^\infty].
\]
For every finite \(p\)-subgroup \(B\subset L\) the formula
\begin{equation}\label{sc:eq-central-projection}
 c_L(b)=(q_A|_{A[p^\infty]})^{-1}(q(b))
\end{equation}
defines a homomorphism \(B\to A[p^\infty]\), natural for
automorphisms of \(L\). It is the identity on \(B\cap Z(L)\).
The homomorphism
\[
                     b\longmapsto b\,c_L(b)^{-1}
\]
has image in \(L_{\mathrm{der}}\) and kills \(B\cap Z(L)\).
\end{lemma}
\begin{proof}
The kernel of \(q_A\) is contained in \(Z(L_{\mathrm{der}})\),
whose order is prime to \(p\). Hence the isogeny is an
isomorphism on \(p\)-primary torsion. Also
\(\pi_0Z(L)\) is a quotient of \(Z(L_{\mathrm{der}})\),
because \(L=A L_{\mathrm{der}}\); it has order prime to \(p\).
Thus every central \(p\)-element lies in \(A\), where the
formula is the identity. All maps used to define \(c_L\) are
intrinsic. The last assertion follows by applying \(q\), and
multiplicativity follows because \(c_L\) has central image.
\end{proof}

\begin{lemma}[Extending a prescribed wild central homomorphism]
\label{sc:wild-central-extension}
Let \(A\) be a complex torus with a finite tame \(W_F\)-action.
Every continuous \(W_F\)-equivariant homomorphism
\[
                    c:P_F\longrightarrow A[p^\infty]
\]
with finite image extends to a continuous cocycle
\(\kappa:W_F\to A[p^\infty]\). Its restriction is exactly
\(c\), not merely conjugate to \(c\). Such an extension has
all its values in a finite \(p\)-subgroup of \(A\). Any two extensions
differ by a tame central cocycle. In particular, if \(c\)
factors through a fixed finite quotient of \(P_F\), no new wild
ramification is introduced by the extension.
\end{lemma}
\begin{proof}
Give \(A[p^\infty]\) the discrete topology. The tame inertia
group \(I_F/P_F\) is pro-prime-to-\(p\), so its positive-degree
continuous cohomology with these coefficients vanishes. This follows
by averaging over finite quotients, whose orders are invertible on
\(p\)-primary torsion. The quotient of \(W_F/P_F\) by tame
inertia is \(\mathbb Z\), of cohomological dimension one.
The Hochschild--Serre sequence therefore gives
\[
 H^2(W_F/P_F,A[p^\infty])=0.
\]
Inflation--restriction for \(P_F\subset W_F\) now makes
\[
 H^1(W_F,A[p^\infty])\longrightarrow
 \Hom_{\mathrm{cts}}(P_F,A[p^\infty])^{W_F}
\]
surjective. Wild inertia acts trivially on the coefficients, so
its coboundaries are trivial: the prescribed restriction is an
actual equality for any representing cocycle. Two such cocycles
have ratio trivial on \(P_F\), proving the comparison assertion.
Continuity makes the set of values on compact inertia finite.
The values on a Frobenius element and their translates under the
finite action, together with the inertia values, generate a finite
\(p\)-subgroup of \(A\) containing every value of \(\kappa\).
The final assertion follows from the exact restriction formula.
\end{proof}

The extension in Lemma~\ref{sc:wild-central-extension} has a
different role from the tame twisted-Levi \(L\)-embedding discussed
in Section~\ref{sc:reverse-reference}. In the reverse construction
we prescribe the entire wild restriction
\[
 \kappa|_{P_F}=c_L\circ\phi|_{P_F},
 \qquad A=Z(L)^\circ,
\]
using Lemma~\ref{sc:central-p-projection}. After factorization
through a tame embedding, division by \(\kappa\) therefore
replaces the wild image by
\(b\mapsto b\,c_L(b)^{-1}\). This kills the top-depth
central image and gives the proper quotient used for termination.
A tame embedding preserves the wild restriction when the parameter
is pulled back; it does not by itself supply this prescribed central
factor. Remark~\ref{sc:kaletha-reverse-embedding} explains the
complementary embedding input.

\subsection{The reference Levi and its tame \texorpdfstring{\(L\)}{L}-embedding}
\label{sc:reverse-reference}

\begin{proposition}[Rational reference and tame factorization of the extension]
\label{sc:parameter-levi-reference}
Let \(M\) be tame, let \(\phi:W_F\times\SL_2(\C)\to{}^LM\)
be an admissible parameter, and assume
\(p\nmid|W_{\mathrm{abs}}(M)|\). For a positive ramification
index \(r\), put
\[
 B=\phi(P_F),\qquad K=\phi(I_F^r),\qquad
 L=Z_{\widehat M}(K),\qquad \mathcal L_r=L\,\phi(W_F).
\]
Wild images here mean dual-group components. Then \(L\) is a
connected Levi, \(B\subset L\), and \(\mathcal L_r\) has
a tame splitting preserving a pinning of \(L\). Thus there is
a tame quasi-split group \(H^*\), a tame \(L\)-embedding
\[
 \iota:{}^LH^*\longrightarrow{}^LM,
 \qquad \iota(\widehat H^*)=L,\quad
 \iota({}^LH^*)=\mathcal L_r,
\]
through which the full parameter \(\phi\) factors. Inside a
quasi-split inner form \(M^*\) there is a rational quasi-split
twisted Levi reference \(j:H^*\hookrightarrow M^*\) with
this full dual root-data marking, including the ambient common
\(Z\). If \(\phi\) is discrete, this reference inclusion is
elliptic and
\begin{equation}\label{sc:eq-reverse-fixed-centers}
       (Z(L)^{W_F})^\circ=(Z(\widehat M)^{W_F})^\circ.
\end{equation}
\end{proposition}
\begin{proof}
Lemma~\ref{sc:p-centralizers} makes \(B\) abelian and
toral, and makes \(L\) a connected Levi. In particular \(B\)
centralizes \(K\) and is contained in \(L\). Normality of
\(I_F^r\) shows that \(\phi(W_F)\) normalizes \(L\).
The \(\SL_2\)-image centralizes \(B\), hence lies in \(L\).
The finite outer action on \(L\) is tame because the action
of \(\phi(P_F)\) on \(L\) is inner.

The pinning-preserving splitting, trivial on the dual-group part
of wild inertia, is the lemma in \cite[\S2.3.1]{DatFintzen}. Apply
\cite[Proposition~2.1.5]{DatFintzen} to the resulting stable
Levi-center embedding in the dual group to obtain the stated
rational quasi-split reference. That proposition constructs the
full dual root datum, not only an inclusion of centers. Taking
the quotient by the given ambient \(Z\) provides the compatible
cover diagram as in Proposition~\ref{sc:reference-root-data}.

Finally \((Z(L)^{W_F})^\circ\) centralizes the full parameter:
the Weil action on this center is induced by \(\phi\), and the
\(\SL_2\)-image lies in \(L\). If \(\phi\) is discrete,
this torus is contained in \((Z(\widehat M)^{W_F})^\circ\).
The reverse inclusion always holds. Equation
\eqref{sc:eq-reverse-fixed-centers} follows, and
Lemma~\ref{sc:elliptic-fixed-centers} proves ellipticity.
\end{proof}

\begin{remark}[Kaletha's twisted-Levi embeddings and reverse factorization]
\label{sc:kaletha-reverse-embedding}
For an already realized tame twisted Levi inclusion, existence of
a compatible tame \(L\)-embedding is Kaletha's twisted-Levi
embedding lemma, \cite[Lemma~5.2.8]{Kaletha}.
The lemma does not impose regularity on a depth-zero seed. Its proof
removes a central \(2\)-cocycle obstruction and chooses the
correction to be tame. For the arbitrary-characteristic scope of
this paper, the Dat--Fintzen splitting used in the preceding proof
remains available; see also the construction from
\cite[Corollary~6.7]{KalethaCovers} in
Theorem~\ref{sc:generic-stage-theorem}.

Once the rational reference \(j:H^*\hookrightarrow M^*\)
and its full dual marking have been obtained, a tame embedding
\(\iota:{}^LH^*\to{}^LM\) with that marking automatically
has the image required in Proposition~\ref{sc:parameter-levi-reference}.
Here is the comparison argument. Identify \(\widehat H^*=L\)
and put \(\xi(w)=\iota(1,w)\). The automorphisms of \(L\)
induced by \(\phi(w)\) and \(\xi(w)\) have the same
outer class. Thus, for each \(w\), there is \(l_w\in L\)
such that
\[
 \Ad(\phi(w))|_L=\Ad(l_w\xi(w))|_L.
\]
The element \((l_w\xi(w))^{-1}\phi(w)\) has trivial Weil
projection and centralizes \(L\). Since a Levi subgroup satisfies
\(Z_{\widehat M}(L)=Z(L)\subset L\), it follows that
\(\phi(w)\in L\xi(w)\), and hence
\[
              \iota({}^LH^*)=L\phi(W_F)=\mathcal L_r.
\]
The \(\SL_2\)-image lies in \(L\), so the full parameter
factors through \(\iota\). No coherent choice of the individual
\(l_w\) is needed for this image comparison.

This gives an alternative proof of the embedding and factorization
part once the correctly marked rational reference is available.
The reference itself is still supplied by the rational-realization
input above. The prescribed central factor in Section
\ref{sc:reverse-central-factor}, and the common-\(Z\) enhancement
and relative-rigidification comparisons, are separate ingredients.
\end{remark}

\subsection{One reverse step for an arbitrary cuspidal rigid parameter}
\label{sc:reverse-step-proof}

\begin{theorem}[Depth-lowering factorization with rigid enhancements]
\label{sc:reverse-step-theorem}\label{prop:large-p-peeling}\label{prop:inverse-onestep-descent}
Let \(M\) be a tame connected reductive group over \(F\), with
\(p\nmid|W_{\mathrm{abs}}(M)|\), and fix its ambient rigid
inner twist represented using \(Z\). Let
\((\phi_M,\rho_M^+)\) be any relevant cuspidal rigidly enhanced
parameter in the convention of Definition~\ref{sc:rigid-cuspidality},
of positive depth \(r\). There exist a rational elliptic tame
twisted Levi \(H\subset M\), a compatible tame \(L\)-embedding,
a full quotient character \(\bar\chi_H\) of depth \(r\), and
a relevant cuspidal rigidly enhanced parameter
\((\phi_H,\rho_H^+)\) of depth less than \(r\), such that
\begin{equation}\label{sc:eq-proved-reverse-step}
 \phi_M=\kappa\,\iota(\phi_H),\qquad
 [\kappa]=\mathcal L_{D_H}(\bar\chi_H),\qquad
 \rho_H^+=\rho_M^+\circ\alpha,
\end{equation}
after compatible choices of representatives. Here
\[
 \widehat H=Z_{\widehat M}(\phi_M(I_F^r)),\qquad
 \alpha:\mathcal S_{\phi_H}^{+,Z}\xrightarrow{\sim}
                     \mathcal S_{\phi_M}^{+,Z}.
\]
If \(\widehat H\ne\widehat M\), the pullback of
\(\bar\chi_H\) is \(M\)-generic of depth \(r\) in the full
Yu sense. If \(\widehat H=\widehat M\), this is a central
step with unchanged group, enhancement, and rigidification.
In both cases the wild restriction of \(\phi_H\) factors
through a proper quotient of the finite wild image of \(\phi_M\).
The same \(Z\) works throughout, and the relative marked
realization is retained.
\end{theorem}
\begin{proof}
The parameter is discrete by the definition of the target. Set
\(B=\phi_M(P_F)\), \(K=\phi_M(I_F^r)\), and
\(L=Z_{\widehat M}(K)\). The group \(K\) is nontrivial,
whereas \(\phi_M(I_F^{r+})=1\). Apply
Proposition~\ref{sc:parameter-levi-reference} and initially
work with its quasi-split reference \(H^*\). Put
\(\phi_L=\iota^{-1}\phi_M\); since the embedding is tame,
its wild restriction is the original homomorphism into \(B\subset L\).
Write \(A=Z(L)^\circ\).

Apply the intrinsic central projection of
Lemma~\ref{sc:central-p-projection} to \(B\). The composite
\begin{equation}\label{sc:eq-reverse-wild-factor}
 c:P_F\xrightarrow{\phi_L}B\xrightarrow{c_L}A[p^\infty]
\end{equation}
is equivariant for the pinned \(W_F\)-action on \(A\):
equivariance for conjugation by \(\phi_L\) follows from
naturality, and inner conjugation in \(L\) acts trivially on
\(A\) and \(L_{\mathrm{ab}}\). Lemma
\ref{sc:wild-central-extension} extends it to a full cocycle
\(\kappa:W_F\to A[p^\infty]\). Define
\[
                   \phi_H=\kappa^{-1}\phi_L
\]
on the Weil group, leaving the \(\SL_2\)-part unchanged.
This is a homomorphism because \(\kappa\) is central for the
source action; Frobenius semisimplicity and admissibility follow
from the power calculation in
Theorem~\ref{sc:generic-stage-theorem}, applied in reverse.

The group \(K\) is central in \(L\), so \(c_L\) is the
identity on \(K\). Hence \(\kappa|_{I_F^r}=\phi_M|_{I_F^r}\)
as actual leading maps, and \(\phi_H(I_F^r)=1\).
Both factors are trivial on \(I_F^{r+}\). More strongly, on all
wild inertia the new homomorphism is the composite
\begin{equation}\label{sc:eq-residual-wild-map}
 P_F\longrightarrow B\longrightarrow L_{\mathrm{der}},
       \qquad b\longmapsto b\,c_L(b)^{-1}.
\end{equation}
The second arrow is a homomorphism killing the nontrivial subgroup
\(K\). It follows that the new wild image is a proper quotient
of \(B\), and that its depth is strictly less than \(r\).
The latter also follows directly from the finite ramification
filtration on this quotient. The full cocycle \(\kappa\) has
depth exactly \(r\), since its leading restriction is \(K\).

Since \(A=\widehat{D_{H^*}}\), the torus LLC gives a full
character \(\bar\chi_{H^*}\) of \(D_{H^*}(F)\). The
positive-depth filtration exactness used in
Lemma~\ref{sc:graded-genericity} shows that its pullback to
\(H^*(F)\) has exact depth \(r\). Every root outside \(L\)
is nontrivial on \(K\), by the definition of its full
centralizer. Thus that lemma proves genericity, including
\(\mathrm{GE2}\), when \(L\) is proper. This uses a full
quotient character constructed from \(\kappa\), not an
extension chosen only from a graded character.

Any element centralizing \(\phi_M\), or just its Weil
restriction, centralizes \(K\) and therefore lies in \(L\).
Inside \(L\), the factors \(\kappa\) cancel from the
centralizer equation. Consequently both ordinary centralizer
equalities of Lemma~\ref{sc:ordinary-centralizers} hold.
Equation~\eqref{sc:eq-reverse-fixed-centers} proves that
\(\phi_H\) is discrete. The common-cover identification gives
\(\alpha\), and the equal Weil centralizers and unchanged
unipotent data preserve cuspidality in the reverse direction as
well as in the forward direction.

It remains to realize the relevance sector inside the prescribed
ambient rigid form. Transport \(\rho_M^+\) back by \(\alpha\)
and apply Theorem~\ref{sc:enhancement-realization} to the
reference \(j:H^*\hookrightarrow M^*\). Its source central
character determines a rigid class represented with the same
\(Z\), and a marked embedding of its inner form \(H\) into
\(M\). Keep this reduction in the relative record. The
character \(\bar\chi_{H^*}\) transports to \(D_H(F)\) via
the rational abelian-quotient identification of the inner twist.
The embedding is elliptic. It is tame: the connected center of
\(H\) is a tame torus, and over a tame extension splitting
that center and \(M\), its centralizer \(H\) is a split Levi.
Genericity and depth follow from the same root-data and quotient
character calculation on this actual inner form.

If \(L=\widehat M\), take the identity reference embedding
and the identity \(L\)-embedding from the start. The cocycle
is then central in \(\widehat M\), so all centralizers,
enhancement labels, and relative rigidifications are unchanged.
The depth and wild-image conclusions follow from the same
projection calculation. This proves every assertion.
\end{proof}

\subsection{Termination, reconstruction, and the scope of exhaustion}
\label{sc:proved-exhaustion}

\begin{corollary}[A terminating reverse construction]
\label{sc:terminating-reverse}\label{thm:large-p-exhaustion}\label{prop:reconstruct-galois-yu-datum}
Under the assumptions of Theorem~\ref{sc:reverse-step-theorem},
iteration terminates at a depth-zero cuspidal rigid seed. It
produces an elliptic tame twisted Levi tower, full generic quotient
characters at strictly increasing positive depths when read
forwards, compatible tame \(L\)-embeddings, and relative rigid
records using the original common \(Z\). Forward reconstruction
recovers the prescribed enhanced parameter and its ambient rigid
class. Thus every cuspidal rigidly enhanced parameter in this prime
range is obtained from a quotient-character datum of the specified
kind, allowing the final central twist already present in
Definition~\ref{sc:def-data}.
\end{corollary}
\begin{proof}
Let \(B_0\) be the initial finite wild image. At every positive
step, \eqref{sc:eq-residual-wild-map} gives a proper quotient
of the preceding wild image. Its order decreases by at least a
factor of \(p\), so at most \(\log_p|B_0|\) positive steps
are possible. Equivalently, the kernels on \(P_F\) strictly
increase and every later wild homomorphism factors through the
initial finite quotient. There is no new ramification-break
assumption in this termination proof.

After each step the remaining wild image lies in the derived
group of the new dual group. The intersection of that derived
group with its center has order prime to \(p\). Therefore a
later nontrivial top-depth image cannot be central in the same
group. A central step can occur only as the first reverse step;
it becomes the optional final central twist on reading the
construction forwards. Every other step has a proper Levi and
strictly decreases depth. Consecutive relative realizations
compose by Proposition~\ref{sc:relative-composition}, giving
an actual nested rational tower in the original ambient form.

At termination the wild restriction is trivial. All transported
enhancements remain cuspidal and relevant, so the endpoint is an
intrinsic depth-zero rigid seed. Reversing each identity
\eqref{sc:eq-proved-reverse-step} recovers the parameter, and
reversing each \(\alpha\) recovers its original enhancement.
The retained relative reductions give a compatible marked realization
in the original ambient rigid class. Choose compatible building data
on this tame tower. Genericity at compatible points follows from
Lemma~\ref{sc:graded-genericity}.
\end{proof}

The reverse procedure chooses tame splittings, cocycle extensions and
rational references. Two extensions of a fixed central wild cocycle
differ by a tame cocycle. Section~\ref{sc:equivalence} compares the
resulting changes when the rational identifications are supplied.
Theorem~\ref{sc:exact-fibers} compares arbitrary reverse extractions
through their entire wild restriction and reduced enhanced seed.

\begin{corollary}[Construction and exhaustion for quotient-character data]
\label{sc:proved-construction-exhaustion}
Let \(G\) be tame and let \(p\nmid|W_{\mathrm{abs}}(G)|\).
Fix an ambient rigid inner form and a finite central subgroup \(Z\)
representing its rigidification \(\zeta_G\). Write \(\Dpres\) for the
strict rationally marked presentations with a relevant depth-zero
cuspidal rigid seed and compatible relative records
(Definition~\ref{sc:classification-data}). There is a surjection
\[
 \mathrm{Fwd}:\Dpres\twoheadrightarrow
       \Phi_{\mathrm{cusp}}^{+,Z}(G,\zeta_G).
\]
Each proper stage has the properties of Theorem~\ref{sc:forward-theorem}.
A finite reverse extraction constructs a preimage of any prescribed
enhanced class, with its original ambient relevance character.
\end{corollary}
\begin{proof}
At each proper step, Theorem~\ref{sc:generic-stage-theorem}
supplies the tame embedding, full centralizer equalities,
admissibility, and discreteness. The relative records give the prescribed relevance compatibility.
A final central twist leaves all centralizers and enhancements
unchanged; if it has a new highest depth, its leading restriction
gives that depth. The depth-zero case is the supplied seed.
Theorem~\ref{sc:forward-theorem} therefore applies.
Theorem~\ref{sc:reverse-step-theorem} and
Corollary~\ref{sc:terminating-reverse} construct a preimage in the
references used for extraction. In other coordinates the transports
are \(T_{i+1}E_iT_i^{-1}\), and their composite with the transported
seed gives the same final label since \(T_d=1\).
\end{proof}

\section{Rational presentation moves and the Weyl-group obstruction}
\label{sc:equivalence}

Local changes of characters, seeds and rational references preserve
the enhanced parameter. We prove a converse within an aligned
presentation and exhibit two rationally inequivalent presentations of
one enhanced parameter.
\subsection{The input class and its normalization}
\label{sc:classification-input}

Fix a tame ambient rigid inner twist \(\mathfrak g=(G^*,G,\xi,z)\),
its class \(\zeta_G\), and its common finite central subgroup \(Z\).
Assume \(p\nmid|W_{\mathrm{abs}}(G)|\). In each presentation retain
an explicit quasi-split rational reference for the actual marked tower,
as in Theorem~\ref{sc:simultaneous-references}. In that reference
use the inverse-image maps \(E_i\) of the common cover. A presentation
may equally be expressed in a second reference by retaining the
comparison to the first. Its seed is transported by \(T_0\), its
edge maps are \(E_i'=T_{i+1}E_iT_i^{-1}\), and \(T_d=1\), as in
Theorem~\ref{sc:coherent-origin-transport}. All these data travel
under an actual marked isomorphism. This convention allows a separate
choice for each presentation.

\begin{definition}[Strict rationally marked quotient-character data]
\label{sc:classification-data}
Write \(\Dpres\) for the following
presentations of the data of Definition~\ref{sc:def-data}.
The actual rational tower is elliptic, tame, and strict, with
\(d\) proper inclusions and \(G^d=G\). Its proper-step depths
are \(0<r_0<\cdots<r_{d-1}\). The datum contains the full
quotient characters \(\bar\chi_0,\ldots,\bar\chi_d\), a
depth-zero cuspidal rigid seed, and its compatible marked relative
records. It also supplies compatible tame \(L\)-embeddings,
cocycle representatives, and full dual root-data and cover markings.
The latter must be the markings of the actual rational tower, as
in Proposition~\ref{sc:reference-root-data}. Use the enhancement
transports of its stated reference, with the coherent rule above for
any subsequent change of coordinates.

The final character \(\bar\chi_d\) is an arbitrary full central
factor, with no restriction on its depth. In particular, when
\(d>0\) it may have depth at most \(r_{d-1}\).
For \(d=0\), the datum is a depth-zero rigid seed followed by
one such quotient-character twist on \(G\). Set
\[
 b_0=0,\qquad b_i=r_{i-1}\quad(1\le i\le d),
 \qquad A_i=Z(\widehat{G^i})^\circ,\quad
 C_i=Z(\widehat{G^i}).
\]
Thus the parameter \(\phi_i\) before applying the \(i\)-th
character has depth at most \(b_i\), and has depth \(b_i\)
for \(i>0\).

Tori, splitting fields, and compatible building points used to witness
genericity are auxiliary choices. The intrinsic depth-zero seed carries
no additional torus or vertex label.
\end{definition}

The relative realization is part of a presentation. In particular the
distinct source relevance sectors in Appendix~\ref{sc:rank-one} are
recorded before taking any quotient.

\subsection{Elementary rational presentation moves}
\label{sc:exact-moves}

\begin{definition}[Rational equivalence of marked presentations]
\label{sc:exact-equivalence}
On \(\Dpres\), let
\(\sim_{\mathrm{Yu}}^{+,Z}\) be the equivalence relation generated
by the following moves and their inverses. Every move carries its
full markings and relative records as indicated.
\begin{enumerate}[label=\textup{(E\arabic*)},leftmargin=3em]
\item \emph{Isomorphisms and representatives.}
Allow isomorphisms of enhanced seeds, isomorphisms of the relative
reductions, and simultaneous rational conjugacy in the fixed ambient
form. Allow actual reference-group isomorphisms with their duals,
characters, covers, and rigid records transported as in Proposition
\ref{sc:reference-isomorphism-covariance}; the ambient marking
and common \(Z\) are fixed. Also allow simultaneous changes of
dual coordinates and cocycle representatives with the prescribed
lifted dual conjugacies. Replace auxiliary genericity witnesses
while keeping the preceding substantive data fixed. A source outer
isomorphism carries its enhancement by pullback; it does not simply
erase the source marking.

\item \emph{A depth-zero quotient twist at the seed.}
For a character \(\bar\eta:D_0(F)\to\C^\times\) of depth
at most zero, replace
\[
 (\Phi_0^+,\bar\chi_0)
       \longmapsto
 (T_{\mathcal L_{D_0}(\bar\eta)}\Phi_0^+,
                         \bar\chi_0\bar\eta^{-1}).
\]
Here \(T\) is the intrinsic central twist, with unchanged
enhancement class and relative record. The operation also applies
when \(d=0\).

\item \emph{A neighboring-character transfer.}
For \(1\le i\le d\), take a full character
\(\bar\eta:D_i(F)\to\C^\times\) of depth at most
\(b_i=r_{i-1}\). Replace just
\begin{equation}\label{sc:eq-neighbor-move}
 (\bar\chi_{i-1},\bar\chi_i)
       \longmapsto
 (\bar\chi_{i-1} f_{i-1}^*\bar\eta,
                         \bar\chi_i\bar\eta^{-1}).
\end{equation}
Keep the seed enhancement, all embeddings, and all relative
records fixed. Use the functorial cocycle representatives, up to
the representative changes in \textup{(E1)}.

\item \emph{Compensated tame embedding change.}
Replace the Weil parts by \(\iota_i'=\delta_i\iota_i\),
where \(\delta_i\in Z^1(W_F,C_i)\) is tame and the
connected-dual inclusions are unchanged. Keep all quotient characters
and the final factor fixed. Set
\[
 v_d=1,\qquad v_i=\delta_i^{-1}\nu_i(v_{i+1}),
 \qquad \Phi_0^{+\prime}=T_{v_0}\Phi_0^+.
\]
Here \(\nu_i:C_{i+1}\to C_i\) is the full-center inclusion.
The enhancement and relative records are unchanged. These are
full-center twists; they need not be quotient-torus characters.

\item \emph{Calibrated change of rational reference.}
After actual reference isomorphisms, write
\(j_i'=\Ad(b_i^{\mathrm{ref}})j_i\), with the full marked
root-data diagrams identified. Use exactly the relative-record
change and the characters \(T_i\) of Theorem
\ref{sc:coherent-origin-transport}: replace the seed by
\(T_0\Phi_0^+\) and the edge transports by
\(E_i'=T_{i+1}E_iT_i^{-1}\). The ambient correction is
\(T_d=1\). Any accompanying tame change of Weil parts is
compensated by \textup{(E4)}. No arbitrary character with the
same restriction to the relevance center may replace \(T_i\).
\end{enumerate}
For \textup{(E1)} and \textup{(E5)}, a rational embedding record
may replace a rigid reduction only through the mutually inverse
constructions of Proposition~\ref{sc:marked-extraction}.
\end{definition}

This list does not permit twisting a seed while holding its
compensating factors fixed, changing an embedding without its
compensation, or identifying different relative reductions merely
because their underlying source groups are isomorphic.

\begin{lemma}[The generators are valid and preserve the output]
\label{sc:generators-valid}
Each move in Definition~\ref{sc:exact-equivalence} is a
move between valid data and preserves the enhanced output, up to
the specified ambient dual conjugacy. It preserves the prescribed
ambient rigid class and the transported marked realization.
\end{lemma}
\begin{proof}
The assertions for \textup{(E1)} are covariance and representative
independence from Theorem~\ref{sc:rigid-preservation}, Lemma
\ref{sc:coboundary}, and Proposition
\ref{sc:reference-isomorphism-covariance}. The recursive
construction does not use a choice of genericity witness once the
full quotient character and marked tower are fixed. In the present
prime range Lemma~\ref{sc:graded-genericity} identifies the
genericity condition at compatible points.

For \textup{(E2)}, the seed twist is tame, and the corresponding
factor in \(\kappa_0\) cancels it. Central twisting leaves
the enhancement groups, cuspidality, and relevance unchanged.
For \textup{(E3)}, let \(a\) represent
\(\mathcal L_{D_i}(\bar\eta)\). The new parameter at
level \(i\) is \(a\phi_i\); the factor \(a^{-1}\)
in the next character cancels it. All later parameters are
unchanged, including the final one when \(i=d\). The twist
at level \(i\) is central there, so its refined group and
enhancement class are unchanged.

We also check validity of this neighboring move. The character
\(\bar\chi_{i-1}f_{i-1}^*\bar\eta\) still has depth
\(r_{i-1}\): the added factor has depth at most that number,
and its coroot restrictions for roots of \(G^i\) are trivial.
It therefore cannot cancel any of the generic leading root
pairings of the proper step. If \(i<d\), the change to
\(\bar\chi_i\) has depth at most \(r_{i-1}<r_i\),
so its depth and generic leading term are unchanged. There is no
genericity requirement on the final ambient-central character.
Admissibility, discreteness, and the unchanged relevance sectors
then follow from the stage theorem and central twisting. The
relative reductions are the same.

Theorem~\ref{sc:tame-embedding-compensation} proves all
assertions for \textup{(E4)}. Theorem
\ref{sc:coherent-origin-transport} proves them for
\textup{(E5)}, including the nontrivial enhancement correction.
The latter is essential: using a raw edge map in each rational
origin would not prove this lemma.
\end{proof}

\subsection{Refactorization bounds and neighboring transfers}
\label{sc:refactorization-bounds}

\begin{proposition}[An exact fixed-presentation refactorization criterion]
\label{sc:exact-refactorization}
Fix a strict tower, its depths, \(L\)-embeddings, rational
reference presentation, and relative reductions. For two tuples
of full quotient characters define
\begin{equation}\label{sc:eq-classification-xi}
 \bar\xi_i=\prod_{j=i}^{d}
           f_{i,j}^*(\bar\chi_j\bar\chi_j'^{-1}).
\end{equation}
The following condition is equivalent to passage between the data
by moves \textup{(E2)}--\textup{(E3)}, up to representatives:
\begin{equation}\label{sc:eq-refactorization-bounds}
 \depth(\bar\xi_i)\le b_i\quad(0\le i\le d),
 \qquad
 \Phi_0^{+\prime}\simeq
 T_{\mathcal L_{D_0}(\bar\xi_0)}\Phi_0^+.
\end{equation}
The seed equality includes its enhancement and relative record.
In particular, the full quotient characters, rather than only
their pullbacks to \(G^i(F)\), occur in this criterion.
\end{proposition}
\begin{proof}
The cumulative characters satisfy
\[
 \bar\chi_i'=\bar\chi_i\bar\xi_i^{-1}f_i^*\bar\xi_{i+1}
       \quad(i<d),\qquad
 \bar\chi_d'=\bar\chi_d\bar\xi_d^{-1}.
\]
First apply \textup{(E2)} with \(\bar\eta=\bar\xi_0\),
then \textup{(E3)} with \(\bar\eta=\bar\xi_i\) in the
order \(i=1,\ldots,d\). The displayed identities show
that the result is the primed tuple and seed. The preceding lemma
shows that every intermediate datum is valid, including when a
temporary final factor has small depth.

Conversely, one neighboring transfer changes only the \(i\)-th
cumulative character in \eqref{sc:eq-classification-xi},
by \(\bar\eta\); the seed move changes only the zeroth
one. Characters of depth at most \(b_i\) form a subgroup,
and their seed twists multiply. Thus any sequence gives exactly
\eqref{sc:eq-refactorization-bounds}. Coboundaries are
handled by \textup{(E1)}.
\end{proof}

\begin{lemma}[Automatic bounds and separation of embedding changes]
\label{sc:automatic-refactorization-bounds}
Use the identified rational presentation, connected inclusions, and
bottom-dual conjugator of Theorem~\ref{sc:realized-comparison-generation},
and assume that its lift carries the endpoint label. Then a compensated
tame embedding change followed by the fixed-presentation refactorization
of Proposition~\ref{sc:exact-refactorization} relates the two data.
The cumulative characters satisfy \eqref{sc:eq-refactorization-bounds}.
\end{lemma}
\begin{proof}
First use \textup{(E4)} to change the embeddings of the first
datum to those of the second. This produces a valid intermediate
datum with the first character tuple, a possibly full-center tame
twist of its seed, and the same enhanced output in the given coordinates.
With the embeddings now equal, absorb the stipulated bottom-dual
conjugator into the first seed and its label as in the proof of
Theorem~\ref{sc:realized-comparison-generation}. The two full endpoint
maps then agree pointwise.

Writing \(\phi_i\) and \(\phi_i'\) for the intermediate
parameters in this common presentation, solve backwards:
\begin{equation}\label{sc:eq-connected-backward-comparison}
 v_d=\kappa_d\kappa_d'^{-1},\qquad
 v_i=\kappa_i'^{-1}\kappa_i\nu_i(v_{i+1}),\qquad
 \phi_i'=v_i\phi_i.
\end{equation}
The last equality follows inductively from injectivity of the
\(L\)-embeddings. Every \(v_i\) now takes values in
\(A_i\): the terminal ratio does, all \(\kappa_i\)'s
do, and \(\nu_i(A_{i+1})\subset A_i\). Its torus-LLC
class is exactly \(\mathcal L_{D_i}(\bar\xi_i)\).
Since both \(\phi_i\) and \(\phi_i'\) have trivial dual parts
on \(I_F^{b_i+}\), their ratio is trivial there.
For \(i=0\) this means triviality on \(P_F\). Tame torus depth preservation
therefore gives \(\depth(\bar\xi_i)\le b_i\).

The common-cover enhancement maps are isomorphisms. Central
twisting leaves them unchanged, so equality of the final enhanced
labels determines the required equality of bottom labels. Their
central characters determine the same bottom rigid class; Theorem
\ref{sc:one-step-realization} for the composite reference
embedding identifies the relative records modulo the isomorphisms
in \textup{(E1)}. This proves the seed clause of
\eqref{sc:eq-refactorization-bounds} and hence the assertion.

The order of the argument matters. Before removing embedding
differences, the recursion contains \(\delta_i^{-1}\) in
the full center \(C_i\); it need not be a cocycle in
\(A_i\), and it need not correspond to a quotient character.
The tame embedding move removes precisely this discrepancy.
\end{proof}

\subsection{The terminal factor and reduced conventions}
\label{sc:terminal-normalization}

Call a datum \emph{terminally normalized} if
\(\bar\chi_d=1\) or \(\depth(\bar\chi_d)>b_d\).
For \(d=0\), this permits a positive-depth final twist or
no final twist.

\begin{proposition}[Terminal normalization]
\label{sc:terminal-normalization-prop}\label{def:reduced-howe}
Every datum is equivalent to a terminally normalized one without
changing its proper tower or depths. For two terminally normalized
data with the same enhanced output, the final character is trivial
in one datum if and only if it is trivial in the other. Under the
additional alignment hypotheses of Lemma~\ref{sc:automatic-refactorization-bounds},
their remaining character comparison, after compensation of the
embeddings, satisfies \eqref{sc:eq-refactorization-bounds}.
\end{proposition}
\begin{proof}
If \(d>0\) and \(\depth(\bar\chi_d)\le b_d\), use
\textup{(E3)} at \(i=d\) with
\(\bar\eta=\bar\chi_d\). This absorbs the terminal
factor into \(\bar\chi_{d-1}\) and makes the new final
factor trivial. Its genericity is preserved by Lemma
\ref{sc:generators-valid}. If \(d=0\), use
\textup{(E2)} to absorb a depth-zero final factor into the seed.

The wild centralizer flag recovers the proper transition depths.
For a terminally normalized datum the final parameter has depth
\(b_d\) if the terminal factor is trivial, and strictly larger
depth if it is nontrivial: the parameter before that factor has
depth at most \(b_d\), so it cannot cancel a higher-depth
central factor. The same assertion holds with \(b_d=0\)
when \(d=0\). Equality of the outputs therefore determines
the trivial/nontrivial alternative. The last assertion is Lemma
\ref{sc:automatic-refactorization-bounds}.
\end{proof}

The conditions now have the same form as the refactorization
conditions F0--F2 recalled in \cite[\S3.5]{Kaletha}: the terminal
alternative agrees, cumulative factors have the preceding-depth
bound, and the bottom object is twisted by the cumulative factor.
Here the bottom object is an intrinsic enhanced parameter and the
characters are specified on full quotient tori. Pulling the
characters back to the groups gives the corresponding character
identities and depth bounds, by Lemma
\ref{sc:graded-genericity}. This is a precise comparison of
the formulas, not an invocation of representation-side classification.
The calibrated reference and tame embedding moves are additional
parameter-presentation operations whose compensation is explicit.

\subsection{A rational obstruction to the converse}

\begin{proposition}[A ramified torus with trivial rational Weyl group]
\label{sc:rational-weyl-obstruction}
For \(F=\mathbb Q_7\), split \(G=\SL_2\), and \(Z=1\), two strict
one-step presentations in the same rational reference have conjugate
enhanced parameters and are inequivalent under moves \textup{(E1)--(E5)}
of Definition~\ref{sc:exact-equivalence}.
\end{proposition}
\begin{proof}
Put \(E=F(\varpi_E)\), \(\varpi_E^2=7\), and
\(T=\operatorname{Res}^1_{E/F}\mathbb G_m\). Let \(j:T\hookrightarrow G\)
be the multiplication embedding on the \(F\)-space \(E\). Use \(j\) as reference
and take the neutral source and ambient rigidifications. The map
\(t\mapsto(t-1)/\varpi_E\bmod\varpi_E\) identifies
\(T(F)_{1/2}/T(F)_{1/2+}\) with \((\mathbb F_7,+)\).
Surjectivity follows by taking
\(t=(1+a\varpi_E)/(1-a\varpi_E)\) for \(a\in\mathbb Z_7\), whose residue
in this quotient is \(2a\). Extend a nontrivial character of this quotient to
a full character \(\theta\) of \(T(F)\), using divisibility of \(\C^\times\)
on the discrete quotient by \(T(F)_{1/2+}\). Its leading character has
order \(7\). The normed coroot on \(E^\times\) is \(x\mapsto x/\sigma(x)\)
up to inversion, with \(\sigma\) the nonidentity element of \(\mathrm{Gal}(E/F)\).
It surjects on the displayed quotient by the same calculation.
Hence \(\theta\) and \(\theta^{-1}\) are both generic.

Take the torus parameter \(w\mapsto(1,w)\) with trivial enhancement as seed in
both presentations. Use \(\theta\) and \(\theta^{-1}\) on the same embedded
torus, with trivial terminal factor. Let \(\epsilon:W_F\to\mathbb Z/2\) cut out \(E/F\) and let
\(n\in\PGL_2(\C)\) be an order-two Weyl representative acting on
\(\widehat T\) by inversion. Use the tame embedding
\((t,w)\mapsto(t n^{\epsilon(w)},w)\) in both presentations. If \(\kappa\)
represents \(\theta\), the two Weil parameters have dual-group parts
\[
       \kappa(w)n^{\epsilon(w)},\qquad
       \kappa(w)^{-1}n^{\epsilon(w)} .
\]
They are conjugate by \(n\). In either case the nontrivial leading
\(7\)-group cuts the ordinary centralizer down to the torus before
imposing the remaining Weil equations. An element above the nontrivial
element of \(\mathrm{Gal}(E/F)\) then gives \(t=t^{-1}\). Thus the full
centralizer is \(\mu_2\), the \(\SL_2\)-factor is trivial, and the transported
trivial enhancement is cuspidal and fixed by \(n\).

On the rational side any \(F\)-linear implementer of inversion of \(T\) has
the form \(x\mapsto a\sigma(x)\), \(a\in E^\times\). Indeed, composing
with \(\sigma\) gives a map commuting with multiplication by \(T\), which
generates \(E\) as an \(F\)-algebra. The determinant of \(x\mapsto a\sigma(x)\)
is \(-N_{E/F}(a)\). Since a norm of valuation zero comes from a unit and
reduces to a square modulo \(7\), it cannot equal \(-1\).
Hence \(N_G(jT)(F)/jT(F)=1\).

For this tower the rational conjugacy class of the embedded torus
together with its leading character is invariant under the five moves.
In (E1) reparametrizing the source transports both the character
and its embedding. On the image the character stays the same up to
ambient rational conjugacy. The seed twist (E2) is trivial on the
positive-depth quotient. In (E3) the prospective transfer comes from
\(D_G=1\). The embedding compensation (E4) keeps the quotient
character fixed. The origin change (E5) keeps the realized embedding
and its character fixed. The two leading characters are distinct
inverses and the rational Weyl group is trivial.
\end{proof}

In this example the dual flag is the same for both presentations and
the dual Weyl element acts on its bottom term by inversion. Its action
on the embedded character cannot be realized by a rational conjugacy.
This is the obstruction to deducing rational refactorization from
recovery of the dual flag, even when the actual reference embedding
has already been fixed. Appendix~\ref{sc:rank-one} concerns the
separate question of preservation of intrinsic AMS groups.

\section{Classification by wild restriction and depth-zero enhancement}
\label{sc:wild-classification}

In the ramified torus example of
Proposition~\ref{sc:rational-weyl-obstruction}, two presentations
give conjugate enhanced parameters, although the Weyl element
relating them has no representative in the rational normalizer of
the marked torus. To decide whether two presentations give the same
enhanced class, we first compare their restrictions to wild inertia.
Once both restrictions have been identified with the same homomorphism
\(\omega:P_F\to\hG\), the remaining conjugating elements lie in
\(L_\omega=Z_{\hG}(\omega(P_F))\). Fixing a Weil section carrying
\(\omega\) then gives a reduced depth-zero parameter for the group
\(H_\omega\) with dual group \(L_\omega\). We compare these reduced
parameters, with their individual enhancements, up to
\(L_\omega\)-conjugacy. The reduction puts the two original seeds
in a common group where this comparison makes sense.

\subsection{The wild restriction and its depth-zero parameters}

For underlying Weil parameters the reduction is due to
Dat--Helm--Kurinczuk--Moss \cite[\S3]{DHKM} and Dat--Fintzen
\cite[\S4.3.1, diagram (4.3.2)]{DatFintzen}. We describe it for
cuspidal rigid enhancements relevant to the prescribed ambient form.

Let \(\omega:P_F\to\hG\) be continuous with finite image and assume
\(u\mapsto(\omega(u),u)\) extends to an admissible Weil parameter of
the quasi-split group \(G^*\). Write \(\widetilde\omega(u)=(\omega(u),u)\)
and let \(\operatorname{pr}:{}^LG\to W_F\) be the projection. Define
\begin{align}
 L_\omega&=Z_{\hG}(\omega(P_F)),\notag\\
 \mathcal L_\omega
 &=\bigl\{x\in{}^LG \mid
     \Ad(x)\circ\widetilde\omega
        =\widetilde\omega\circ\Ad(\operatorname{pr}(x))|_{P_F}\bigr\}.
        \label{sc:eq-wild-transporter}
\end{align}
An extension of \(\widetilde\omega\) gives surjectivity of
\(\mathcal L_\omega\to W_F\). Its kernel is \(L_\omega\), since \(P_F\)
acts trivially on \(\hG\). Lemma~\ref{sc:p-centralizers} makes
\(L_\omega\) a connected Levi. Conjugation by \(\mathcal L_\omega\)
gives an outer Weil action on \(L_\omega\) and an action on its center.
The outer action is tame since \(\omega(P_F)\subset Z(L_\omega)\).
Write \(H_\omega\) for the corresponding quasi-split group with
\(\widehat H_\omega=L_\omega\). Call \([\omega]\)
\emph{elliptic} when
\begin{equation}\label{sc:eq-wild-ellipticity}
       (Z(L_\omega)^{W_F})^\circ=(Z(\hG)^{W_F})^\circ .
\end{equation}
Let \(\Omega_{\rm ell}\) be the set of \(\hG\)-conjugacy classes of
extendible wild homomorphisms satisfying \eqref{sc:eq-wild-ellipticity}.

For each such class choose a representative, a pinning of \(L_\omega\),
and a pinning-preserving section \(\sigma_\omega:W_F\to\mathcal L_\omega\)
with \(\sigma_\omega(u)=(\omega(u),u)\). Such sections exist by
\cite[\S2.2.5(i)]{DatFintzen}. Equivalently one may first take the tame
section of \cite[\S2.2.4]{DatFintzen} and multiply it by a central extension of
\(\omega\). Indeed its \(p\)-primary image lies in \(Z(L_\omega)^\circ\),
and Lemma~\ref{sc:wild-central-extension} extends its equivariant
homomorphism on \(P_F\) to a Weil cocycle in that center.
For \(\omega\ne1\), the prescribed wild part makes \(\sigma_\omega\)
different from the tame sections at the individual Yu steps. Choose a
marked rational reference \(j_\omega:H_\omega\hookrightarrow G^*\)
realizing the inclusion \(L_\omega\subset\hG\), using the quasi-split
centralizer of a rational Levi-center embedding
\cite[Proposition~2.1.5]{DatFintzen}, as in
Proposition~\ref{sc:parameter-levi-reference}. Use the original \(Z\)
and its dual inverse-image cover. The section defines the isomorphism over \(W_F\)
\begin{equation}\label{sc:eq-wild-section}
 J_\omega:\LG{H_\omega}=L_\omega\rtimes_{\Ad\sigma_\omega} W_F
       \xrightarrow{\ \sim\ }\mathcal L_\omega\subset\LG G,\qquad
                    (h,w)\longmapsto h\sigma_\omega(w).
\end{equation}

\begin{theorem}[Classification by the whole wild restriction]
\label{sc:wild-bijection}
With these choices there is a bijection
\begin{equation}\label{sc:eq-wild-classification}
 \coprod_{[\omega]\in\Omega_{\rm ell}}\
 \coprod_{\eta\in j_{\omega,*}^{-1}(\zeta_G)}
       \Phi_{0,\mathrm{cusp}}^{+,Z}(H_\omega,\eta)
       \ \xrightarrow{\ \sim\ }\ \Phi_{\mathrm{cusp}}^{+,Z}(G,\zeta_G).
\end{equation}
The notation \(\Phi_{0,\mathrm{cusp}}^{+,Z}(H_\omega,\eta)\) means the
\(L_\omega\)-conjugacy classes of depth-zero cuspidal parameters with their
irreducible common-cover enhancements relevant to the rigid class \(\eta\).
The indexing map \(j_{\omega,*}\) extends rigid torsors along \(j_\omega\),
from \(H^1_{\mathrm{rig}}(F,Z\to H_\omega)\) to
\(H^1_{\mathrm{rig}}(F,Z\to G^*)\). In the summand \([\omega]\), the map is composition
with \(J_\omega\) on the parameter and the common-cover identification
on the enhancement. Its image is precisely the stratum of wild
restriction \([\omega]\).
\end{theorem}

\begin{proof}
Let \(\psi\) have depth zero for \(H_\omega\). Then
\(\Phi=J_\omega\psi\) has wild restriction \(\omega\) in the chosen
coordinates. Conversely, conjugate any parameter with this wild class
to have restriction exactly \(\omega\). Its Weil part satisfies
\eqref{sc:eq-wild-transporter}, and its \(\SL_2\)-part belongs to \(L_\omega\) because
it centralizes wild inertia. Hence \(J_\omega^{-1}\Phi\) exists and
is trivial on wild inertia in its dual part. For admissibility, take
a positive power \(N\) of Frobenius killing both pinned actions.
Writing \(\sigma_\omega(\Fr)^N=(z,\Fr^N)\), we have \(z\in Z(L_\omega)\)
because \(z\) centralizes the Levi. If \(\psi(\Fr)^N=(a,\Fr^N)\), then
\(\Phi(\Fr)^N=(az,\Fr^N)\). Semisimplicity of \(a\) is equivalent to that
of \(az\), since \(z\) is central and semisimple. Continuity, finiteness of
the inertia images, and algebraicity on \(\SL_2\) hold in both
directions. The two operations are inverse on parameter maps with
this fixed wild restriction.

An ambient element centralizing \(\Phi|_{W_F}\) centralizes \(\omega\)
pointwise, so it is in \(L_\omega\). Inside \(L_\omega\), \(J_\omega\) is
the identity on the connected dual group and intertwines the two
Weil actions. Thus the ordinary Weil centralizers of \(\psi\) and \(\Phi\)
are equal, as are their full parameter centralizers. Taking inverse
images in the common cover identifies both groups
\(Q^{+,Z}\), their lifted unipotent elements and the refined component
groups. The local systems and their cuspidality thus agree by the
argument of Proposition~\ref{sc:rigid-cuspidality-preserved}.
Equation \eqref{sc:eq-wild-ellipticity} makes discreteness equivalent
on the two sides. Conversely, if \(\Phi\) is discrete then
\((Z(L_\omega)^{W_F})^\circ\) centralizes \(\Phi\) and contains
\((Z(\hG)^{W_F})^\circ\), so equality is necessary. Thus every target
occurs among the listed wild classes.

For conjugacy, two maps with the same \emph{literal} wild restriction
can be conjugate in \(\hG\) only by an element \(h\) satisfying
\(h\omega(u)h^{-1}=\omega(u)\) for every \(u\). Thus the remaining
conjugacy group is precisely \(L_\omega\), which is the conjugacy
group already used in \(\Phi_{0,\mathrm{cusp}}^{+,Z}(H_\omega,\eta)\).
Conjugation by a lift of \(h\) acts on the common cover independently
of the lift, since the kernel is central. Hence the inverse descends to
enhanced conjugacy classes, including independence of the conjugator
used to align \(\omega\).

Finally let \(\rho^+\) be the transported irreducible label. Its central
character pulled back to \(\mathcal Z_{H_\omega}\) via the relevance map of
Section~\ref{sc:rigid-comparison} determines exactly one source class
\(\eta\) by rigid duality. The restriction of that character along
\(t_{j_\omega}\) is the ambient character by
Proposition~\ref{sc:rigid-relevance}. Consequently
\(j_{\omega,*}\eta=\zeta_G\) is equivalent to ambient relevance.
The distinct source sectors in \eqref{sc:eq-wild-classification} are disjoint.
Conversely every label in those sectors gives a relevant ambient
label. This proves both directions of the bijection.
\end{proof}

For each prescribed enhanced ambient class in the \([\omega]\)-stratum,
the choices \(j_\omega,\sigma_\omega\) and the representative \(\omega\)
therefore determine a unique reduced enhanced class and its source
relevance sector.
If the section is changed with the same pinning, its ratio with the
old section is a cocycle in \(Z(L_\omega)\) trivial on \(P_F\).
The inverse depth-zero parameter changes by the inverse tame central
twist, while its enhancement group and label stay fixed. A different
initial reference for the rigid forms can change the indexing by
\(\eta\). Formula \eqref{sc:eq-wild-classification} uses the displayed reference
consistently. For \(\omega=1\) one can choose \(H_\omega=G^*\),
\(j_\omega=\mathrm{id}\) and \(\sigma_\omega(w)=(1,w)\). The corresponding
map is the identity on the depth-zero sector.

\subsection{The fibers on strict tower presentations}
\label{sc:fiber-classification}
The character and embedding factors of a strict presentation determine a wild homomorphism
and a change of \(L\)-group coordinates. In those coordinates it suffices
to compare depth-zero enhanced seeds.

For a central cocycle \(\kappa\) at any level write
\(m_\kappa(h,w)=(\kappa(w)h,w)\). This is an automorphism of that
\(L\)-group over \(W_F\). For a presentation \(\mathcal D\) form the
homomorphism, using no seed,
\begin{equation}\label{sc:eq-tower-embedding}
 J_{\mathcal D}=m_{\kappa_d}\circ\iota_{d-1}\circ
   m_{\kappa_{d-1}}\circ\cdots\circ\iota_0\circ m_{\kappa_0}:
               \LG{G^0}\hookrightarrow\LG G .
\end{equation}
For \(d=0\) use \(J_{\mathcal D}=m_{\kappa_0}\). Define
\(\omega_{\mathcal D}(u)\) by
\(J_{\mathcal D}(1,u)=(\omega_{\mathcal D}(u),u)\).
On \(P_F\) this is the product of the \(\kappa_i|_{P_F}\) embedded in the
connected center of the bottom dual Levi, including the terminal factor.

\begin{lemma}\label{sc:wild-bottom}
The connected restriction of \(J_{\mathcal D}\) identifies \(\widehat{G^0}\)
with \(L_{\omega_{\mathcal D}}\), and \(J_{\mathcal D}\) identifies
\(\LG{G^0}\) with \(\mathcal L_{\omega_{\mathcal D}}\).
\end{lemma}
\begin{proof}
The seed is wild-trivial, so the final parameter restricts to
\(\omega_{\mathcal D}\) on \(P_F\). Proposition~\ref{sc:wild-flag-recovery}
at \(s=0\) now gives \(L_{\omega_{\mathcal D}}=\widehat{G^0}\) in the
given coordinates, including the case \(d=0\). Tameness of \(G^0\) gives
\((h,w)(1,u)(h,w)^{-1}=(1,wuw^{-1})\) in \(\LG{G^0}\).
Applying \(J_{\mathcal D}\) puts its image in
\(\mathcal L_{\omega_{\mathcal D}}\). The image covers \(W_F\) and
already contains the kernel \(L_{\omega_{\mathcal D}}\), proving equality.
The wild class is elliptic because the constructed parameter is discrete,
by Theorem~\ref{sc:forward-theorem} and the necessity of
\eqref{sc:eq-wild-ellipticity} proved above.
\end{proof}

Fix the choices of Theorem~\ref{sc:wild-bijection}. For each \(\mathcal D\) choose
\(g_{\mathcal D}\in\hG\) carrying \(\omega_{\mathcal D}\) to the representative
\(\omega\) of its class. Define the isomorphism
\begin{equation}\label{sc:eq-align-seed}
 B_{\mathcal D}=J_\omega^{-1}\circ\Ad(g_{\mathcal D})\circ
     J_{\mathcal D}:\LG{G^0}\xrightarrow{\ \sim\ }\LG{H_\omega}.
\end{equation}
It sends \((1,u)\) to \((1,u)\) for every \(u\in P_F\). Its connected
restriction lifts to the common covers by the supplied inclusion
into the ambient cover and conjugation by a lift of \(g_{\mathcal D}\) in \(\widehat{G/Z}\).
Thus it takes the depth-zero seed, including its irreducible label,
to a depth-zero enhanced class \(R_\omega(\mathcal D)\) for \(H_\omega\).
Use the seed in the reference with the inverse-image transports of
Section~\ref{sc:classification-input}. Equivalently, for
transported edge maps, undo \(T_0\) on the re-expressed seed before
using the common-cover map. If \(g_{\mathcal D}\) is replaced by
\(h g_{\mathcal D}\) then \(h\in L_\omega\), and the new \(B_{\mathcal D}\)
differs by the inner conjugacy \(\Ad(h)\). Thus \(R_\omega(\mathcal D)\)
is well defined as a class. Its character on \(\mathcal Z_{H_\omega}\) supplies
the \(\eta\) indexing the inner summand of \eqref{sc:eq-wild-classification}.

\begin{theorem}[Exact fibers on strict tower presentations]
\label{sc:exact-fibers}\label{prop:howe-independence}
For \(\mathcal D,\mathcal D'\in\Dpres\), write
\(\mathcal D\sim_{\mathrm{sp}}\mathcal D'\) when both conditions hold.
\begin{enumerate}[label=\textup{(\roman*)},leftmargin=2.5em]
\item The two wild homomorphisms in \eqref{sc:eq-tower-embedding} are
\(\hG\)-conjugate.
\item For their common representative \(\omega\), the two reduced
depth-zero parameters \(R_\omega(\mathcal D),R_\omega(\mathcal D')\)
are \(L_\omega\)-conjugate as enhanced parameters. A lift of the
conjugating element must carry one full common-cover label to the other.
\end{enumerate}
These conditions hold if and only if the forward enhanced outputs are
conjugate. Their conjunction is independent of the common choices of
representative, pinning, section and indexing reference. The rational
equivalence of Definition~\ref{sc:exact-equivalence} implies
\(\sim_{\mathrm{sp}}\) and can be strictly finer.
\end{theorem}
\begin{proof}
The definition of \(B_{\mathcal D}\) gives
\[
 J_\omega B_{\mathcal D}=\Ad(g_{\mathcal D})J_{\mathcal D}.
\]
The same identity holds on common-cover transports by restriction of
the lifted inclusions and conjugation. In re-expressed coordinates
the intervening \(T_i\)'s cancel and \(T_d=1\). Thus \(J_\omega\)
takes \(R_\omega(\mathcal D)\), with its individual enhancement, to
the enhanced conjugacy class of \(\mathrm{Fwd}(\mathcal D)\).
Distinct wild classes have disjoint images. For a common wild class,
injectivity in Theorem~\ref{sc:wild-bijection} gives necessity of \textup{(ii)},
and the displayed identity gives sufficiency. This reasoning works
with every common choice for the stratum, proving the independence
assertion. Lemma~\ref{sc:generators-valid} gives the implication from
rational equivalence. Proposition~\ref{sc:rational-weyl-obstruction}
proves that the converse can fail.
\end{proof}

For example, in Proposition~\ref{sc:rational-weyl-obstruction} take
\(\omega=\kappa|_{P_F}\), take \(\sigma_\omega(w)=(\kappa(w)n^{\epsilon(w)},w)\),
and align the inverse wild map by \(g=n\). Equation \eqref{sc:eq-align-seed} gives
the identity for \(\theta\) and torus inversion on the connected dual
group for \(\theta^{-1}\), with identity on the chosen Weil section.
Both send the trivial enhanced seed to the same class. At a fixed full
depth-zero parameter, two inequivalent irreducibles of its refined
component group give distinct enhanced classes. Indeed, a conjugator
fixing the parameter belongs to its centralizer and acts on the refined
group by an inner automorphism after lifting. This is why the second
test compares the individual labels.

\begin{corollary}[Classification of strict presentations by their reduced enhanced seeds]
\label{sc:classification-bijection}
Forward construction induces a bijection
\begin{equation}\label{sc:eq-classification-bijection}
 \Dpres/\!\sim_{\mathrm{sp}}\ \xrightarrow{\ \sim\ }\
        \Phi_{\mathrm{cusp}}^{+,Z}(G,\zeta_G).
\end{equation}
Any terminating reverse extraction gives its inverse class. The
bijection also holds on terminally normalized presentations with the
restriction of \(\sim_{\mathrm{sp}}\).
\end{corollary}
\begin{proof}
Corollary~\ref{sc:proved-construction-exhaustion} gives surjectivity and
Theorem~\ref{sc:exact-fibers} identifies the fibers. Starting with an
enhanced parameter, forward construction on any extracted tower recovers
that parameter and its original label by Corollary~\ref{sc:terminating-reverse}.
Starting with a tower, any extraction of its output has the same wild
class and reduced enhanced seed, by Theorem~\ref{sc:exact-fibers}.
Finally, Proposition~\ref{sc:terminal-normalization-prop} supplies a
terminally normalized representative using rational moves, which preserve
both invariants.
\end{proof}
\appendix
\section{A rank-one obstruction to AMS group preservation}
\label{sc:rank-one}\label{sec:ams}

For a parameter of \(J\), consider the intrinsic AMS convention
\begin{equation}\label{sc:eq-ams}
 \cS_\phi=\pi_0p_J^{-1}
 \bigl(Z_{\widehat J}(\phi)Z(\widehat J)/Z(\widehat J)\bigr),
 \qquad p_J:\widehat J_{\mathrm{sc}}\to\widehat J_{\mathrm{ad}}.
\end{equation}
This is the definition in \cite[\S6, (83)--(86)]{AMS}.
The two groups \(H\) and \(G\) use different adjoint quotients and
different simply connected covers.

\begin{proposition}[Failure for an elliptic torus in \(\SL_2\)]
\label{sc:sl2-counterexample}
Let \(p>2\), let \(E/F\) be unramified quadratic, set
\(G=\SL_2\), and let \(H=\operatorname{Res}^{1}_{E/F}\mathbb G_m\)
be an elliptic maximal torus in \(G\). There is a tame positive-depth
torus step with an intrinsic depth-zero seed for \(H\), full
centralizer cut-down, and discrete output, for which
\[
 \cS_{\phi_H}=1,\qquad \cS_{\phi_G}\simeq\mu_4.
\]
For the split group \(G\), the output has two relevant cuspidal
enhancements, while the fixed torus seed has only one enhancement.
\end{proposition}
\begin{proof}
We have \(\widehat H=\C^\times\) with Frobenius acting by
inversion, and \(\widehat G=\PGL_2(\C)\). Identify the torus
with
\[
 T=\{[\operatorname{diag}(t,1)]:t\in\C^\times\}
       \subset\PGL_2(\C).
\]
Let \(n\) be the class of
\(\left(\begin{smallmatrix}0&1\\-1&0\end{smallmatrix}\right)\).
It has order two in \(\PGL_2\) and conjugates \(t\) to \(t^{-1}\).
If \(\nu:W_F\to\Z\) sends \(\Fr\) to \(1\), then
\[
 \iota(t,w)=(t n^{\nu(w)},w)
\]
defines a tame \(L\)-embedding \({}^LH\to{}^LG\).

Take the depth-zero torus parameter \(\phi_H(w)=(1,w)\), with
trivial \(\SL_2\)-part. Since \(\widehat H_{\mathrm{sc}}=1\),
its enhancement group is trivial, and its unique enhancement is cuspidal.
Choose a positive-depth character of \(H(F)\) with nontrivial leading
restriction and let \(\kappa\) be its full torus-LLC cocycle.
Such characters exist, for example from a nontrivial character of a
positive norm-one unit quotient. Its leading image \(K\subset T\)
is a nontrivial finite \(p\)-group, and hence contains an element
of order greater than two.

The centralizer in \(\PGL_2\) of that element is \(T\).
The nontrivial Weyl element is inversion and does not fix \(K\)
pointwise. Thus the full cut-down and the absolute pointwise stabilizer
condition both hold. The normed coroot on units is, in a splitting
coordinate, \(a\mapsto a/\sigma(a)\). On a positive graded unit
quotient it is \(b\mapsto b-\sigma(b)\), surjective onto the
trace-zero line because \(p\ne2\). The chosen leading character
is therefore nontrivial on the normed coroot at its exact depth, giving
the rank-one rootwise genericity condition. Inversion does not fix its
nontrivial leading \(p\)-character, giving the stabilizer condition too.

Define
\[
 \phi_G(w)=(\kappa(w)n^{\nu(w)},w),
 \qquad \phi_G|_{\SL_2}=1.
\]
This is an admissible parameter: it is a cocycle product, its inertia
image is finite, and its Frobenius dual-group element has square one,
since \(n\) acts on \(T\) by inversion. In fact the dual-group
image is finite. An element centralizing \(\phi_G\) must belong
to \(T\) by cut-down. Commuting with
\(\kappa(\Fr)n\) then imposes \(t=t^{-1}\), so
\[
 C:=Z_{\PGL_2}(\phi_G)
     =\{1,[\operatorname{diag}(-1,1)]\}\simeq\mu_2.
\]
The same group is \(Z_{\widehat H}(\phi_H)\). Since
\(Z(\PGL_2)=1\), the finiteness of \(C\) proves discreteness.

The covering \(\SL_2(\C)\to\PGL_2(\C)\) sends
\(\operatorname{diag}(a,a^{-1})\) to the torus coordinate \(a^2\).
The inverse image of \(C\) is therefore
\[
 \{\operatorname{diag}(a,a^{-1}):a^4=1\}\simeq\mu_4.
\]
It is finite, so \eqref{sc:eq-ams} gives
\(\cS_{\phi_G}=\mu_4\). On the torus side the adjoint quotient
and its simply connected cover are trivial, so
\(\cS_{\phi_H}=1\), as claimed.

The split relevance character on
\(Z(\widehat G_{\mathrm{sc}})=\mu_2\subset\mu_4\) is trivial.
Exactly two characters of \(\mu_4\) restrict trivially to this
\(\mu_2\): the trivial character and its character of order two.
They are cuspidal, since the \(\SL_2\)-part is trivial and the
AMS Weil-centralizer group is finite; its identity component is trivial.
This also follows directly from the observation following
\cite[Definition~6.9]{AMS}.
\end{proof}

This example has full ordinary-centralizer equality, genericity, and
discrete output, but different intrinsic AMS groups. The common-cover
convention of Section~\ref{sc:rigid-comparison} gives the
following description of the same step.

\subsection{The rank-one example with its rigidifications retained}

In Proposition~\ref{sc:sl2-counterexample}, choose
\(Z=Z(\SL_2)=\mu_2\). Then
\(\widehat{G/Z}=\SL_2(\C)\), and the inverse-image torus
\(\widehat{H/Z}=\C^\times\) maps to \(\widehat H\)
by \(a\mapsto a^2\). Since the ordinary centralizer is
\(C=\mu_2\), we obtain
\begin{equation}\label{sc:eq-rigid-rank-one}
 S_{\phi_H}^{+,Z}=\mu_4=S_{\phi_G}^{+,Z},\qquad
 \cS_{\phi_H}^{+,Z}=\mu_4=\cS_{\phi_G}^{+,Z}.
\end{equation}
The common cover thus gives the refined group \(\mu_4\) at both
stages. The intrinsic AMS groups for the same example are \(1\) and
\(\mu_4\), as computed in Proposition~\ref{sc:sl2-counterexample}.

The central relevance groups in this example are
\[
 \mathcal Z_H^{+,Z}=\mu_4,\qquad
 \mathcal Z_G^{+,Z}=\mu_2,
\]
with \(t_{G,H}\) the inclusion. A fixed torus rigidification
prescribes a character of all of \(\mu_4\), selecting exactly
one enhancement. The split ambient rigidification prescribes
triviality only on \(\mu_2\), selecting the two characters
\(\lambda_0(a)=1\) and \(\lambda_2(a)=a^2\).
Rigid duality for the torus identifies these with two distinct
source rigidification classes mapping to the split ambient class.
They are already the two pure classes
\[
 H^1(F,H)=F^\times/N_{E/F}E^\times\simeq\Z/2\Z,
 \qquad H^1(F,\SL_2)=1.
\]
Relative to a reference embedding, these are the two rational
conjugacy classes of embeddings within its stable conjugacy class.
They do not change the abstract torus or its torus-LLC character,
but they must remain distinct in the relative seed record. A
choice of reference embedding fixes the origin of this description.

\subsection{Comparison with rigid cuspidal-support conventions}
\label{sc:ds-interface}
Under the identifications
\[
 (\varphi,u,\rho)_{\mathrm{DS}}=(\phi|_{W_F},\pi_J^Z(u^+),\rho^+),
 \qquad H^+_{\mathrm{DS}}=Q_{J,\phi}^{+,Z},
\]
the enhancement \(\rho^+\) gives the equivariant local system on the
\(H^+_{\mathrm{DS}}\)-orbit of the canonical unipotent lift \(u^+\)
as in \cite{DillerySchwein}. Their comparison with AMS conventions in
\cite[Remark~5.2.3]{DillerySchwein} uses prescribed central-character
and rigidification choices. It does not assert an intrinsic AMS
group isomorphism across a twisted Levi step. Their relevance theorem
for cuspidal support concerns parabolic Levi subgroups; the rational
realization of the twisted Levi steps used here is proved separately
in Section~\ref{sec:relative-realization}.

Enlargement of \(Z\) is compatible with inflation and the comparison
maps of Proposition~\ref{sc:relative-composition} and
Theorem~\ref{sc:coherent-origin-transport}. Descent of an
enhancement through a smaller cover requires its kernel to act
trivially. The classification fixes \(Z\); it does not identify
arbitrary enhancement sectors for different covers.
\section{Quotient-character normalization and label comparisons}\label{sec:quotient-normalization}

Kaletha's quotient-character normalization uses
\cite[Definition~3.7.1 and Lemmas~3.3.2, 3.5.1, and~3.7.2]{Kaletha}.
The argument precedes the regularity restriction on the depth-zero
seed and permits a choice of full characters of the quotient tori
after modifying the remaining datum. The following lemma describes
the ambiguity in extending the positive-depth character.

\begin{lemma}[A local quotient-character normalization]
\label{sc:normalization}
Let \(x\) be a building point for \(J\), and suppose
\[
 q_J(J(F)_{x,t})=D_J(F)_t\quad(t>0)
\]
with kernel \(J_{\mathrm{der}}(F)\cap J(F)_{x,t}\).
Let \(\chi\) be a smooth character of positive depth \(r\) at \(x\)
which is trivial on \(J_{\mathrm{der}}(F)\cap J(F)_{x,0+}\).
Then there is a smooth character \(\bar\chi\) of \(D_J(F)\), of depth
\(r\), such that
\[
 \epsilon=\chi(\bar\chi\circ q_J)^{-1}
 \quad\hbox{is trivial on }J(F)_{x,0+}.
\]
All choices of \(\bar\chi\) agreeing with the prescribed character on
\(D_J(F)_{0+}\) form a torsor under the characters of
\(D_J(F)/D_J(F)_{0+}\).
\end{lemma}
\begin{proof}
The restriction of \(\chi\) descends to a character \(\chi_+\) of
\(D_J(F)_{0+}\). Choose an open subgroup \(U\subset D_J(F)_{0+}\)
on which it is trivial, taking \(U=D_J(F)_{r+}\) if desired.
The character of the subgroup \(D_J(F)_{0+}/U\) extends to the abelian
group \(D_J(F)/U\), since \(\C^\times\) is divisible. Its pullback
is smooth and extends \(\chi_+\). The positive-depth filtration
surjectivity proves that its depth is exactly \(r\). The claimed
property of \(\epsilon\) follows by construction. The quotient of
two such extensions is precisely a character trivial on \(D_J(F)_{0+}\).
\end{proof}

\begin{remark}[The endpoint of a normalization]
At the level of Yu characters, the depth-zero residual factors can be moved
down the tower in the normalization procedure. At the bottom they can alter
the seed. For an intrinsic parameter datum, transporting the seed under an
arbitrary residual character is not supplied by a representation-side
normalization theorem. Section~\ref{sc:seed-twist} treats residual
factors which are given by central cocycles. We therefore make
quotient-character data primary and do not claim to normalize every
more general intrinsic datum without further input.
\end{remark}

\begin{lemma}[A second extension ambiguity]\label{sc:extension-ambiguity}
Suppose a fixed character of \(J(F)\), trivial on
\(J_{\mathrm{der}}(F)\), has extensions \(\bar\chi,\bar\chi'\)
to \(D_J(F)\) through \(q_J\). Their ratio is a character of
\(D_J(F)/q_J(J(F))\). Conversely, multiplication by any smooth character
of this quotient preserves the pullback.
\end{lemma}
\begin{proof}
The two pullbacks agree exactly when their ratio is trivial on the image
of \(q_J\).
\end{proof}

This ambiguity disappears when \(q_J:J(F)\to D_J(F)\) is surjective.
In general it must be recorded or shown to disappear in the final
parameter construction. A difference of quotient-torus classes might
become trivial after a further map into a full dual center; no assertion
that it necessarily survives such a map is made here.

\subsection{Coherence on irreducible labels}
Let \(B\) be a finite group and let \(\alpha,\beta\) be cocycles.
An equality \([\alpha]=[\beta]\) provides existence of a cochain
implementing an isomorphism between twisted group algebra presentations.
If two such cochains differ by \(\eta\in\Hom(B,\C^\times)\),
the induced maps on irreducibles can differ by tensoring with \(\eta\).
For \(B=\Z/2\Z\) and trivial cocycles, the automorphism sending
the nontrivial basis element to its negative exchanges the two simple
modules. Thus a multiplier-class equality does not select an absolute
label.

For enhancement transport defined by projective intertwiners, the
following criterion gives independence of the comparison path on
irreducible labels.

\begin{lemma}[Coherence sufficient on irreducible labels]
\label{sc:inner-coherence}
Suppose comparison maps are algebra isomorphisms and any two comparison
paths between the same objects differ by an inner automorphism of the
endpoint algebra. They induce the same map on isomorphism classes of
irreducible modules.
\end{lemma}
\begin{proof}
If the automorphism is \(a\mapsto uau^{-1}\), then multiplication
by the invertible operator representing \(u\) intertwines the original
module action with its twist. Hence the module isomorphism class is unchanged.
\end{proof}

\section*{Acknowledgments}
I used ChatGPT to explore proof strategies, develop and revise parts of
the arguments, and edit the exposition. I take responsibility for the
statements and proofs in this paper.


\begin{thebibliography}{AMS18}
\bibitem[AMS18]{AMS}
Anne-Marie Aubert, Ahmed Moussaoui, and Maarten Solleveld,
\emph{Generalizations of the Springer correspondence and cuspidal Langlands
parameters}, Manuscripta Math. \textbf{157} (2018), 121--192.
The text incorporating the 2025 erratum is available at
\url{https://www.math.ru.nl/~solleveld/genSpringer.pdf}.

\bibitem[DF26]{DatFintzen}
Jean-Fran\c{c}ois Dat and Jessica Fintzen,
\emph{Parametrization and reduction to depth zero of
\(\overline{\mathbb Z}[1/p]\)-blocks of tame \(p\)-adic groups},
\href{https://arxiv.org/abs/2602.12333v1}{arXiv:2602.12333v1},
12 February 2026, 105 pp. See also
\url{https://webusers.imj-prg.fr/~jean-francois.dat/recherche/publis/tamegroups.pdf}.
See Proposition~2.1.5, \S\S2.2.4--2.2.5, Lemma~2.2.9(iv),
the lemma in \S2.3.1, Lemma~2.4.1, and diagram~(4.3.2).

\bibitem[DHKM25]{DHKM}
Jean-Fran\c{c}ois Dat, David Helm, Robert Kurinczuk, and Gilbert Moss,
\emph{Moduli of Langlands parameters}, J. Eur. Math. Soc.
\textbf{27} (2025), 1827--1927,
\href{https://doi.org/10.4171/JEMS/1599}{doi:10.4171/JEMS/1599}.

\bibitem[Dil23]{Dillery}
Peter Dillery, \emph{Rigid inner forms over local function fields},
Adv. Math. \textbf{430} (2023), Art.~109204,
\href{https://doi.org/10.1016/j.aim.2023.109204}{doi:10.1016/j.aim.2023.109204}.

\bibitem[Dil24]{DilleryLimits}
Peter Dillery, \emph{Isocrystals and limits of rigid local Langlands
correspondences}, to appear in Algebra \& Number Theory.
For Proposition~5.1 and \S5.4, see
\href{https://arxiv.org/abs/2312.09195v2}{arXiv:2312.09195v2} (2024).

\bibitem[DS24]{DillerySchwein}
Peter Dillery and David Schwein,
\emph{A stacky generalized Springer correspondence and rigid
enhancements of \(L\)-parameters}, preprint (2024),
\href{https://arxiv.org/abs/2308.11752v2}{arXiv:2308.11752v2}.

\bibitem[Fin21]{FintzenTypes}
Jessica Fintzen, \emph{Types for tame \(p\)-adic groups},
Ann. of Math. (2) \textbf{193} (2021), no.~1, 303--346,
\href{https://doi.org/10.4007/annals.2021.193.1.4}{doi:10.4007/annals.2021.193.1.4}.

\bibitem[GGP12]{GGP}
Wee Teck Gan, Benedict H. Gross, and Dipendra Prasad,
\emph{Symplectic local root numbers, central critical \(L\)-values,
and restriction problems in the representation theory of classical groups},
in \emph{Sur les conjectures de Gross et Prasad. I},
\href{https://numdam.org/item/AST_2012__346__1_0/}{Ast\'erisque
\textbf{346} (2012), 1--109}.

\bibitem[HM08]{HM}
Jeffrey Hakim and Fiona Murnaghan, \emph{Distinguished tame supercuspidal
representations}, Int. Math. Res. Pap. IMRP \textbf{2008}, no.~2,
Art.~ID rpn005, 166 pp.

\bibitem[Kal16]{KalethaRigid}
Tasho Kaletha, \emph{Rigid inner forms of real and \(p\)-adic groups},
Ann. of Math. (2) \textbf{184} (2016), no.~2, 559--632,
\href{https://doi.org/10.4007/annals.2016.184.2.6}{doi:10.4007/annals.2016.184.2.6}.
The refined centralizer and relevance convention are in \S5.4;
see \url{https://www.math.uni-bonn.de/people/kaletha/ri4.pdf}.

\bibitem[Kal18]{KalethaIsocrystals}
Tasho Kaletha, \emph{Rigid inner forms vs isocrystals},
J. Eur. Math. Soc. \textbf{20} (2018), no.~1, 61--101,
\href{https://doi.org/10.4171/JEMS/759}{doi:10.4171/JEMS/759}.

\bibitem[Kal19]{Kaletha}
Tasho Kaletha, \emph{Regular supercuspidal representations}, J. Amer.
Math. Soc. \textbf{32} (2019), no.~4, 1071--1170,
\href{https://doi.org/10.1090/jams/925}{doi:10.1090/jams/925}.

\bibitem[Kal21]{KalethaPackets}
Tasho Kaletha, \emph{Supercuspidal \(L\)-packets}, to appear in
\href{https://annals.math.princeton.edu/articles/22753}{Ann. of Math.}
The proposition numbering follows
\href{https://arxiv.org/abs/1912.03274v2}{arXiv:1912.03274v2} (2021).

\bibitem[Kal26]{KalethaCovers}
Tasho Kaletha, \emph{On \(L\)-embeddings and double covers of tori
over local fields}, Amer. J. Math. \textbf{148} (2026), no.~5, 1229--1300,
\href{https://doi.org/10.1353/ajm.2026.a1002870}{doi:10.1353/ajm.2026.a1002870}.
References to Lemma~6.4, Fact~6.5, and Corollary~6.7 use the numbering of
\href{https://arxiv.org/abs/1907.05173v2}{arXiv:1907.05173v2}.

\bibitem[KalE]{KalethaErrata}
Tasho Kaletha, \emph{Errata and addenda}, entry for
\emph{Rigid inner forms of real and \(p\)-adic groups},
correction to Corollary~5.4 (irrelevant to the nonarchimedean case),
\url{https://www.math.uni-bonn.de/people/kaletha/errata.pdf}.

\bibitem[Yu09]{YuTori}
Jiu-Kang Yu, \emph{On the local Langlands correspondence for tori},
in \emph{Ottawa Lectures on Admissible Representations of Reductive
\(p\)-adic Groups}, Fields Inst. Monogr. \textbf{26}, Amer. Math.
Soc., Providence, RI, 2009, pp.~177--183,
\href{https://doi.org/10.1090/fim/026/07}{doi:10.1090/fim/026/07}.
See Theorem~7.10 for the tame filtration comparison.

\bibitem[Yu01]{Yu2001}
Jiu-Kang Yu, \emph{Construction of tame supercuspidal representations},
J. Amer. Math. Soc. \textbf{14} (2001), no.~3, 579--622.
\end{thebibliography}
\end{document}